\documentclass[12pt]{article}

\usepackage{tikz, enumitem, caption, subcaption, float}
\usepackage{indentfirst, hyperref}
\usepackage{amsmath,amsfonts,amssymb,amsthm,epsfig,epstopdf,titling,url,array, mathdots}
\usepackage[margin=1in]{geometry}
\usepackage{mathtools}

\usepackage{nccmath}
\newtagform{roman}[]()

\usetikzlibrary{snakes}
\usetikzlibrary{decorations.markings}
\tikzset{unlabelled/.style={black, draw=black, circle, fill, scale=0.5}}
\tikzset{labelled/.style={black, draw=black, circle, scale=0.7, font=\Large}}

\tikzset{->-/.style={decoration={
			markings,
			mark=at position #1 with {\arrow{>}}},postaction={decorate}}}
\tikzset{dashedarc/.style={decoration={
			markings,
			mark=at position #1 with {\arrow{>}}},postaction={decorate}}}
		
\usepackage{algorithm, algorithmic}

\begin{document}

\newtheorem{theorem}{\hspace{5mm}Theorem}[section]
\newtheorem{corollary}[theorem]{\hspace{5mm}Corollary}
\newtheorem{lemma}[theorem]{\hspace{5mm}Lemma}
\newtheorem{proposition}[theorem]{\hspace{5mm}Proposition}
\newtheorem{definition}[theorem]{Definition\hspace{5mm}}

\newcommand{\lemmaautorefname}{Lemma}
\newcommand{\corollaryautorefname}{Corollary}
\newcommand{\propositionautorefname}{Proposition}

\newcommand{\pf}{{\bf Proof: }}

\newenvironment{AMS}{}{}
\newenvironment{keywords}{}{}

\title{Obstructions for local tournament orientation completions}

\author{Kevin Hsu\ and\ Jing Huang 
     \thanks{Department of Mathematics and Statistics,
      University of Victoria, Victoria, B.C., Canada V8W 2Y2; 
       huangj@uvic.ca; Research supported by NSERC}}
\date{}

\maketitle

\begin{abstract}
Orientation completion problems commonly generalize orientation and recognition
problems for graph and digraph classes, and parallel the study of representation
extension problems.

The orientation completion problem for a class of oriented graphs asks whether
a given partially oriented graph can be completed to an oriented graph in 
the class by orienting the unoriented edges of the partially oriented graph.
Orientation completion problems have been studied recently for several classes of
oriented graphs, yielding both polynomial time solutions as well as NP-completeness
results.

Local tournaments are a well-structured class of oriented graphs that generalize
tournaments and their underlying graphs are intimately related to proper 
circular-arc graphs. According to Skrien, a connected graph can be oriented as 
a local tournament if and only if it is a proper circular-arc graph. 
Proper interval graphs are precisely the graphs which can be oriented as acyclic 
local tournaments. It has been proved that the orientation completion problems for 
the classes of local tournaments and acyclic local tournaments are both polynomial 
time solvable.

In this paper we characterize the partially oriented graphs that can be completed 
to local tournaments by determining the complete list of obstructions. These are in 
a sense minimal partially oriented graphs that cannot be completed to 
local tournaments. 
The result may be viewed as an extension of the well-known forbidden subgraph 
characterization of proper circular-arc graphs obtained by Tucker.
\end{abstract}

\section{Introduction}

One of the fundamental problems in graph theory is to ask whether a given graph 
has an orientation that satisfies a prescribed property and to find such an 
orientation if it exists. A celebrated theorem of Robbins \cite{robbins} which 
answers a question of this type states that a graph has a strong orientation if 
and only if it is 2-edge-connected (i.e., has no bridge). This theorem not only
characterize the graphs which have strong orientations but also provides a 
certificate (i.e., a strong orientation) to show a graph is 2-edge-connected. 

Orientation completion problems (defined below) commonly generalize orientation
problems as well as recognition problems for certain classes of graphs and digraphs 
\cite{bhz,hhh,huang_locp}. They also parallel the study of representation 
extension problems 
\cite{cfk,cggkl,klavik,kkkw,kkorssv,kkos,kkosv,ks,zeman}. 

We consider graphs, digraphs, and partially oriented graphs in this paper. For
graphs we assume that they do not contain loops or multiple edges (i.e., they are
{\em simple}), and for digraphs we assume they do not contain loops or two arcs
joining the same pair of vertices (i.e., they are {\em oriented} graphs).

A {\em partially oriented graph} is a mixed graph $H$ obtained from some graph $G$
by orienting the edges in a subset of the edge set of $G$. The graph $G$ is called
the {\em underlying graph} of $H$. We denote $H$ by $(V,E\cup A)$ where $E$ is
the set of (non-oriented) edges and $A$ is the set of arcs in $H$.
We use $uv$ to denote an edge in $E$ with endvertices $u, v$ and use $(u,v)$ to
denote an arc in $A$ with {\em tail} $u$ and {\em head} $v$.
In either case we say that $u, v$ are {\em adjacent} in $H$.
The partially oriented graph $H$ is {\em connected} if its underlying graph $G$ is.

A class $\cal C$ of graphs is called {\em hereditary} if it is closed under taking 
induced subgraphs, that is, if $G \in \cal C$ and $G'$ is an induced subgraph of 
$G$ then $G' \in \cal C$. Similarly, a class of digraphs is {\em hereditary} if
it is closed under taking induced subdigraphs. We extend this concept to
partially oriented graphs.

Let $H = (V,E\cup A)$ and $H' = (V',E'\cup A')$ be partially oriented graphs.
We says that $H$ {\em contains} $H'$ (or $H'$ is {\em contained} in $H$) if 
$V' \subseteq V$ and for all $u, v \in V'$,
\begin{itemize}
\item $u$ and $v$ are adjacent in $H'$ if and only if they are adjacent in $H$;
\item if $(u,v) \in A'$ then $(u,v) \in A$;
\item if $uv \in E'$, then $uv \in E$, or $(u,v) \in A$, or $(v,u) \in A$.
\end{itemize}.
Equivalently, $H'$ is contained in $H$ if and only if it is obtained 
from $H$ by deleting some vertices, followed by replacing some arcs $(u,v)$ with 
edges $uv$.

We note that, in case when $H$ and $H'$ are both graphs or digraphs, $H$
contains $H'$ if and only if $H$ contains $H'$ as an induced subgraph
or as an induced subdigraph. We call a class $\cal C$ of partially oriented graphs 
{\em hereditary} if $H \in \cal C$ and $H'$ is contained in $H$ then 
$H' \in \cal C$.
 
Fix a hereditary class $\cal C$ of oriented graphs,
the {\em orientation completion problem} for $\cal C$ asks whether a given 
partially oriented graph $H = (V,E\cup A)$ can be completed to an oriented graph 
in $\cal C$ by orienting the edges in $E$. 
The hereditary property of $\cal C$ ensures that if a partially oriented graph $H$ 
can be completed to an oriented graph in $\cal C$ then every partially 
oriented graph that is contained in $H$ can also be completed to 
an oriented graph in $\cal C$. Therefore the partially oriented graphs which can be
completed to oriented graphs in $\cal C$ form a hereditary class. 

Orientation completion problems have been studied for several classes of
oriented graphs, including local tournaments, local transitive tournaments, and
acyclic local tournaments, cf. \cite{bhz,huang_locp}. 
A {\em local tournament} is an oriented 
graph in which the in-neighbourhood as well as the out-neighbourhood of each vertex
induces a tournament. When the the in-neighbourhood as well as 
the out-neighbourhood of each vertex induces a transitive tournament, 
the oriented graph is called a {\em local transitive tournament}. 
If a local tournament does not contain a directed cycle then it is called 
{\em acyclic}. These three classes of oriented graphs are nested; 
the class of local tournaments properly contains local transitive tournaments, 
which in turn as a class properly contains acyclic local tournaments.
It has been proved in \cite{bhz} that the orientation completion problem is
polynomial time solvable for local tournaments and for acyclic local tournaments,
but NP-complete for local transitive tournaments.

Any hereditary class of graphs or digraphs admits a characterization by forbidden
subgraphs or subdigraphs. The forbidden subgraphs or subdigraphs consists of 
minimal graphs or digraphs which do not belong to the class and they are sometimes 
called obstructions for the class. It turns out this is also the case for 
a hereditary class of partially oriented graphs and in particular for the class of 
partially oriented graphs which can be completed to local tournaments.   

We call a partially oriented graph $X = (V,E\cup A)$ an {\em obstruction for 
local tournament orientation completions} (or simply, an \textit{obstruction}) if
the following three properties hold:

\begin{enumerate}
\item $X$ cannot be completed to a local tournament;
\item For each $v \in V$, $X-v$ can be completed to a local tournament;
\item For each $(u,v) \in A$, the partially oriented graph obtained from
           $X$ by replacing $(u,v)$ with the edge $uv$ can be completed to
           a local tournament.
\end{enumerate}

Thus an obstruction $X$ is a partially oriented graph which cannot be completed to 
a local tournament and is minimal in the sense that if $X'$ is contained
in $X$ and $X' \neq X$ then $X'$ can be completed to a local tournament. 

Clearly, every obstruction must be connected. If an obstruction contains no arc 
then it is a forbidden subgraph for the class of graphs which can be oriented as 
local tournaments (or equivalently, proper circular-arc graphs, 
cf. Theorem~\ref{Tucker}). We shall prove that if an obstruction
contains arcs then it contains exactly two arcs. 

The {\em dual} of an obstruction $X$ is the one obtained from $X$ by reversing 
the arcs (if any) in $X$. 
If an obstruction $X$ does not contain arcs then the dual of $X$ is $X$ itself.
Clearly, the dual of an obstruction is again an obstruction. 
%Obstructions are present under critical containment in any partially oriented 
%graph that can not be completed to a local tournament.

If $H$ can be completed to a local tournament, then every partially oriented graph 
contained in $H$ can also be completed to a local tournament so $H$ does not contain
an obstruction. On the other hand, suppose that $H$ cannot be completed to a local 
tournament. By deleting vertices and replacing arcs with edges in $H$ as long as 
the resulting partially oriented graph still cannot be completed to a local 
tournament we obtain an obstruction that is contained in $H$. Therefore,
a partially oriented graph $H$ cannot be completed to a local tournament if and
only if it contains an obstruction.

As mentioned above, the orientation completion problems for the class of acyclic 
local tournaments and for the class of general local tournaments are both 
polynomial time solvable. This suggests the possibility of finding all obstructions
for the completions of either classes. Indeed, the obstructions for acyclic
local tournament orientation completions have been found in \cite{hhh}. 
In this paper we will find all obstructions for local tournament orientation 
completions and hence characterize partially oriented graphs by obstructions
which can be completed to local tournaments. Specifically, we prove the following:

\begin{theorem} \label{main}
A partially oriented graph can be completed to a local tournament if and only if
it does not contain any of obstructions $X$ listed below:
\begin{itemize}
\item $X$ has no arcs and $\overline{X}$ is a $C_{2k}$ with $k \geq 3$, a 
      $C_{2k+1}+K_1$ with $k \geq 1$, or a graph in Figure~\ref{tuckerlist};
\item $X$ has exactly two arcs, and $X$ or its dual is a graph in 
      Figures~\ref{dividing-cut-vertex-figure}--\ref{1non-dividing8};
\item $X$ has exactly two arcs and $\overline{U(X)}$ is a graph in 
      Figures~\ref{disconnectedfigure}--\ref{c5} (with arcs in $X$ or its dual
      being specified in the figures).
\end{itemize}
\end{theorem}

\section{Preliminary results} \label{Preliminary Results}

A {\em proper circular-arc} graph is the intersection graph of
a family of circular-arcs on a circle where no circular-arc contains another.
Skrien \cite{skrien} proved that a connected graph is a proper circular-arc graph
if and only if it can be oriented as a local tournament. It follows that if
a partially oriented graph $H$ can be completed to a local tournament then each
component of the underlying graph of $H$ is a proper circular-arc graph.
	
Tucker \cite{tucker} found all minimal graphs which are not proper circular-arc 
graphs.
	
\begin{theorem} \label{Tucker} \cite{tucker} 
A graph $G$ is a proper circular-arc graph if and only if $G$ does not contain 
$C_k + K_1$ ($k \geq 4$) or tent $+ K_1$ and $\overline{G}$ does not contain
$C_{2k}$ ($k \geq 3$), $C_{2k+1}+K_1$ ($k \geq 1$), or any of the graphs in 
Figure~\ref{tuckerlist} as an induced subgraph.
\qed
\end{theorem}

\begin{corollary} \label{noarc}
Let $X$ be an obstruction for local tournament orientation completions.
If $X$ has no arc then $\overline{X}$ is a $C_{2k}$ with $k \geq 3$, a
$C_{2k+1}+K_1$ with $k \geq 1$, or a graph in Figure~\ref{tuckerlist}.
\end{corollary}
\pf Since $X$ is an obstruction for local tournament orientation completions that
has no arc, it is a forbidden subgraph for proper circular arc graphs.
Since $X$ is connected, $X$ is not a $C_k + K_1$ ($k \geq 4$) or tent $+ K_1$.
Hence, by Theorem~\ref{Tucker}, $\overline{X}$ is a $C_{2k}$ with $k \geq 3$, a
$C_{2k+1}+K_1$ with $k \geq 1$, or a graph in Figure~\ref{tuckerlist}.
\qed

\begin{figure}[H]
	\centering
	\captionsetup[subfigure]{labelformat=empty}

	\begin{subfigure}{0.33\textwidth}
		\centering
		\begin{tikzpicture}
			\node[unlabelled] (a) at (1, 1.73)		{};
			\node[unlabelled] (b) at (0.5, 0.865)	{};
			\node[unlabelled] (c) at (1.5, 0.865)	{};
			\node[unlabelled] (d) at (0, 0)			{};
			\node[unlabelled] (e) at (1, 0)			{};
			\node[unlabelled] (f) at (2, 0)			{};
			\draw 	(a) edge (b)	(a) edge (c)	(b) edge (c)
			(b) edge (d)	(b) edge (e)	(c) edge (e)
			(c) edge (f)	(d) edge (e)	(e) edge (f);
		\end{tikzpicture}
\subcaption{tent}
	\end{subfigure}%
	\begin{subfigure}[b]{0.33\textwidth}
		\centering
		\begin{tikzpicture}
			\node[unlabelled] (a) at (1, 1.7)	{};
			\node[unlabelled] (b) at (0, 1)	{};
			\node[unlabelled] (c) at (1, 1)	{};
			\node[unlabelled] (d) at (2, 1)	{};
			\node[unlabelled] (e) at (0, 0)	{};
			\node[unlabelled] (f) at (1, 0)	{};
			\node[unlabelled] (g) at (2, 0)	{};
			\draw	(a) edge (c)	(b) edge (c)	(b) edge (e)
			(c) edge (d)	(c) edge (f)	(d) edge (g)
			(e) edge (f)	(f) edge (g);
		\end{tikzpicture}
%\subcaption{$T_4$}
	\end{subfigure}%
	\begin{subfigure}[b]{0.33\textwidth}
		\centering
		\begin{tikzpicture}
			\node[unlabelled] (a) at (0, 1)	{};
			\node[unlabelled] (b) at (1, 1)	{};
			\node[unlabelled] (c) at (2, 1)	{};
			\node[unlabelled] (d) at (0, 0)	{};
			\node[unlabelled] (e) at (1, 0)	{};
			\node[unlabelled] (f) at (2, 0)	{};
			\node[unlabelled] (g) at (1.5, 0.5)	{};
			\draw	(a) edge (b)	(b) edge (c)	(a) edge (d)
			(b) edge (e)	(c) edge (f)	(d) edge (e)
			(e) edge (f)	(b) edge (g)	(g) edge (e);
		\end{tikzpicture}
%\subcaption{$T_5$}
	\end{subfigure}
	\par\bigskip
	\par\bigskip

	\begin{subfigure}{0.5\textwidth}
		\centering
		\begin{tikzpicture}
			\node[unlabelled] (a) at (0, 0)	{};
			\node[unlabelled] (b) at (1, 0)	{};
			\node[unlabelled] (c) at (2, 0)	{};
			\node[unlabelled] (d) at (3, 0)	{};
			\node[unlabelled] (e) at (4, 0)	{};
			\node[unlabelled] (f) at (2, 1)	{};
			\node[unlabelled] (g) at (2, 2)	{};
			\draw	(a) edge (b)	(b) edge (c)
			(c) edge (d)
			(d) edge (e)	(b) edge (f)	(d) edge (f)
			(f) edge (g);
		\end{tikzpicture}
%\subcaption{$T_6$}
	\end{subfigure}%
	\begin{subfigure}{0.5\textwidth}
		\centering
		\begin{tikzpicture}
			\node[unlabelled] (a) at (0, 0)	{};
			\node[unlabelled] (b) at (1, 0)	{};
			\node[unlabelled] (c) at (2, 0)	{};
			\node[unlabelled] (d) at (3, 0)	{};
			\node[unlabelled] (e) at (4, 0)	{};
			\node[unlabelled] (f) at (2, 1)	{};
			\node[unlabelled] (g) at (2, 2)	{};
			\draw	(a) edge (b)	(b) edge (c)	(c) edge (d)
			(d) edge (e)	(c) edge (f)	(f) edge (g);
		\end{tikzpicture}
%\subcaption{$T_7$}
	\end{subfigure}
	\caption{}
	\label{tuckerlist}
\end{figure}

In view of Corollary~\ref{noarc} we only need to find obstructions that contain 
arcs. By definition the underlying graph of any obstruction that contains arcs 
is a proper circular-arc graph and hence local tournament orientable. 

Let $G = (V,E)$ be a graph and $Z(G) = \{(u,v): uv \in E\}$ be the set of all 
ordered pairs $(u,v)$ such that $uv \in E$. Note that each edge $uv \in E$ 
gives rise to two ordered pairs $(u,v), (v,u)$ in $Z(G)$. Suppose that
$(u, v)$ and $(x, y)$ are two ordered pairs of $Z(G)$. We say $(u,v)$ 
\textit{forces} $(x,y)$ and write $(u,v) \Gamma (x,y)$ if one of the following 
conditions is satisfied:
	
\begin{itemize}
\item $u = x$ and $v = y$;
\item $u = y$, $v \neq x$, and $vx \notin E$;
\item $v = x$, $u \neq y$, and $uy \notin E$.
\end{itemize}
	
\noindent We say that $(u,v)$ \textit{implies} $(x,y)$ and write 
$(u,v) \Gamma^* (x,y)$ if there exists a sequence of pairs 
$(u_1, v_1), (u_2, v_2), \dots, (u_k, v_k) \in Z(G)$ such that 
$$(u, v) = (u_1, v_1) \Gamma (u_2, v_2) \Gamma \dots \Gamma (u_k, v_k) = (x, y).$$
\noindent We will call such a sequence a $\Gamma$-\textit{sequence} from 
$(u, v)$ to $(x, y)$. It is easy to verify that $\Gamma^*$ is an equivalence 
relation on $Z(G)$. 

We say a path $P$ {\em avoids} a vertex $u$ if $P$ does not contain $u$ or any
neighbour of $u$.

\begin{proposition}\label{complement gamma}
Let $G$ be a graph and $u, v, w$ be vertices Suppose that $P$ is a path of length
$k$ connecting $v, w$ that avoids $u$ in $\overline{G}$. If $k$ is even, then 
$(u,v)\Gamma^*(u,w)$. Otherwise, $(u,v)\Gamma^*(w,u)$.
\end{proposition}
\pf Denote $P: v_0v_1 \dots v_k$ where $v_0 = v$ and $v_k = w$. Since $P$ avoids 
$u$ in $\overline{G}$, $(u,p_i) \Gamma (p_{i+1},u)$ for each $0 \leq i \leq k-1$. 
If $k$ is even, then 
$$(u,v) = (u,p_0) \Gamma (p_1,u)\Gamma (u,p_2)\Gamma\dots\Gamma(u,p_k) = (u,w).$$ 
Otherwise, 
$$(u,v) = (u,p_0)\Gamma(p_1,u)\Gamma(u,p_2)\Gamma\dots\Gamma(p_k,u) = (w,u).$$
\qed

\begin{proposition}[] \cite{huang} \label{Gamma-class}
Let $G$ be a graph and $D = (V,A)$ be a local tournament orientation of $G$.
Suppose that $(u, v) \Gamma^* (x, y)$ for some $(u, v), (x, y) \in Z(G)$. Then
$(u,v) \in A$ if and only if $(x,y) \in A$.
\qed
\end{proposition}

Regardless whether or not $G$ is local tournament orientable, the relation 
$\Gamma^*$ on $Z(G)$ induces a partition of the edge set of $G$ into 
\textit{implication classes} as follows: two edges $uv, xy$ of $G$ 
are in the same implication class if and only if $(u,v) \Gamma^* (x,y)$ or 
$(u,v) \Gamma^* (y,x)$. An implication class is called {\em trivial} if it has
only one edge and {\em non-trivial} otherwise. An edge $uv$ of $G$ is called 
{\em balanced} if $N[u] = N[v]$ and {\em unbalanced} otherwise. Clearly, 
any balanced edge forms a trivial implication class and the unique edge in any 
trivial implication class is balanced.

The following theorem characterizes the implication classes of a local tournament 
orientable graph and describes all possible local tournament orientations of such
a graph.
	
\begin{theorem}[]\label{structure} \cite{huang} 
Let $G = (V,E)$ be a connected graph and let $H_1, H_2, \dots, H_k$ be the 
components of $\overline{G}$. Suppose that $G$ is local tournament orientable and 
$F$ is an implication class of $G$. Then $F$ is one of the following types:
\begin{itemize}
\item $F$ is trivial;
\item $F$ consists of all unbalanced edges of $G$ within $H_i$ for some $i$;
\item $F$ consists of all edges of $G$ between $H_i$ and $H_j$ for some $i \neq j$.
\end{itemize}
Moreover, suppose that $F_1, F_2, \dots, F_\ell$ are the implication classes of 
$G$. For each $1 \leq i \leq \ell$, let $A_i$ be the equivalence class of 
$\Gamma^*$ containing $(u,v)$ for some $uv \in F_i$ and let 
$A = \cup_{i=1}^{\ell} A_i$. Then $D = (V,A)$ is a local tournament orientation
of $G$.
\qed
\end{theorem}

Let $H = (V,E\cup A)$ be a partially oriented graph and $(a,b), (c,d)$ be arcs of
$H$. We say that the two arcs $(a,b), (c,d)$ are {\em opposing} in $H$ if 
$(a,b) \Gamma^* (d,c)$. For convenience we also call an arc of $H$ {\em balanced}
if the corresponding edge is balanced. Clearly, if $(a,b), (c,d)$ are opposing 
then neither of them is balanced.

\begin{proposition} \label{oppo} 
Suppose that $H$ is a partially oriented graph whose underlying graph $U(H)$ is 
local tournament orientable. Then $H$ can be completed to a local tournament
if and only if it does not contain opposing arcs.
\end{proposition}
\pf If $H$ contains opposing arcs, then by Proposition~\ref{Gamma-class} it cannot 
be completed to a local tournament. On the other hand, suppose that $H$ does not 
contain opposing arcs. Let $F_1, F_2, \dots, F_{\ell}$ be the implication classes 
of $U(H)$. For each $1 \leq i \leq \ell$, if no edge in $F_i$ is oriented then let 
$A_i$ be an equivalence class of $\Gamma^*$ containing $(u,v)$ for some 
$uv \in F_i$; otherwise let $A_i$ be the equivalence class of $\Gamma^*$ containing
$(u,v)$ where $uv \in F_i$ and $(u,v)$ is an arc. With $A = \cup_{i=1}^{\ell} A_i$,
Theorem~\ref{structure} ensures that $D = (V,A)$ is a local tournament completion
of $H$.
\qed

The next theorem is fundamental in determining whether a partially oriented 
graph whose underlying graph is local tournament orientable is an obstruction.
	
\begin{theorem} \label{different implication classes theorem}
Let $X$ be a partially oriented graph whose underlying graph $U(X)$ is local
tournament orientable. Then $X$ is an obstruction if and only if $X$ contains 
exactly two arcs $(a,b), (c,d)$ which are opposing and, for every vertex 
$v \in V(X) \setminus \{a,b,c,d\}$, the arcs $(a,b), (c,d)$ are not opposing in
$X-v$ (that is, the edges $ab, cd$ belong to different implication classes in
$U(X-v)$). Moreover, any $\Gamma$-sequence connecting $(a,b)$ and $(d,c)$ must
include all vertices of $X$.
\end{theorem}
\pf For sufficiency, suppose that $(a,b), (c,d)$ are the only arcs and they are 
opposing in $X$ and that, for every vertex $v \in V(X) \setminus \{a,b,c,d\}$, 
the arcs $(a,b), (c,d)$ are not opposing in $X-v$. Since $X$ contains opposing 
arcs, it cannot be completed to a local tournament by Proposition~\ref{oppo}.
Let $v$ be a vertex in $X$. Since $U(X)$ is local tournament orientable, 
$U(X-v)$ is also local tournament orientable. If $v \in \{a,b,c,d\}$, then $X-v$
contains at most one arc and hence no opposing arcs. If $v \notin \{a,b,c,d\}$, 
then the only two arcs in $X-v$ are not opposing by assumption. Hence $X-v$ can 
be completed to a local tournament by Proposition~\ref{oppo}. Therefore $X$ is 
an obstruction. 

Conversely, suppose that $X$ is an obstruction. By Proposition~\ref{oppo} $X$ 
must contain opposing arcs. Let $(a,b), (c,d)$ be opposing arcs in $X$. If $X$ 
contains an arc $(x,y)$ that is distinct from $(a,b), (c,d)$, then replacing
the arc $(x,y)$ by the edge $xy$ gives a partially orientable graph in which 
$(a,b), (c,d)$ are still opposing and hence cannot be completed to a local 
tournament. This contradicts the assumption that $X$ is an obstruction. So 
$(a,b), (c,d)$ are the only arcs in $X$. Since $X$ is an obstruction, for every
$v \in V(X)$, $X-v$ can be completed to a local tournament and hence
by Proposition~\ref{oppo} contains no opposing arcs. This implies in particular 
that if $v \in V(X) \setminus \{a,b,c,d\}$, the arcs $(a,b), (c,d)$ are not 
opposing in $X-v$.

The second part of the theorem follows from the fact that deleting any vertex 
results in a graph that contains no $\Gamma$-sequence connecting 
$(a,b)$ and $(d,c)$.
\qed

Let $v$ be a vertex and $(x,y)$ be an arc in a partially oriented graph $H$ where 
$v \notin \{x,y\}$. We call $v$ the {\em $(x,y)$-balancing} vertex if $v$ is 
the only vertex adjacent to exactly one of $x,y$; when the arc $(x,y)$ does not 
need to be specified, we simply call $v$ an {\em arc-balancing} vertex.
Each obstruction has at most two arc-balancing vertices as it contains at most 
two arcs.

A vertex of a graph $G$ is called a {\em cut-vertex} of $G$ if $G-v$ has more 
components than $G$. For a partially oriented graph $H$, a cut-vertex of $U(H)$ is 
also called a {\em cut-vertex} of $H$. 

\begin{proposition}\label{vertex classification}
Let $X$ be an obstruction with opposing arcs $(a,b),(c,d)$ and let 
$v \notin \{a,b,c,d\}$. Then $v$ is an arc-balancing vertex, or a cut-vertex of 
$U(X)$, or a cut-vertex of $\overline{U(X)}$.
\end{proposition}
\pf Assume that $v$ is not a cut-vertex of $U(X)$ or of $\overline{U(X)}$ as 
otherwise we are done. We show that $v$ must be an arc-balancing vertex. 
Since $ab, cd$ are in the same implication class of $U(X)$, 
by Theorem \ref{structure} 
$ab, cd$ are unbalanced edges either contained in a component or between two 
components of $\overline{U(X)}$. Since $v$ is not a cut-vertex of 
$\overline{U(X)}$, each component of $\overline{U(X-v)}$ is a component of 
$\overline{U(X)}$ except possibly missing $v$. It follows that $ab, cd$ are 
contained in some component or between two components of $\overline{U(X-v)}$. 
Since $v$ is not a cut-vertex of $U(X)$, $U(X-v)$ is connected.
If $ab, cd$ are both unbalanced edges in $U(X-v)$, then they remain in the same 
implication class of $U(X-v)$ and hence $(a,b),(c,d)$ are still opposing in
$X-v$, which contradicts the assumption that $X$ is an obstruction.
So one of $ab, cd$ is balanced in $U(X-v)$, which means that $v$ is
$(a,b)$-balancing or $(c,d)$-balancing.
\qed

An {\em arc-balancing triple} in a partially oriented graph $H$ is a set of three 
vertices in which one balances an arc between the other two. 

\begin{corollary}\label{at most 6 non-cut-vertices}
Let $X$ be an obstruction with opposing arcs $(a,b),(c,d)$. Suppose that $U(X)$ has 
no cut-vertices. Then $\overline{U(X)}$ contains at most six non-cut-vertices. 
In the case when $\overline{U(X)}$ has six non-cut-vertices, the six 
non-cut-vertices form two disjoint arc-balancing triples.
\end{corollary}   
\pf Let $v$ be a non-cut-vertex of $\overline{U(X)}$. By assumption $v$ is not 
a cut-vertex of $U(X)$ and thus, by Proposition~\ref{vertex classification}, it
is either in $\{a,b,c,d\}$ or an arc-balancing vertex. There are at most two
arc-balancing vertices so $\overline{U(X)}$ contains at most six non-cut-vertices.
When $\overline{U(X)}$ has six non-cut-vertices, among the six non-cut-vertices
two are arc-balancing vertices and the other four are incident with arcs. Hence
the six non-cut-vertices form two disjoint arc-balancing triples. 
\qed                    

A {\em proper interval graph} is the intersection graph of a family of intervals
in a line where no interval contains another. Proper interval graphs form a
prominent subclass of proper circular-arc graphs and play an important role in
the orientation completion problem for local tournaments. It is proved in
\cite{hh} that a graph is a proper interval graph if and only if it can be
oriented as an acyclic local tournament.

A {\em straight enumeration} of a graph $G$ is a vertex ordering $\prec$ such that
for all $u \prec v \prec w$, if $uw$ is an edge of $G$, then both $uv$ and $vw$ are
edges. This property is referred to as the {\em umbrella property} of the vertex 
ordering.  

\begin{proposition}\label{umbrella property} \cite{huang}
A graph is a proper interval graph if and only if it has a straight enumeration.
\qed
\end{proposition}

\begin{proposition}\label{pig_gamma_sequence}
Let $G$ be a connected proper interval graph and let $\prec$ be a straight
enumeration of $G$. Suppose that $(u,v) \Gamma^* (x,y)$. Then $u \prec v$ if and
only if $x \prec y$.
\end{proposition}
\pf It suffices to show that if $u \prec v$ and $(u,v) \Gamma (x,y)$ then
$x \prec y$. So assume that $(u,v) \Gamma (x,y)$. Then one of the following holds:
\begin{itemize}
\item $u = x$ and $v = y$;
\item $u = y$, $v \neq x$, and $vx \notin E(G)$;
\item $v = x$, $u \neq y$, and $uy \notin E(G)$.
\end{itemize}
Clearly, $x \prec y$ when $u = x$ and $v = y$.
Suppose that $u = y$, $v \neq x$, and $vx \notin E(G)$. If $u \prec x \prec v$,
then it violates the umbrella property because $uv \in E(G)$ but $xv \notin E(G)$.
If $u \prec v \prec x$, then it again violates the umbrella property because
$ux \in E(G)$ but $vx \notin E(G)$. Hence we must have $x \prec u = y$.
The proof for the case when $v = x$, $u \neq y$, and $uy \notin E(G)$ is similar.
\qed

Let $H$ be a partially oriented graph whose underlying graph $U(H)$ is a proper 
interval graph. Suppose that $\prec$ is a straight enumeration of $U(H)$. We call 
an arc $(u,v)$ of $H$ \textit{positive} (with respect to $\prec$) if $u \prec v$ 
and \textit{negative} otherwise. If $H$ does not contain negative arcs, then $H$ 
can be completed to an acyclic local tournament by replacing all edges of $H$ with 
positive arcs. Similarly, if $H$ does not contain positive arcs then it can also
be completed to an acyclic local tournament. It follows that if $X$ is an
obstruction such that $U(X)$ is a proper interval graph, then the two arcs in $X$
must be {\em opposite} (i.e., one is positive and the other is negative).

A vertex in a graph is {\em universal} if it is adjacent to every other vertex.

\begin{theorem} \cite{huang} \label{pigstructure}
Suppose that $G = (V,E)$ is a connected proper interval graph that is not a
complete graph. Then $\overline{G}$ has a unique non-trivial component $H$. 
If $F$ is an implication class of $G$, then $F$ is one of the following types: 
\begin{itemize}
\item $F$ is trivial;
\item $F$ consists of all unbalanced edges within $H$;
\item $F$ consists of all edges of $G$ between $H$ and a universal vertex of $G$.
\end{itemize}
In particular, if $G$ contains no universal vertex, then $G$ has a unique
non-trivial implication class.
\qed
\end{theorem}

\begin{proposition}\label{universal vertex maker}
Let $G$ be a connected proper interval graph and let $v_1, v_2,\dots, v_n$ be a 
straight enumeration of $G$. Suppose that $v_{\alpha}$ is a cut-vertex of 
$\overline{G}$. Then $\alpha \in \{1, n\}$ and $G-v_{\alpha}$ contains a vertex 
that is adjacent to every vertex except $v_{\alpha}$ in $G$.
\end{proposition}
\pf Since $\overline{G}$ has a cut-vertex, $G$ is not a complete graph and by 
Theorem~\ref{pigstructure}, $\overline{G}$ has a unique non-trivial component 
$H$. Thus the cut-vertex $v_{\alpha}$ of $\overline{G}$ is in fact a cut-vertex of 
$H$. Again by Theorem~\ref{pigstructure}, $H-v_{\alpha}$ has at most one 
non-trivial component. 
Hence $H$ contains a vertex $v_{\beta}$ that is only adjacent to 
$v_{\alpha}$ in $\overline{G}$, that is, in $G$ it is adjacent to every vertex 
except $v_{\alpha}$. If $\alpha < \beta$, then $\alpha = 1$ as otherwise 
we have $1 < \alpha < \beta$ and $v_{\beta}$ is adjacent to $v_1$ but not to 
$v_{\alpha}$, a contradiction to the umbrella property of the straight enumeration.
Similarly, if $\beta < \alpha$, then $\alpha = n$ as otherwise 
$\beta < \alpha < n$ and $v_{\beta}$ is adjacent to $v_n$ but not to $v_{\alpha}$,
also a contradiction to the umbrella property of the straight enumeration.
Therefore, $\alpha \in \{1, n\}$.
\qed

\section{Obstructions with cut-vertices}

Our goal is to find all obstructions for local tournament orientation completions
that contains arcs. By Theorem \ref{different implication classes theorem}
each of them contains exactly two arcs which are opposing and its underlying 
graph is a connected proper circular arc graph (i.e., local tournament orientable).
In this section, we examine such obstructions that contain cut-vertices. 
Since they contain cut-vertices, their underlying graphs are necessarily proper 
interval graphs and thus have straight enumerations according to
Proposition~\ref{umbrella property}.
	
Let $X$ be an obstruction that contains two arcs and let $\prec$ be a straight 
enumeration of $U(X)$. Suppose that $v$ is a cut-vertex of $X$. Then $v$ is neither 
the first nor the last vertex in $\prec$ and moreover, for all $u, w$ with
$u \prec v \prec w$, $uw$ is not an edge in $U(X)$. We call $v$ a {\em dividing} 
cut-vertex if one of the two arcs in $X$ is incident with a vertex preceding $v$ 
and the other is incident with a vertex succeeding $v$ in $\prec$; if $v$ is not 
dividing then it is called {\em non-dividing}.

\subsection{Dividing cut-vertices}

\begin{figure}[H]
	\captionsetup[subfigure]{labelformat=empty}
	\centering
	\useshortskip
	\usetagform{roman}
	\begin{equation}
	\begin{subfigure}{.5\textwidth}
	\centering
	\begin{tikzpicture}
	\node[unlabelled]	(1) at (0,0)	{};
	\node[unlabelled]	(2) at (1,0)	{};
	\node				(4) at (2,0)	{$\dots$};
	\node[unlabelled]	(6) at (3,0)	{};
	\node[unlabelled]	(7) at (4,0)	{};
	\draw	(2)--(4)	(4)--(6); 
	\draw[->-=.5]	(1) to (2);
	\draw[->-=.5]	(7) to (6);
	\end{tikzpicture}
	%\subcaption{$n \geq 3$}
	\end{subfigure}
	\end{equation}
	\begin{equation}
	\begin{subfigure}{.5\textwidth}
	\centering
	\begin{tikzpicture}
	\node[unlabelled]	(0) at (-1,0)	{};
	\node[unlabelled]	(1) at (0,0)	{};
	\node[unlabelled]	(2) at (1,0)	{};
	\node[unlabelled]	(3) at (2,0)	{};
	\node				(4) at (3,0)	{$\dots$};
	\node[unlabelled]	(6) at (4,0)	{};
	\node[unlabelled]	(7) at (5,0)	{};
	\draw	(0)--(1)	(2)--(3)	(3)--(4)	(4)--(6);
	\draw	(1) edge[bend left=45]	(3);
	\draw[->-=.5]	(1) to (2);
	\draw[->-=.5]	(7) to (6);
	\end{tikzpicture}
	%\subcaption{$n \geq 5$}
	\end{subfigure}
	\end{equation}
	\begin{equation}
	\begin{subfigure}{.5\textwidth}
	\centering
	\begin{tikzpicture}
	\node[unlabelled]	(0) at (-1,0)	{};
	\node[unlabelled]	(1) at (0,0)	{};
	\node[unlabelled]	(2) at (1,0)	{};
	\node[unlabelled]	(3) at (2,0)	{};
	\node				(4) at (3,0)	{$\dots$};
	\node[unlabelled]	(5) at (4,0)	{};
	\node[unlabelled]	(6) at (5,0)	{};
	\node[unlabelled]	(7) at (6,0)	{};
	\node[unlabelled]	(8) at (7,0)	{};
	\draw	(0)--(1)	(2)--(3)	(3)--(4)	(4)--(5)
	(5)--(6)	(7)--(8);
	\draw	(1) edge[bend left=45]	(3)
	(5) edge[bend left=45]	(7);
	\draw[->-=.5]	(1) to (2);
	\draw[->-=.5]	(7) to (6);
	\end{tikzpicture}
	%\subcaption{$n \geq 7$}
	\end{subfigure}%
	\end{equation}
	\caption{Obstructions with dividing cut-vertices.}
	\label{dividing-cut-vertex-figure}
\end{figure}
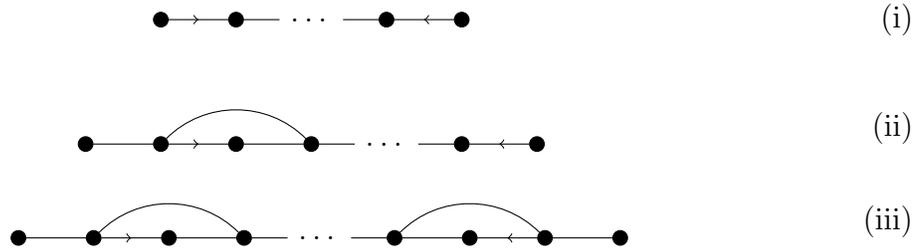

In this subsection, we focus on the obstructions that contain dividing cut-vertices.
We will show that they consist of the three infinite classes in 
Figure \ref{dividing-cut-vertex-figure} and their duals. In each of these graphs, 
the dots in the middle represent a path of length $\geq 0$; when the length of 
the path is 0 the two vertices beside the dots are the same vertex. 
	
\begin{lemma}\label{dividing cut-vertex lemma}
Let $X$ be an obstruction that contains a dividing cut-vertex and let 
$\prec:\ v_1, v_2, \dots, v_n$ be a straight enumeration of $U(X)$. Suppose that
$v_c$ is the first dividing cut-vertex in $\prec$. Then, either $c = 2$ and
$v_1,v_2$ are the endvertices of an arc, or $c =4$ and $v_2,v_3$ are the end 
vertices of an arc. In the case when $c = 4$, $v_1, v_2, v_3, v_4$ 
induce in $U(X)$ the following graph:
		
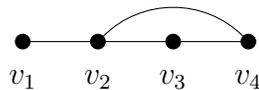
\begin{figure}[H]
			\centering
			\begin{tikzpicture}
			\node[unlabelled]	(1) at (0,0)	{};
			\node				(L1) at (0,-0.5){$v_1$};
			\node[unlabelled]	(2) at (1,0)	{};
			\node				(L2) at (1,-0.5){$v_2$};
			\node[unlabelled]	(3) at (2,0)	{};
			\node				(L3) at (2,-0.5){$v_3$};
			\node[unlabelled]	(4) at (3,0)	{};
			\node				(L4) at (3,-0.5){$v_4$};
			\draw	(1)--(2)	(3)--(4)	(2)--(3)
			(2) edge[bend left = 45] (4);
			\end{tikzpicture}
		\end{figure}
\end{lemma}

\pf By considering the dual of $X$ if necessary we may assume that $(v_j,v_k)$
and $(v_s,v_t)$ are the two arcs in $X$ where $j < k \leq c \leq t < s$. 
By Theorem \ref{different implication classes theorem}, there is a
$\Gamma$-sequence from $(v_t,v_s)$ to $(v_j,v_k)$ that includes all vertices
of $X$. Let
\[(v_t,v_s)=(u_1,w_1)\Gamma (u_2,w_2)\Gamma \dots \Gamma(u_q,w_q)=(v_j,v_k)\]
be a shortest such a sequence. Since $v_t \prec v_s$, we have
$u_i \prec w_i$ for each $i$ by Proposition \ref{pig_gamma_sequence}.
Let $\ell$ be the smallest subscript such that
$u_{\ell+1} \prec w_{\ell+1} = v_c = u_{\ell} \prec w_{\ell}$. 
Such $\ell$ exists because $v_c$ is a cut-vertex dividing $(v_j,v_k)$ and 
$(v_s,v_t)$. We distinguish two cases depending on whether or not $k=c$. 
Suppose first $k=c$.  
Note that $(u_{\ell},w_{\ell}) \Gamma (v_j,v_k)$. Thus the choice of 
the $\Gamma$-sequence implies 
$(u_{\ell+1},w_{\ell+1}) = (u_q,w_q) = (v_j,v_k)$. Since the $\Gamma$-sequence 
includes all vertices of $X$, $v_j$ is the only vertex preceding $v_c$ in $\prec$,
that is, $c = 2$ (and $(v_1,v_2)$ is an arc in $X$). 

Suppose now that $k < c$. Thus $j < k < c$. We claim that $v_j, v_k, v_c$ are 
consecutive vertices in $\prec$ (i.e., $j+1 = k = c-1$). Suppose that 
$k > j+1$. Since $v_j, v_k$ are adjacent, $v_{j+1}$ cannot be a cut-vertex of 
$U(X)$. Since $v_{j+1}$ is not the first or the last vertex in $\prec$, 
Proposition \ref{universal vertex maker} ensures that $v_{j+1}$ cannot be 
a cut-vertex of 
$\overline{U(X)}$. By Proposition \ref{vertex classification}, $v_{j+1}$ is an 
arc-balancing vertex. Clearly, $v_{j+1}$ is not $(v_j,v_k)$-balancing.
So it must be $(v_s,v_t)$-balancing. Since $j+1 < c$ and $v_c$ is a cut-vertex, 
$v_{j+1}$ has no neighbours succeeding $v_c$. It follows that $v_t = v_c$.
Since $v_{j+1}v_c$ is an edge and $j+1 < k < c$, $v_kv_c$ is an edge by
the umbrella property. Again, since $v_c$ is a cut-vertex, $v_k$ cannot be 
adjacent to $v_s$. This contradicts the fact that $v_{j+1}$ is arc-balancing 
for the arc between $v_s, v_t$. Hence $j+1 = k$, i.e., $v_j$ and $v_k$ are 
consecutive vertices in $\prec$.

Suppose $c > k+1$. Neither of $v_k, v_{k+1}$ can be a cut-vertex of $U(X)$ as 
otherwise it would be a dividing cut-vertex preceding $v_c$, a contradiction to 
the choice of $v_c$. Since $v_{k+1}$ is not the first or the last vertex in 
$\prec$, it is not a cut-vertex of $\overline{U(X)}$ according to 
Proposition \ref{universal vertex maker}. 
By Proposition \ref{vertex classification}, 
$v_{k+1}$ is an arc-balancing vertex. Since $v_k$ is not a cut-vertex of $U(X)$,
$v_{k-1} = v_j$ is adjacent to $v_{k+1}$. So $v_{k+1}$ is adjacent to both 
$v_j, v_k$ and hence not arc-balancing for the the arc between them. So 
$v_{k+1}$ is arc-balancing for the arc between $v_s, v_t$. Similarly as above
we have $v_c = v_t$ and $v_{k+1}$ is adjacent to $v_t$ but not to $v_s$. 
If $c > k+2$, then $v_{k+2}$ is adjacent to $v_c$ by the umbrella property
and the fact $v_{k+1}$ is adjacent to $v_c$. Thus $v_{k+2}$ is adjacent to 
$v_c = v_t$ but not to $v_s$, a contradiction to that $v_{k+1}$ is arc-balancing 
to the arc between $v_s, v_t$. If $c = k+2$, since $v_{k+1}$ is not
a cut-vertex of $U(X)$, $v_k$ is adjacent to $v_{k+2} = v_c$. Thus $v_k$ is 
adjacent to $v_t = v_c$ but not to $v_s$, a contradiction again to the fact
that $v_{k+1}$ is arc-balancing to the arc between $v_s, v_t$. Hence
$c = k+1$, i.e., $v_k, v_c$ are consecutive vertices in $\prec$.
Therefore $v_j, v_k, v_c$ are consecutive in $\prec$. 

Since $v_c$ is the first dividing cut-vertex in $\prec$, $v_k$ cannot be 
a cut-vertex of $U(X)$ and hence $v_j, v_c$ are adjacent in $X$.
We claim that there exists a vertex preceding $v_j$ in $\prec$ which is adjacent 
to $v_j$ but not to $v_k$. First, observe that if no vertex is adjacent to exactly one of $v_j, v_k$, then $v_j$ and $v_k$ would share the same closed neighbourhood.
In this case, the arc between $v_j$ and $v_k$ would be balanced, a contradiction. 
Hence, there is at least one vertex adjacent to exactly one of $v_j, v_k$. 
Clearly, such a vertex must precede $v_j$ in $\prec$ and hence is adjacent to 
$v_j$ but not to $v_k$. Assume that $v_p$ is such a vertex closest to $v_j$.

We show that $v_p$ and $v_j$ are consecutive in $\prec$, that is,
$p = j-1$. If $p < j-1$, then $v_{j-1}$ cannot be a cut-vertex of $U(X)$ 
because $v_p$ is adjacent to $v_j$. On the other hand, 
by Proposition \ref{universal vertex maker}, $v_{j-1}$ is not a cut-vertex of 
$\overline{U(X)}$. It follows from Theorem \ref{vertex classification} that 
$v_{j-1}$ is an arc-balancing vertex. 
The choice of $v_p$ implies that $v_{j-1}$ is adjacent to both $v_j, v_k$ so
it does not balance the arc between $v_j$ and $v_k$. Hence, $v_{j-1}$ is 
an arc-balancing vertex for the arc between $v_s$ and $v_t$. 
By definition it is the unique vertex adjacent to exactly one of $v_s$ and $v_t$.
This also implies $v_c = v_t$. But then $v_k$ is also a vertex adjacent to 
$v_t$ but not to $v_s$, a contradiction. Hence $p = j-1$. 

Since $(u_{\ell+1},w_{\ell+1}) \Gamma (u_{\ell},w_{\ell})$ and 
$u_{\ell+1} \prec w_{\ell+1} = v_c = u_{\ell} \prec w_{\ell}$ (i.e.,
$u_{\ell+1}$ is a vertex preceding and adjacent to $v_c$ but not adjacent to
$w_{\ell}$), $u_{\ell+1}$ can only be $v_{c-1}$ or $v_{c-2}$, 
Since $\Gamma$-sequence is chosen to be the shortest from $(v_t,v_s)$ to
$(v_j,v_k)$, we must have $u_{\ell+1} = v_{c-2}$. 
It follows that 
\[(v_t,v_s)=(u_1,w_1)\Gamma (u_2,w_2)\Gamma \dots (u_\ell,w_\ell)
\Gamma (v_{c-2},v_c) \Gamma (v_{c-3},v_{c-2}) \Gamma (v_{c-2},v_{c-1}) =(v_j,v_k)\]
is a shortest $\Gamma$-sequence. 
The $\Gamma$-sequence must contain all vertices of $X$, which means 
$v_{c-3},v_{c-2},v_{c-1}$ are all the vertices preceding $v_c$.
Therefore $c = 4$ and $v_1, v_2, v_3, v_4$ induce in $U(X)$ the graph in 
the statement. 
\qed

We can now apply Lemma \ref{dividing cut-vertex lemma} to prove the following:

\begin{theorem} \label{div-cut}
Let $X$ be an obstruction that contains a dividing cut-vertex with respect to
a straight enumeration. Then $X$ or its dual belongs to one of the three infinite 
classes in Figure \ref{dividing-cut-vertex-figure}.
\end{theorem}

\pf Let $\prec:\ v_1, v_2, \dots, v_n$ be a straight enumeration of $U(X)$ and
let $v_c$ and $v_d$ be the first and last dividing cut-vertices respectively 
with respect to $\prec$. By considering the dual of $X$ if necessary we assume 
that $(v_j,v_k)$ and $(v_s,v_t)$ are the arcs in $X$ where 
$j < k \leq c \leq d \leq t < s$.

Suppose $c = 2$ and $d = n-1$. Since $v_c = v_2$ is a cut-vertex, $v_2$ is
the only neighbour of $v_1$. Similarly, $v_{n-1}$ is the only neighbour of $v_n$.
If $v_p$ is adjacent to $v_q$ for some $2 \leq p < q-1 \leq n-1$, then 
it is easy to see that the partially oriented graph obtained from $X$ 
by deleting $v_{p+1}, \dots, v_{q-1}$ cannot be completed to local tournament 
orientation, a contradiction to the assumption $X$ is an obstruction.
Hence $X$ belongs to Figure \ref{dividing-cut-vertex-figure}(i). 

Suppose that $c \neq 2$. Then $c = 4$ by Lemma \ref{dividing cut-vertex lemma}. 
If $d = n-1$, then a similar argument as above shows that $X$ belongs to 
Figure \ref{dividing-cut-vertex-figure}(ii). On the other hand if 
$d \neq n-1$, then again by Lemma \ref{dividing cut-vertex lemma} we must have
$d = n-3$. In this case $X$ belongs to 
Figure \ref{dividing-cut-vertex-figure}(iii). 
\qed

\subsection{Only non-dividing cut-vertices}

In this subsection, we will determine the rest of obstructions that contain 
cut-vertices, i.e., those containing only non-dividing cut-vertices. 

\begin{lemma}\label{non-dividing cut-vertex is second}
Let $X$ be an obstruction and $\prec: v_1, v_2, \dots, v_n$ be a straight 
enumeration of $U(X)$. Suppose that $v_c$ is a non-dividing cut-vertex. 
Then $c=2$ or $c=n-1$. Moreover, if $v_c$ is incident with both arcs then 
$n = 4$.
\end{lemma}
\pf Let $(v_j,v_k)$ ($j < k$) and $(v_s,v_t)$ ($s > t$) be the arcs in 
$X$. Since $v_c$ is non-dividing, either $c \leq \mbox{min}\{j,t\}$ or 
$c \geq \mbox{max}\{k,s\}$. Suppose that $c \leq \mbox{min}\{j,t\}$.

Let $(v_j, v_k) = (u_1,w_1), \dots,(u_q,w_q) = (v_t, v_s)$ be a $\Gamma$-sequence 
of $U(X)$ between $(v_j, v_k)$ and $(v_t, v_s)$. By Theorem 
\ref{different implication classes theorem}, the sequence must include
all vertices of $X$. Let $\alpha$ be the smallest subscript such that one of 
$u_{\alpha}, w_{\alpha}$ precedes $v_c$ (and hence the other vertex is $v_c$ 
since $v_c$ is a cut-vertex). Similarly, let $\beta$ be the largest subscript 
such that one of $u_\beta, w_\beta$ precedes $v_c$ (and hence the other vertex is 
$v_c$). Then it is easy to verify that 
$(u_1,w_1), \dots,(u_\alpha,w_\alpha), (u_{\beta+1},w_{\beta+1}), \dots, 
(u_q,w_q)$ is a $\Gamma$-sequence between $(v_j, v_k)$ and $(v_t, v_s)$. 
Since this sequence contains a unique vertex preceding $v_c$ and includes all 
vertices of $X$, we must have $c = 2$. A similar argument shows that if 
$c \geq \mbox{max}\{k,s\}$ then $c = n-1$. 

Suppose $v_c$ is incident with both arcs. Then either $c = j = t = 2$
or $c = k = s = n-1$. If $c = j = t = 2$, then 
$(v_j,v_k) \Gamma (v_1,v_j) \Gamma (v_t,v_s)$ and by Theorem 
\ref{different implication classes theorem}, 
$v_1, v_j = v_t, v_k, v_s$ are all the vertices of $X$ so $n = 4$.
A similar argument shows that $X$ has exactly four vertices if $c = k = s = n-1$.
\qed

The following theorem deals with the case when $v_2$ and $v_{n-1}$ are both
non-dividing cut-vertices of $U(X)$.

\begin{theorem}\label{2 non-dividing cut-vertices}
Let $X$ be an obstruction and $\prec: v_1, v_2, \dots, v_n$ be a straight 
enumeration of $U(X)$. Suppose that $v_2$ and $v_{n-1}$ are the two cut-vertices 
of $U(X)$, both non-dividing. Then $X$ or its dual is one of the two graphs in 
Figure~\ref{2non-dividing}.
	
	\begin{figure}[H]
		\centering
		\begin{subfigure}{0.5\textwidth}
			\centering
			\begin{tikzpicture}
			\node[unlabelled]	(1) at (0,0)	{};
			\node[unlabelled]	(2) at (1,0)	{};
			\node[unlabelled]	(3) at (2,0)	{};
			\node[unlabelled]	(4) at (3,0)	{};
			\node[unlabelled]	(5) at (4,0)	{};
			\draw	(1) -- (2)	(4) -- (5)
			(2) edge[bend left = 45] (4);
			\draw[->-=.5]	(2) to (3);
			\draw[->-=.5]	(4) to (3);
			\end{tikzpicture}
		\end{subfigure}%
		\begin{subfigure}{0.5\textwidth}
			\centering
			\begin{tikzpicture}
			\node[unlabelled]	(1) at (0,0)	{};
			\node[unlabelled]	(2) at (1,0)	{};
			\node[unlabelled]	(3) at (2,0)	{};
			\node[unlabelled]	(4) at (3,0)	{};
			\node[unlabelled]	(5) at (4,0)	{};
			\node[unlabelled]	(6) at (5,0)	{};
			\draw	(1) -- (2)	(5) -- (6)	(3) -- (4)
			(2) edge[bend left = 45] (4)
			(3) edge[bend left = 45] (5)
			(2) edge[bend left = 55] (5);
			\draw[->-=.5]	(2) to (3);
			\draw[->-=.5]	(5) to (4);
			\end{tikzpicture}
		\end{subfigure}%
\caption{Obstructions with two non-dividing cut-vertices.
\label{2non-dividing}}
	\end{figure}
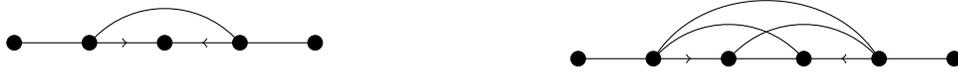
\end{theorem}
\pf Since both $v_2$ and $v_{n-1}$ are non-dividing cut-vertices, $n \geq 5$. Hence
by Lemma \ref{non-dividing cut-vertex is second} each of $v_2$ and $v_{n-1}$ 
is incident with at most one arc.

We show that $v_1$ and $v_n$ are arc-balancing vertices. By symmetry we only
prove that $v_1$ is arc-balancing. Clearly $v_1$ is not a cut-vertex of $U(X)$ and 
is not incident with an arc.
By Proposition \ref{vertex classification}, it can only be an arc-balancing 
vertex or a cut-vertex of $\overline{U(X)}$. 
Assume that $v_1$ is a cut-vertex of $\overline{U(X)}$. 
By Proposition \ref{universal vertex maker}, some vertex $v$ is 
adjacent to every vertex in $X$ except $v_1$. Since $v_{n-1}$ is the only 
neighbour of $v_n$ in $U(X)$, $v = v_{n-1}$.
Since the vertex $v = v_{n-1}$ is adjacent to $v_2$, by the umbrella property, 
the vertices $v_i$ with $2 \leq i \leq n-1$ induce a clique in $U(X)$.
Thus the vertices $v_i$ with $3 \leq i \leq n-2$ have the same closed 
neighbourhood in $U(X)$ and hence cannot contain both endvertices of any arc. 
It follows that each arc is incident with $v_2$ or $v_{n-1}$.  
From the above we know that each of $v_2$ and $v_{n-1}$ is incident with at most 
one arc. It is not possible that $v_2$ and $v_{n-1}$ are incident with the same
arc (as otherwise the endvertices of the other arc have the same closed 
neighbourhood). Hence $v_2$ and $v_{n-1}$ are incident with different arcs. 
We see that $v_1$ is an arc-balancing vertex. 

By taking the dual of $X$ if necessary we assume $(v_2,v_k)$ and $(v_{n-1},v_t)$ 
are the two arcs in $X$ where $3 \leq k, t \leq n-2$.
Then $v_1$ is the $(v_2,v_k)$-balancing vertex and $v_n$ is the 
$(v_{n-1},v_t)$-balancing vertex. No vertex $v_i$ with $2 < i < n-1$ is 
a cut-vertex of $U(X)$ or $\overline{U(X)}$ and hence each must be incident with 
an arc of $X$ by Proposition \ref{vertex classification}. Hence $v_k$ and $v_t$ 
are the only vertices between $v_2$ and $v_{n-1}$ in $\prec$. It is now easy to 
verify that $X$ is one of the two graphs in Figure~\ref{2non-dividing}.
\qed

It remains to consider the case when $X$ has only one cut-vertex and it is 
non-dividing. By Lemma~\ref{non-dividing cut-vertex is second} and reversing 
the straight enumeration $\prec$ if necessary we will assume $v_2$ is this vertex.

\begin{lemma}\label{only 1 cut-vertex lemma}
Let $X$ be an obstruction and $\prec: v_1, v_2, \dots, v_n$ be a straight 
enumeration of $U(X)$. Suppose that $v_2$ is the only cut-vertex and it is 
non-dividing. Then, the following statements hold:
\begin{enumerate}[label=(\alph*)]
\item For each $i \geq 3$, $v_i$ is an arc-balancing vertex or incident with an arc;
\item For some $i \geq 3$, $v_i$ is adjacent to every vertex except for $v_1$. 
      Moreover, there are at most two such vertices, each incident with exactly 
      one arc;
\item The number of vertices in $X$ is between 4 and 8 (i.e., $4 \leq n \leq 8$).
\end{enumerate}
\end{lemma}
\pf For (a), if each $v_i$ with $i \geq 3$ is an arc-balancing vertex or incident 
with an arc then we are done. Otherwise, by Proposition \ref{vertex classification},
some $v_i$ with $i \geq 3$ is a cut-vertex of $\overline{U(X)}$. According to 
Proposition \ref{universal vertex maker}, $v_i = v_n$ and there is a vertex 
adjacent to every vertex except $v_n$ in $U(X)$. Such a vertex can only be $v_2$. 
Since $v_{n-1}$ is not a cut-vertex of $U(X)$, $v_n$ is adjacent to $v_{n-2}$.
Since $v_2$ is not adjacent to $v_n$, $n-2 > 2$ (i.e., $n > 4$) and hence 
by Lemma \ref{non-dividing cut-vertex is second}, there is an arc which is not
incident with $v_2$. This arc must have endvertices strictly between $v_2$ and 
$v_n$ in $\prec$. Therefore $v_n$ is an arc-balancing vertex, which contradicts 
our assumption.

Statement (b) holds if $v_1$ is a cut-vertex of $\overline{U(X)}$. Indeed, by 
Proposition \ref{universal vertex maker} there is a vertex $v_i$ which is adjacent 
to every vertex except $v_1$ and it is clear that $i \geq 3$. So assume $v_1$ is 
not a cut-vertex of $\overline{U(X)}$. Since $v_2$ is the only cut-vertex and it is
non-dividing, $v_1$ is neither a cut-vertex of $U(X)$ nor incident with an arc, and
hence must be an arc-balancing vertex by Proposition \ref{vertex classification}.
Without loss of generality, assume $v_1$ balances an arc between $v_2$ and $v_j$
for some $j > 2$. If $v_j = v_n$ or $v_j$ is adjacent to $v_n$, then $v_j$ is 
adjacent to every vertex except $v_1$ and we are done. Otherwise, $j<n$ and $v_j$ 
is not adjacent to $v_n$. For each $j < k < n$, $v_k$ is not a cut-vertex of $U(X)$
by assumption so $v_{k-1}$ must be adjacent to $v_{k+1}$. Since $v_j$ is not 
adjacent to $v_n$, $j < n-2$ and thus $n > j+2 > 4$. By statement (a), each vertex 
$v_i$ with $i \geq 3$ is an arc-balancing vertex or incident with an arc. Since 
$v_1$ is arc-balancing and $v_2$ is incident with an arc, there are at most four 
vertices $v_i$ with $i \geq 3$. Hence $n \leq 6$ and therefore $n = 6$. It is now 
easy to see that $v_4$ is adjacent to every vertex except $v_1$.

Suppose $v_i$ with $i \geq 3$ is a vertex adjacent to every vertex except $v_1$. 
Clearly $v_i$ is not an arc-balancing vertex and hence by (a) it is incident 
with an arc. We show by contradiction that $v_i$ is incident with exactly one arc. 
So suppose that $v_i$ is incident with both arcs of $X$. 
Let $v_s$ and $v_t$ denote the other endvertices of the two arcs. 
We first show that either $s=2$ or $t=2$. 
By Theorem \ref{different implication classes theorem}, the edges $v_iv_s, v_iv_t$
 belong to different implication classes in $U(X-v_1)$. Since $v_i$ is 
an isolated vertex in $\overline{U(X-v_1)}$, each of $v_s, v_t, v_i$ belongs to 
a different component of $\overline{U(X-v_1)}$ by Theorem~\ref{pigstructure}.
In particular, one of $v_s, v_t$ is an isolated vertex in $\overline{U(X-v_1)}$. 
Without loss of generality, assume $v_s$ is such a vertex. Thus, $v_s$ is adjacent to every vertex except possibly $v_1$ in $X$. 
If $v_s$ is not adjacent to $v_1$, then $v_s$ and $v_i$ share the same closed 
neighbourhood, so the arc between $v_s$ and $v_i$ is balanced, a contradiction. 
Hence, $v_s$ is adjacent to $v_1$ and $v_s=v_2$. 
Consider $v_t$. Suppose $t < i$. Since $2 = s < t < i$ and $v_i$ is adjacent to $v_s$, the umbrella property implies $v_t$ is adjacent to $v_s$. If $v_t$ is also adjacent to $v_n$, then $v_i$ and $v_t$ have the same closed neighbourhood so the arc 
between them is balanced, a contradiction. Hence, $v_t$ is not adjacent to $v_n$. 
Since $s < t < n$, the umbrella property implies that $v_s$ and $v_n$ are not adjacent. Thus $(v_t, v_i)\Gamma(v_i,v_n)\Gamma(v_s,v_i)$ is a $\Gamma$-sequence 
between the arcs and not containing $v_1$, a contradiction by 
Theorem \ref{different implication classes theorem}. It follows that $i < t$. 
If $v_t$ is non-adjacent to $v_s$, then $(v_i,v_t)\Gamma(v_s,v_i)$ is 
a $\Gamma$-sequence between the arcs and not containing $v_1$, a contradiction. 
Hence, $v_t$ is adjacent to $v_s = v_2$. If $t=n$, then the arc between $v_i$ and 
$v_t$ is balanced by the umbrella property, a contradiction. If $t<n$, then $v_t$ 
is adjacent to $v_n$ because $i<t<n$ and $v_i$ is adjacent to $v_n$, leading to 
a similar contradiction. Therefore $v_i$ is incident with exactly one arc.
Suppose $v_i, v_j$ are two such vertices. By the above, each of them is incident 
with an arc. Moreover, they cannot be incident with the same arc because they 
share the same neighbourhood. Hence, they are each incident with a different arc. 
Since $X$ contains two arcs, there are at most two such vertices.

Finally we prove (c). Clearly, $n \geq 4$. Since there are at most four vertices 
incident with arcs and at most two arc-balancing vertices in $X$, there can be 
at most six vertices $v_i$ with $i \geq 3$ by (a). Therefore $n \leq 8$.
\qed

\begin{theorem} \label{non-dividing4-5}
Let $X$ be an obstruction and $\prec: v_1, v_2, \dots, v_n$ be a straight 
enumeration of $U(X)$. Suppose that $v_2$ is the only cut-vertex and it is 
non-dividing. If $n=4$ or 5, then $X$ or its dual is one of the graphs in 
Figure~\ref{1non-dividing4-5}.
	\begin{figure}[H]
\captionsetup[subfigure]{labelformat=empty}
\begin{subfigure}{.3\textwidth}
		\centering
		\begin{tikzpicture}
		\node[unlabelled] (a1) at (3,0) {};
		\node[unlabelled] (a2) at (2,0) {};
		\node[unlabelled] (x1) at (1,0) {};
		\node[unlabelled] (b1) at (0,0) {};
		\draw	(a1) -- (a2)	(x1) -- (b1);
		\draw[->-=.5]	(a1) to [bend right = 45] (x1);
		\draw[->-=.5]	(x1) to (a2);
		\end{tikzpicture}
\subcaption{(i)}
\end{subfigure}%
                \begin{subfigure}{.3\textwidth}
                        \centering
                        \begin{tikzpicture}
                        \node[unlabelled] (1) at (0,0)  {};
                        \node[unlabelled] (2) at (1,0)  {};
                        \node[unlabelled] (3) at (2,0)  {};
                        \node[unlabelled] (4) at (3,0)  {};
                        \node[unlabelled] (5) at (4,0)  {};
                        \draw   (1) -- (2)      (3) -- (4)
                                        (2) edge [bend left = 45] (4)
                                        (3) edge [bend left = 45] (5);
                        \draw[->-=.5]   (2) to (3);
                        \draw[->-=.5]   (5) to (4);
                        \end{tikzpicture}
                        \subcaption{(ii)}
                \end{subfigure}%
                \begin{subfigure}{.4\textwidth}
                        \centering
                        \begin{tikzpicture}
                        \node[unlabelled] (ai) at (4,0) {};
                        \node[unlabelled] (aj) at (3,0) {};
                        \node[unlabelled] (ak) at (2,0) {};
                        \node[unlabelled] (x1) at (1,0) {};
                        \node[unlabelled] (b1) at (0,0) {};
                        \draw   (aj) -- (ak)    (ak) -- (x1)
                        (x1) -- (b1)
                        (aj) edge [bend right = 45] (x1);
                        \draw[->-=.5]   (ak) to [bend left = 45] (ai);
                        \draw[->-=.5]   (ai) to (aj);
                        \end{tikzpicture}
                        \subcaption{(iii)}
                \end{subfigure}

\caption{Obstructions with a unique non-dividing cut-vertex on 4 or 5 vertices.
\label{1non-dividing4-5}}
	\end{figure}
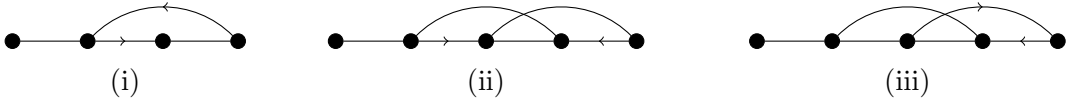
\end{theorem}
\pf Suppose $n = 4$. Since $v_3$ is not a cut-vertex, $v_2$ and $v_4$ are adjacent.
Both $v_3$ and $v_4$ are adjacent to every vertex except for $v_1$ and by 
Lemma \ref{only 1 cut-vertex lemma}(b) they are incident with different arcs.
It is easy to see that $X$ is Figure~\ref{1non-dividing4-5}(i).

Suppose $n = 5$. For each $i=3, 4$, $v_i$ is not a cut-vertex, so $v_{i-1}$ and 
$v_{i+1}$ are adjacent. On the other hand if $v_2$ is adjacent to $v_5$, then 
the umbrella property implies $v_3, v_4, v_5$ are all adjacent to every vertex 
except for $v_1$, contradicting Lemma \ref{only 1 cut-vertex lemma}(b). 
So $v_2$ and $v_5$ are not adjacent. Each of $v_3,v_4$ is adjacent to every vertex
except $v_1$ and by Lemma \ref{only 1 cut-vertex lemma}(b) they are incident with
different arcs. Since $n \neq 4$, $v_2$ is not incident with both arcs according to
Lemma \ref{non-dividing cut-vertex is second}. It follows that $v_5$ must be 
incident with at least one arc. If $v_5$ is incident with exactly one arc, then 
$X$ is or its dual is Figure~\ref{1non-dividing4-5}(ii). Otherwise $v_5$ is 
incident with both arcs and $X$ or its dual is Figure~\ref{1non-dividing4-5}(iii).
\qed

\begin{lemma} \label{pigbalance}
Let $X$ be an obstruction and $\prec: v_1, v_2, \dots, v_n$ be a straight 
enumeration of $U(X)$. Suppose $v_k$ is a $(v_i, v_j)$-balancing vertex. Then, 
either $k < \mbox{min}\{i,j\}$ or $k > \mbox{max}\{i,j\}$. Moreover,
\begin{itemize}
\item If $k < \mbox{min}\{i,j\}$, then no $v_p$ with $p < k$ is adjacent to 
      either one of $v_i, v_j$, and any $v_q$ with $q > \mbox{max}\{i,j\}$ is 
      adjacent to either both or neither of $v_i, v_j$;
\item If $k > \mbox{max}\{i,j\}$, then no $v_p$ with $p > k$ is adjacent to either 
      one of $v_i, v_j$, and any $v_q$ with $q < \mbox{min}\{i,j\}$ is adjacent to 
      either both or neither of $v_i, v_j$.
\end{itemize}
\end{lemma}

\pf First we show that either $k < \mbox{min}\{i,j\}$ or $k > \mbox{max}\{i,j\}$. 
Otherwise, $v_k$ is between $v_i$ and $v_j$. Since $v_iv_j$ is an edge of $U(X)$, 
the umbrella property implies that both $v_i$ and $v_j$ are adjacent to $v_k$, 
a contradiction to the fact that $v_k$ is a $(v_i,v_j)$-balancing vertex. Thus, 
either $k < \mbox{min}\{i,j\}$ or $k > \mbox{max}\{i,j\}$.

By symmetry, it suffices to consider the first case. 
Suppose $k < \mbox{min}\{i,j\}$. If $v_p$ with $p < k$ is adjacent to either one 
of $v_i, v_j$, then it must also be adjacent to $v_k$ by the umbrella property. 
Since $v_k$ is the only vertex adjacent to exactly one of $v_i, v_j$, $v_p$ must be
adjacent to both $v_i$ and $v_j$. By the umbrella property, both $v_i$ and $v_j$ 
are adjacent to $v_k$, a contradiction. On the other hand, since $v_k$ is the only 
vertex adjacent to exactly one of $v_i, v_j$, it is clear that any $v_q$ with 
$q > \mbox{max}\{i,j\}$ is adjacent to either both or neither of $v_i, v_j$.
\qed

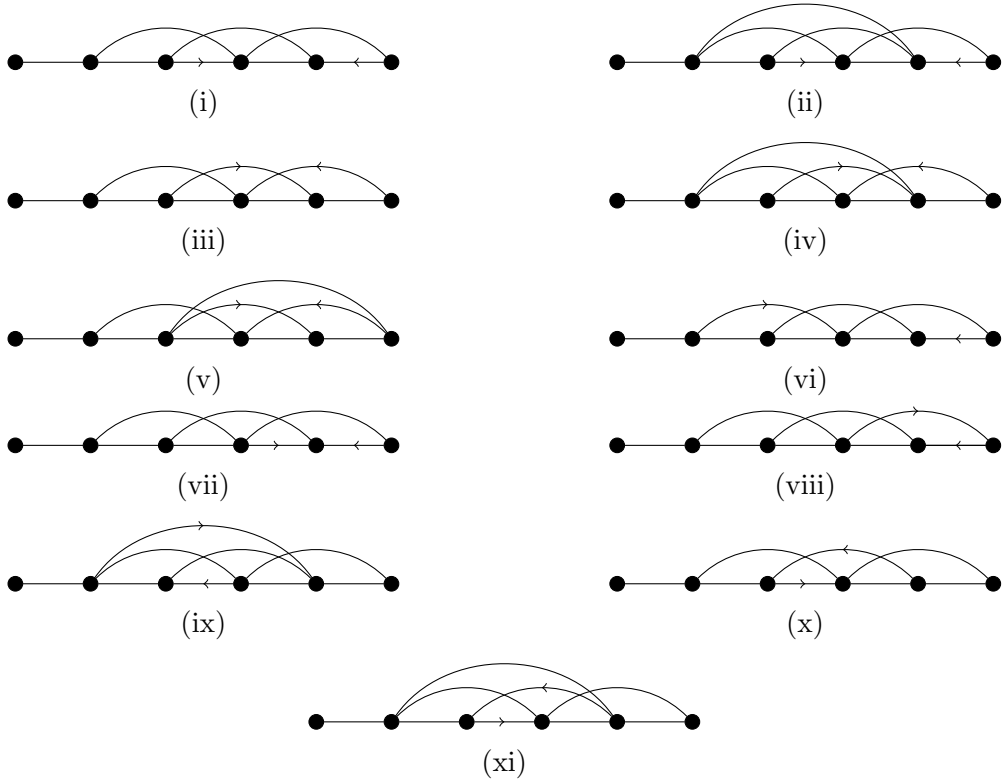
\begin{figure}[H]
		\captionsetup[subfigure]{labelformat=empty}
		\begin{subfigure}[b]{.5\textwidth}
			\centering
			\begin{tikzpicture}
			\node[unlabelled] (a1) at (0,0)	{};
			\node[unlabelled] (x1) at (1,0)	{};
			\node[unlabelled] (b1) at (2,0)	{};
			\node[unlabelled] (b2) at (3,0)	{};
			\node[unlabelled] (b3) at (4,0)	{};
			\node[unlabelled] (b4) at (5,0)	{};
			\draw	(a1) -- (x1)	(x1) -- (b1)	(b2) -- (b3)
			(x1) edge[bend left = 45] (b2)
			(b1) edge[bend left = 45] (b3)
			(b2) edge[bend left = 45] (b4);
			\draw[->-=.5]	(b1) to (b2);
			\draw[->-=.5]	(b4) to (b3);
			\end{tikzpicture}
			\subcaption{(i)}
		\end{subfigure}%
		\begin{subfigure}[b]{.5\textwidth}
			\centering
			\begin{tikzpicture}
			\node[unlabelled] (a1) at (0,0)	{};
			\node[unlabelled] (x1) at (1,0)	{};
			\node[unlabelled] (b1) at (2,0)	{};
			\node[unlabelled] (b2) at (3,0)	{};
			\node[unlabelled] (b3) at (4,0)	{};
			\node[unlabelled] (b4) at (5,0)	{};
			\draw	(a1) -- (x1)	(x1) -- (b1)	(b2) -- (b3)
			(x1) edge[bend left = 45] (b2)
			(b2) edge[bend left = 45] (b4)
			(b1) edge[bend left = 45] (b3)
			(x1) edge[bend left = 55] (b3);
			\draw[->-=.5]	(b1) to (b2);
			\draw[->-=.5]	(b4) to (b3);
			\end{tikzpicture}
			\subcaption{(ii)}
		\end{subfigure}
		
		\begin{subfigure}[b]{.5\textwidth}
			\centering
			\begin{tikzpicture}
			\node[unlabelled] (a1) at (0,0)	{};
			\node[unlabelled] (x1) at (1,0)	{};
			\node[unlabelled] (b1) at (2,0)	{};
			\node[unlabelled] (b2) at (3,0)	{};
			\node[unlabelled] (b3) at (4,0)	{};
			\node[unlabelled] (b4) at (5,0)	{};
			\draw	(a1) -- (x1)	(x1) -- (b1)	(b2) -- (b3)	(b3) -- (b4)
			(b1) -- (b2)
			(x1) edge[bend left = 45] (b2);
			\draw[->-=.5]	(b1) to [bend left = 45] (b3);
			\draw[->-=.5]	(b4) to [bend right = 45] (b2);
			\end{tikzpicture}
			\subcaption{(iii)}
		\end{subfigure}%
		\begin{subfigure}[b]{.5\textwidth}
			\centering
			\begin{tikzpicture}
			\node[unlabelled] (a1) at (0,0)	{};
			\node[unlabelled] (x1) at (1,0)	{};
			\node[unlabelled] (b1) at (2,0)	{};
			\node[unlabelled] (b2) at (3,0)	{};
			\node[unlabelled] (b3) at (4,0)	{};
			\node[unlabelled] (b4) at (5,0)	{};
			\draw	(a1) -- (x1)	(x1) -- (b1)	(b2) -- (b3)	(b1) -- (b2)	(b3) -- (b4)
			(x1) edge[bend left = 45] (b2)
			(x1) edge[bend left = 55] (b3);
			\draw[->-=.5]	(b1) to [bend left = 45] (b3);
			\draw[->-=.5]	(b4) to [bend right = 45] (b2);
			\end{tikzpicture}
			\subcaption{(iv)}
		\end{subfigure}
		
		\begin{subfigure}[b]{.5\textwidth}
			\centering
			\begin{tikzpicture}
			\node[unlabelled] (6) at (5,0)	{};
			\node[unlabelled] (5) at (4,0)	{};
			\node[unlabelled] (4) at (3,0)	{};
			\node[unlabelled] (3) at (2,0)	{};
			\node[unlabelled] (2) at (1,0)	{};
			\node[unlabelled] (1) at (0,0)	{};
			\draw	(1) -- (2)	(2) -- (3)	(3) -- (4)
					(5) -- (6)	(4) -- (5)
					(2) edge [bend left = 45] (4)
					(3) edge [bend left = 55] (6);
			\draw[->-=.5]	(3) to [bend left = 45] (5);
			\draw[->-=.5]	(6) to [bend right = 45] (4);
			\end{tikzpicture}
			\subcaption{(v)}
		\end{subfigure}%
		\begin{subfigure}[b]{.5\textwidth}
			\centering
			\begin{tikzpicture}
			\node[unlabelled] (v) at (0,0)	{};
			\node[unlabelled] (a) at (1,0)	{};
			\node[unlabelled] (b1) at (2,0)	{};
			\node[unlabelled] (b) at (3,0)	{};
			\node[unlabelled] (d) at (4,0)	{};
			\node[unlabelled] (c) at (5,0)	{};
			\draw	(v) -- (a)	(a) -- (b1)	(b1) -- (b)
			(b) -- (d)
			(b1) edge[bend left = 45] (d)
			(b) edge[bend left = 45] (c);
			\draw[->-=.5]	(a) to [bend left = 45] (b);
			\draw[->-=.5]	(c) to (d);
			\end{tikzpicture}
			\subcaption{(vi)}
		\end{subfigure}
		
		\begin{subfigure}[b]{.5\textwidth}
			\centering
			\begin{tikzpicture}
			\node[unlabelled] (a1) at (0,0)	{};
			\node[unlabelled] (x1) at (1,0)	{};
			\node[unlabelled] (b1) at (2,0)	{};
			\node[unlabelled] (b2) at (3,0)	{};
			\node[unlabelled] (b3) at (4,0)	{};
			\node[unlabelled] (b4) at (5,0)	{};
			\draw	(a1) -- (x1)	(x1) -- (b1)
			(b1) -- (b2)
			(x1) edge[bend left = 45] (b2)
			(b1) edge[bend left = 45] (b3)
			(b2) edge[bend left = 45] (b4);
			\draw[->-=.5]	(b2) to (b3);
			\draw[->-=.5]	(b4) to (b3);
			\end{tikzpicture}
			\subcaption{(vii)}
		\end{subfigure}%
		\begin{subfigure}[b]{.5\textwidth}
			\centering
			\begin{tikzpicture}
			\node[unlabelled] (a1) at (0,0)	{};
			\node[unlabelled] (x1) at (1,0)	{};
			\node[unlabelled] (b1) at (2,0)	{};
			\node[unlabelled] (b2) at (3,0)	{};
			\node[unlabelled] (b3) at (4,0)	{};
			\node[unlabelled] (b4) at (5,0)	{};
			\draw	(a1) -- (x1)	(x1) -- (b1)	(b2) -- (b3)	(b3) -- (b4)
			(b1) -- (b2)
			(x1) edge[bend left = 45] (b2)
			(b1) edge[bend left = 45] (b3);
			\draw[->-=.5]	(b4) to (b3);
			\draw[->-=.5]	(b2) to [bend left = 45] (b4);
			\end{tikzpicture}
			\subcaption{(viii)}
		\end{subfigure}
		
		\begin{subfigure}[b]{.5\textwidth}
			\centering
			\begin{tikzpicture}
			\node[unlabelled] (1) at (0,0)	{};
			\node[unlabelled] (2) at (1,0)	{};
			\node[unlabelled] (3) at (2,0)	{};
			\node[unlabelled] (4) at (3,0)	{};
			\node[unlabelled] (5) at (4,0)	{};
			\node[unlabelled] (6) at (5,0)	{};
			\draw	(1) edge (2)	(2) edge (3)	(4) edge (5)
			(5) edge (6)
			(2) edge[bend left = 45] (4)
			(3) edge[bend left = 45] (5)
			(4) edge[bend left = 45] (6);
			\draw[->-=.5]	(2) to [bend left = 55] (5);
			\draw[->-=.5]	(4) to (3);
			\end{tikzpicture}
			\subcaption{(ix)}
		\end{subfigure}%
		\begin{subfigure}[b]{.5\textwidth}
			\centering
			\begin{tikzpicture}
			\node[unlabelled] (a1) at (0,0)	{};
			\node[unlabelled] (x1) at (1,0)	{};
			\node[unlabelled] (b1) at (2,0)	{};
			\node[unlabelled] (b2) at (3,0)	{};
			\node[unlabelled] (b3) at (4,0)	{};
			\node[unlabelled] (b4) at (5,0)	{};
			\draw	(a1) -- (x1)	(x1) -- (b1)	(b2) -- (b3)	(b3) -- (b4)
			(x1) edge[bend left = 45] (b2)
			(b2) edge[bend left = 45] (b4);
			\draw[->-=.5]	(b1) to (b2);
			\draw[->-=.5]	(b3) to [bend right = 45] (b1);
			\end{tikzpicture}
			\subcaption{(x)}
		\end{subfigure}
		
		\begin{subfigure}{\textwidth}
			\centering
			\begin{tikzpicture}
			\node[unlabelled] (a1) at (0,0)	{};
			\node[unlabelled] (x1) at (1,0)	{};
			\node[unlabelled] (b1) at (2,0)	{};
			\node[unlabelled] (b2) at (3,0)	{};
			\node[unlabelled] (b3) at (4,0)	{};
			\node[unlabelled] (b4) at (5,0)	{};
			\draw	(a1) -- (x1)	(x1) -- (b1)	(b2) -- (b3)	(b3) -- (b4)
			(x1) edge[bend left = 45] (b2)
			(b2) edge[bend left = 45] (b4)
			(x1) edge[bend left = 55] (b3);
			\draw[->-=.5]	(b1) to (b2);
			\draw[->-=.5]	(b3) to [bend right = 45] (b1);
			\end{tikzpicture}
			\subcaption{(xi)}
		\end{subfigure}
\caption{Obstructions with a unique non-dividing cut-vertex on 6 vertices.
\label{1non-dividing6}}
	\end{figure}

\begin{theorem} \label{non-dividing6}
Let $X$ be an obstruction and $\prec: v_1, v_2, \dots, v_n$ be a straight 
enumeration of $U(X)$. Suppose that $v_2$ is the only cut-vertex and it is 
non-dividing. If $n = 6$, then $X$ or its dual is one of the graphs in
Figure~\ref{1non-dividing6}.
\end{theorem}
\pf For each $3 \leq i \leq 5$, $v_i$ is not a cut-vertex, so $v_{i-1}$ and
$v_{i+1}$ are adjacent. Now $v_2v_5$ and $v_3v_6$ cannot both be edges in $U(X)$
as otherwise each of $v_3, v_4$ and $v_5$ is adjacent to every vertex except for
$v_1$, contradicting Lemma \ref{only 1 cut-vertex lemma}(b).

We claim that each of $v_4$ and $v_5$ is incident with an arc. Since $v_4$ is
adjacent to every vertex except for $v_1$, it is incident with exactly one arc
by Lemma \ref{only 1 cut-vertex lemma}(b). On the other hand, suppose $v_5$ is
not incident with an arc. By Lemma \ref{only 1 cut-vertex lemma}(a), $v_5$ is
an arc-balancing vertex for some arc. Thus $v_5$ is adjacent to exactly one
endvertex of the arc. It is easy to see that the other endvertex can only be
$v_2$. Since $v_2$ is a cut-vertex, $v_1$ is adjacent to exactly one endvertex
(i.e., $v_2$) of the arc, a contradiction to that $v_5$ is arc-balancing for the
arc. Hence $v_5$ is incident with an arc.

Suppose $v_3$ and $v_6$ are also incident with arcs. Then $v_3, v_4, v_5, v_6$ are
endvertices of the two arcs. Suppose that the two arcs are between $v_3$ and $v_4$
and between $v_5$ and $v_6$ respectively. Then $v_3v_6$ is not an edge of
$U(X)$ as otherwise the arc between $v_5$ and $v_6$ is balanced, a contradiction.
If $v_2v_5$ is not an edge of $U(X)$ then $X$ or its dual is
Figure~\ref{1non-dividing6}(i); otherwise, $X$ or its dual is
Figure~\ref{1non-dividing6}(ii).
Suppose that the two arcs are between $v_3$ and $v_5$ and between $v_4$ and $v_6$
respectively. Then $X$ or its dual is Figure~\ref{1non-dividing6}(iii), (iv) or (v)
depending whether or not $v_2v_5$ and $v_3v_6$ are edges of $U(X)$.
Suppose the two arcs are between $v_3$ and $v_6$ and between $v_4$ and $v_5$
respectively. Then $X$ or its dual is again Figure~\ref{1non-dividing6}(v) (with
$v_5$ and $v_6$ being switched).

Suppose $v_3$ is not incident with an arc. 
By Lemma \ref{only 1 cut-vertex lemma}(a), $v_3$ is an arc-balancing vertex. 
By Lemma \ref{pigbalance}, $v_3$ balances an arc between $v_5$ and $v_6$. 
Since $v_4$ is incident with an arc and $v_3$ is not, the arc incident with $v_4$ 
has the other endvertex being $v_2, v_5, v_6$. These three cases are 
represented by Figure~\ref{1non-dividing6}(vi), (vii) and (viii).

It follows from the above that at least one of $v_3$ and $v_6$ is incident with an 
arc. Thus it remains to consider the case that $v_3$ is incident with an arc
but $v_6$ is not. By Lemma \ref{only 1 cut-vertex lemma}(a), $v_6$ is 
an arc-balancing vertex for some arc. By Lemma \ref{pigbalance}, $v_2$ cannot be 
an endvertex of this arc, so the arc must be between $v_3$ and one of $v_4, v_5$. 
In particular, this implies $v_3v_6$ is not an edge of $U(X)$. It is now easy to 
verify that $X$ or its dual is Figure~\ref{1non-dividing6}(ix), (x) or (xi).
\qed

\begin{lemma}\label{at least 7}
Let $X$ be an obstruction and $\prec: v_1, v_2, \dots, v_n$ be a straight 
enumeration of $U(X)$. Suppose that $v_2$ is the only cut-vertex and it is 
non-dividing. If $n \geq 7$, then $v_2$ is not incident with an arc and
the subgraph of $U(X)$ induced by the vertices $v_i$ with $i \geq 3$ cannot 
contain a copy of $K_5$.
\end{lemma}
\pf Suppose that $v_2$ is incident with an arc. Then $v_1$ is adjacent to exactly
one endvertex of this arc so this arc cannot be balanced by any vertex $v_i$ with
$i \geq 3$. It follows that there is at most one arc-balancing vertex $v_i$ with
$i \geq 3$. By Lemma \ref{only 1 cut-vertex lemma}(a) and the assumption $n \geq 7$
there are at least four vertices $v_i$ with $i \geq 3$ which are incident with
arcs, which is impossible because $v_2$ is such a vertex. 

By Lemma \ref{only 1 cut-vertex lemma}(a) and (c), $n \leq 8$ and each $v_i$ with 
$i \geq 3$ is an arc-balancing vertex or incident with an arc. Since neither of
$v_1,v_2$ is incident with an arc, any set of five vertices $v_i$ with $i \geq 3$ 
must contain an arc-balancing triple and hence cannot induce a copy of $K_5$
in $U(X)$.  
\qed

\begin{theorem} \label{non-dividing7}
Let $X$ be an obstruction and $\prec: v_1, v_2, \dots, v_n$ be a straight 
enumeration of $U(X)$. Suppose that $v_2$ is the only cut-vertex and 
it is non-dividing. If $n = 7$, then $X$ or its dual is one of the graphs in
Figure~\ref{1non-dividing7}.
	\begin{figure}[H]
		\captionsetup[subfigure]{labelformat=empty}
		\begin{subfigure}[b]{.5\textwidth}
			\centering
			\begin{tikzpicture}
			\node[unlabelled] (a1) at (0,0)	{};
			\node[unlabelled] (x1) at (1,0)	{};
			\node[unlabelled] (b1) at (2,0)	{};
			\node[unlabelled] (b2) at (3,0)	{};
			\node[unlabelled] (b3) at (4,0)	{};
			\node[unlabelled] (b4) at (5,0)	{};
			\node[unlabelled] (b5) at (6,0)	{};
			\draw	(a1) -- (x1)	(x1) -- (b1)	(b2) -- (b3)
			(b1) -- (b2)	(b4) -- (b5)	(b3) -- (b4)
			(x1) edge[bend left = 45] (b2)
			(b1) edge[bend left = 45] (b3)
			(b2) edge[bend left = 45] (b4)
			(b3) edge[bend left = 45] (b5)
			(x1) edge[bend left = 55] (b3)
			(b1) edge[bend left = 55] (b4)
			(b2) edge[bend left = 55] (b5);
			%\draw[->-=.5]	(b1) to [bend left = 45] (b3);
			%\draw[->-=.5]	(b4) to [bend right = 45] (b2);
			\draw[->-=.5]	(b1) to (b2);
			\draw[->-=.5]	(b4) to (b3);
			\end{tikzpicture}
			\subcaption{(i)}
		\end{subfigure}%
		\begin{subfigure}[b]{.5\textwidth}
			\centering
			\begin{tikzpicture}
			\node[unlabelled] (a1) at (0,0)	{};
			\node[unlabelled] (x1) at (1,0)	{};
			\node[unlabelled] (b1) at (2,0)	{};
			\node[unlabelled] (b2) at (3,0)	{};
			\node[unlabelled] (b3) at (4,0)	{};
			\node[unlabelled] (b4) at (5,0)	{};
			\node[unlabelled] (b5) at (6,0)	{};
			\draw	(a1) -- (x1)	(x1) -- (b1)	(b2) -- (b3)
			(b1) -- (b2)	(b4) -- (b5)	(b3) -- (b4)
			(x1) edge[bend left = 45] (b2)
			(b1) edge[bend left = 45] (b3)
			(b2) edge[bend left = 45] (b4)
			(b3) edge[bend left = 45] (b5)
			(b1) edge[bend left = 55] (b4)
			(b2) edge[bend left = 55] (b5);
			\draw[->-=.5] (b2) to (b3);
			\draw[->-=.5] (b5) to (b4);
			\end{tikzpicture}
			\subcaption{(ii)}
		\end{subfigure}
		
		\begin{subfigure}[b]{.5\textwidth}
			\centering
			\begin{tikzpicture}
			\node[unlabelled] (a1) at (0,0)	{};
			\node[unlabelled] (x1) at (1,0)	{};
			\node[unlabelled] (b1) at (2,0)	{};
			\node[unlabelled] (b2) at (3,0)	{};
			\node[unlabelled] (b3) at (4,0)	{};
			\node[unlabelled] (b4) at (5,0)	{};
			\node[unlabelled] (b5) at (6,0)	{};
			\draw	(a1) -- (x1)	(x1) -- (b1)	(b2) -- (b3)
			(b1) -- (b2)	(b4) -- (b5)	(b3) -- (b4)
			(x1) edge[bend left = 45] (b2)
			(b1) edge[bend left = 45] (b3)
			(b2) edge[bend left = 55] (b5);
			\draw[->-=.5]	(b2) to [bend left = 45] (b4);
			\draw[->-=.5]	(b5) to [bend right = 45] (b3);
			\end{tikzpicture}
			\subcaption{(iii)}
		\end{subfigure}%
		\begin{subfigure}[b]{.5\textwidth}
			\centering
			\begin{tikzpicture}
			\node[unlabelled] (a1) at (0,0)	{};
			\node[unlabelled] (x1) at (1,0)	{};
			\node[unlabelled] (b1) at (2,0)	{};
			\node[unlabelled] (b2) at (3,0)	{};
			\node[unlabelled] (b3) at (4,0)	{};
			\node[unlabelled] (b4) at (5,0)	{};
			\node[unlabelled] (b5) at (6,0)	{};
			\draw	(a1) -- (x1)	(x1) -- (b1)
			(b1) -- (b2)	(b4) -- (b5)	(b3) -- (b4)	(b2) -- (b3)
			(x1) edge[bend left = 45] (b2)
			(b1) edge[bend left = 45] (b3)
			(b2) edge[bend left = 45] (b4)
			(b3) edge[bend left = 45] (b5)
			(x1) edge[bend left = 55] (b3)
			(b1) edge[bend left = 55] (b4);
			\draw[->-=.5]	(b1) to [bend left = 45] (b3);
			\draw[->-=.5]	(b4) to [bend right = 45] (b2);
			\end{tikzpicture}
			\subcaption{(iv)}
		\end{subfigure}%
		
		\begin{subfigure}[b]{.5\textwidth}
			\centering
			\begin{tikzpicture}
			\node[unlabelled] (a1) at (0,0)	{};
			\node[unlabelled] (x1) at (1,0)	{};
			\node[unlabelled] (b1) at (2,0)	{};
			\node[unlabelled] (b2) at (3,0)	{};
			\node[unlabelled] (b3) at (4,0)	{};
			\node[unlabelled] (b4) at (5,0)	{};
			\node[unlabelled] (b5) at (6,0)	{};
			\draw	(a1) -- (x1)	(x1) -- (b1)	(b2) -- (b3)
			(b1) -- (b2)	(b4) -- (b5)	(b3) -- (b4)
			(x1) edge[bend left = 45] (b2)
			(b1) edge[bend left = 55] (b4)
			(b3) edge[bend left = 45] (b5)
			(x1) edge[bend left = 55] (b3)
			(x1) edge[bend left = 65] (b4);
			\draw[->-=.5]	(b1) to [bend left = 45] (b3);
			\draw[->-=.5]	(b4) to [bend right = 45] (b2);
			\end{tikzpicture}
			\subcaption{(v)}
		\end{subfigure}%
		\begin{subfigure}[b]{.5\textwidth}
			\centering
			\begin{tikzpicture}
			\node[unlabelled] (a1) at (0,0)	{};
			\node[unlabelled] (x1) at (1,0)	{};
			\node[unlabelled] (b1) at (2,0)	{};
			\node[unlabelled] (b2) at (3,0)	{};
			\node[unlabelled] (b3) at (4,0)	{};
			\node[unlabelled] (b4) at (5,0)	{};
			\node[unlabelled] (b5) at (6,0)	{};
			\draw	(a1) -- (x1)	(x1) -- (b1)	(b2) -- (b3)
			(b3) -- (b4)	(b1) -- (b2)
			(x1) edge[bend left = 45] (b2)
			(b2) edge[bend left = 45] (b4)
			(b3) edge[bend left = 45] (b5)
			(x1) edge[bend left = 55] (b3);
			\draw[->-=.5]	(b5) to (b4);
			\draw[->-=.5]	(b1) to [bend left = 45] (b3);
			\end{tikzpicture}
			\subcaption{(vi)}
		\end{subfigure}%

		\begin{subfigure}[b]{.5\textwidth}
			\centering
			\begin{tikzpicture}
			\node[unlabelled] (a1) at (0,0)	{};
			\node[unlabelled] (x1) at (1,0)	{};
			\node[unlabelled] (b1) at (2,0)	{};
			\node[unlabelled] (b2) at (3,0)	{};
			\node[unlabelled] (b3) at (4,0)	{};
			\node[unlabelled] (b4) at (5,0)	{};
			\node[unlabelled] (b5) at (6,0)	{};
			\draw	(a1) -- (x1)	(x1) -- (b1)	
			(b3) -- (b4)	(b4) -- (b5)
			(x1) edge[bend left = 45] (b2)
			(b1) edge[bend left = 45] (b3)
			(b2) edge[bend left = 45] (b4)
			(b3) edge[bend left = 45] (b5)
			(x1) edge[bend left = 55] (b3);
			\draw[->-=.5]	(b1) to (b2);
			\draw[->-=.5]	(b3) to (b2);
			\end{tikzpicture}
			\subcaption{(vii)}
		\end{subfigure}%
		\begin{subfigure}[b]{.5\textwidth}
			\centering
			\begin{tikzpicture}
			\node[unlabelled] (a1) at (0,0)	{};
			\node[unlabelled] (x1) at (1,0)	{};
			\node[unlabelled] (b1) at (2,0)	{};
			\node[unlabelled] (b2) at (3,0)	{};
			\node[unlabelled] (b3) at (4,0)	{};
			\node[unlabelled] (b4) at (5,0)	{};
			\node[unlabelled] (b5) at (6,0)	{};
			\draw	(a1) -- (x1)	(x1) -- (b1)	(b2) -- (b3)
			(b3) -- (b4)	(b4) -- (b5)
			(x1) edge[bend left = 45] (b2)
			(b1) edge[bend left = 45] (b3)
			(b2) edge[bend left = 45] (b4)
			(x1) edge[bend left = 55] (b3);
			\draw[->-=.5]	(b1) to (b2);
			\draw[->-=.5]	(b5) to [bend right = 45] (b3);
			\end{tikzpicture}
			\subcaption{(viii)}
		\end{subfigure}
\caption{Obstructions with a unique non-dividing cut-vertex on 7 vertices.
\label{1non-dividing7}}
	\end{figure}
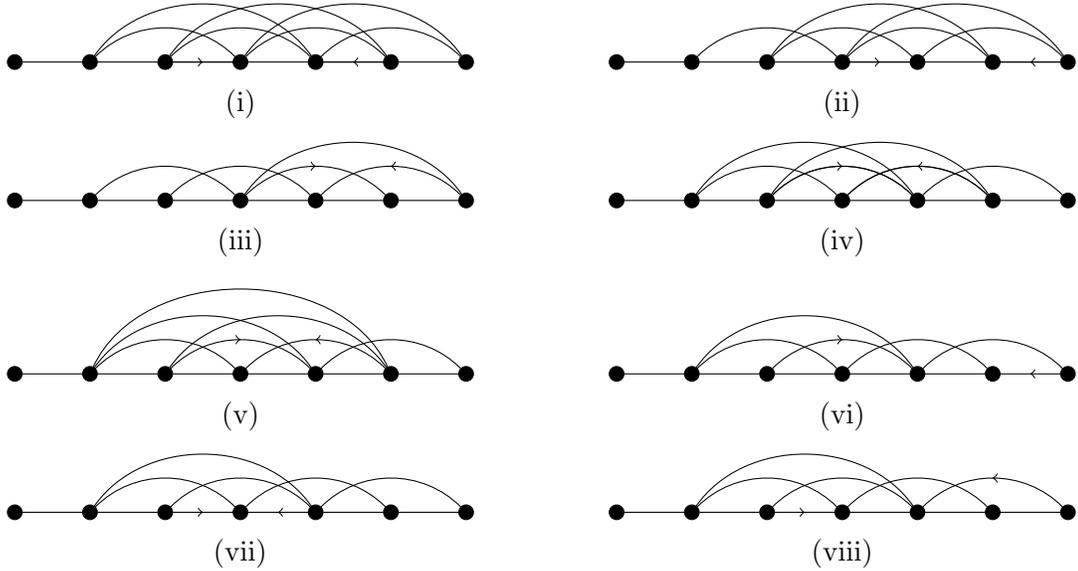
\end{theorem}
\pf First, note that $v_3v_7$ is not an edge in $U(X)$ as otherwise the vertices 
$v_i$ with $i \geq 3$ induce a $K_5$ in $U(X)$, a contradiction to 
Lemma \ref{at least 7}. If $v_2v_6$ and $v_4v_7$ are both edges of $U(X)$, then
each of $v_4, v_5, v_6$ is adjacent to every vertex except for $v_1$, 
contradicting Lemma \ref{only 1 cut-vertex lemma}(c). So, $v_2v_6$ and 
$v_4v_7$ cannot both be edges in $U(X)$.

By Lemma \ref{only 1 cut-vertex lemma}(b), there exists a vertex $v_i$ with 
$i \geq 3$ adjacent to every vertex except for $v_1$. Since $v_3v_7$
is not an edge in $U(X)$, neither $v_3$ nor $v_7$ is such a vertex. It is easy to 
see that if $v_6$ is such a vertex, then $v_5$ is also such a vertex. Hence, 
at least one of $v_4, v_5$ is adjacent to every vertex except for $v_1$.

Suppose $v_4$ is adjacent to every vertex except for $v_1$. This implies in 
particular that $v_4v_7$ is an edge of $U(X)$ and thus $v_2v_6$ is not an edge of 
$U(X)$. By Lemma \ref{only 1 cut-vertex lemma}(b), $v_4$ is incident with exactly 
one arc. Lemma \ref{at least 7} implies the other endvertex of this arc is one of 
$v_3, v_5, v_6$, and $v_7$. Suppose that the other endvertex is $v_3$. If no $v_i$ 
with $i \geq 5$ is an arc-balancing vertex for this arc, then $\{v_5, v_6, v_7\}$ 
must be an arc-balancing triple, a contradiction because these vertices induce 
a clique. Hence for some $i \geq 5$, $v_i$ is an arc-balancing vertex for the arc 
between $v_4$ and $v_3$. By Lemma \ref{pigbalance}, it must be $v_7$. 
Since $v_7$ balances the arc between $v_4$ and $v_3$, we see that $v_3$ must be 
adjacent to $v_6$. Both $v_5$ and $v_6$ are adjacent to $v_i$ for each 
$i \geq 3$ so they cannot be arc-balancing vertices. Hence by  
Lemma \ref{only 1 cut-vertex lemma}(a), both $v_5$ and $v_6$ are incident with 
arcs. This means there is an arc between $v_5$ and $v_6$, which implies that
$v_5$ is adjacent to $v_2$ (as otherwise $v_5$ and $v_6$ have the same closed 
neighbourhood in $U(X)$). Since $v_3v_7$ is not an edge in $U(X)$, $X$ or its dual 
is Figure~\ref{1non-dividing7}(i). Suppose next that there is an arc between 
$v_4$ and $v_5$. Since $v_4$ and $v_5$ cannot have the same closed neighbourhood
in $U(X)$, $v_2v_5$ is not an edge in $U(X)$. Clearly, the arc between $v_4$ and 
$v_5$ is not balanced by any of $v_3, v_6, v_7$, so $\{v_3, v_6, v_7\}$ is 
an arc-balancing triple. By Lemma \ref{pigbalance}, the arc is between $v_6$ and 
$v_7$. It follows that $v_3v_6$ is an edge in $U(X)$. Hence $X$ or its dual is 
Figure~\ref{1non-dividing7}(ii). 

Suppose next that there is an arc between $v_4$ and $v_6$. 
By Lemma \ref{pigbalance}, the arc between $v_4$ and $v_6$ cannot be balanced
by $v_3, v_5, v_7$. Similarly as above, $\{v_3, v_5, v_7\}$ is an arc-balancing 
triple. If $v_7$ balances an arc between $v_3$ and $v_5$, then $v_2v_5$ and 
$v_3v_6$ are edges in $U(X)$, and $X$ or its dual is 
Figure~\ref{1non-dividing7}(i). 
Suppose that $v_3$ balances an arc between $v_5$ and $v_7$. Each vertex except 
$v_3$ is either adjacent to both $v_5, v_7$ or neither. Since $v_2$ is not adjacent
to $v_7$, it is not adjacent to $v_5$. Hence, $X$ or its dual is 
Figure~\ref{1non-dividing7}(ii) or (iii) depending whether or not $v_3v_6$ is 
an edge of $U(X)$. Finally, suppose there is an arc between $v_4$ and $v_7$. 
By Lemma \ref{pigbalance}, none of $v_3, v_5, v_6$ is an arc-balancing vertex for 
this arc. Hence, $\{v_3, v_5, v_6\}$ is an arc-balancing triple. 
By Lemma \ref{pigbalance}, $v_3$ balances an arc between $v_5$ and $v_6$. 
It follows that neither $v_2v_5$ nor $v_3v_6$ can be an edge in $U(X)$. So $X$ or 
its dual is Figure~\ref{1non-dividing7}(iii).

Suppose now $v_4$ is not adjacent to every vertex except for $v_1$. From the above
we know that $v_5$ must be adjacent to every vertex except for $v_1$. So $v_5$ is 
incident with exactly one arc, and the other endvertex of this arc is one of 
$v_3, v_4, v_6$, and $v_7$. We claim that it cannot be $v_6$. Suppose to 
the contrary that there is an arc between $v_5$ and $v_6$. 
By Lemma \ref{pigbalance}, none of $v_3, v_4, v_7$ can be an arc-balancing vertex 
for this arc. Hence, $\{v_3, v_4, v_7\}$ is an arc-balancing triple.
Since neither $v_4v_7$ nor $v_3v_7$ is an edge of $U(X)$, the second arc can
only be between $v_3$ and $v_4$ but it is not balanced by $v_7$, a contradiction.
Hence, there is an arc between $v_5$ and one of $v_3, v_4$, and $v_7$. 

Suppose first that there is an arc between $v_5$ and $v_3$. Assume that this
arc is balanced by a vertex. By Lemma \ref{pigbalance}, it is balanced by $v_7$. 
It follows that $v_3v_6$ is an edge in $U(X)$. Since $v_6$ is adjacent to every 
vertex $v_i$ with $i \geq 3$, which are where all endvertices of arcs are, it 
cannot be an arc-balancing vertex. It follows that $v_6$ is incident with an arc. 
We claim that the other endvertex of this arc is $v_4$. Indeed, if it is not
$v_4$, then $v_4$ would be the arc-balancing vertex for this arc, a contradiction 
by Lemma \ref{pigbalance}. Thus, $X$ or its dual is Figure~\ref{1non-dividing7}(iv) 
or (v) depending whether or not $v_2v_6$ is an edge of $U(X)$. 
Assume now that the arc between $v_5$ and $v_3$ is not balanced by any vertex. 
In this case, $\{v_4, v_6, v_7\}$ is an arc-balancing triple. 
By Lemma \ref{pigbalance}, either $v_4$ balances an arc between $v_6$ and $v_7$, 
or $v_7$ balances an arc between $v_4$ and $v_6$. In the first case, $v_3v_6$ 
cannot be an edge of $U(X)$, as that would imply $v_3v_7$ is also an edge, 
a contradiction. Hence, $X$ or its dual is Figure~\ref{1non-dividing7}(vi). 
In the second case, $v_2v_6$ must be an edge of $U(X)$ and $X$ or its dual is 
Figure~\ref{1non-dividing7}(v). 

Suppose there is an arc between $v_5$ and $v_4$. We claim that $v_3$ is not
arc-balancing vertex. Indeed, if it is, then it must balance an arc between $v_6$ 
and $v_7$. Thus, $v_4v_7$ is an edge of $U(X)$, a contradiction. Hence, $v_3$ is 
incident with an arc. The other endvertex of this arc is $v_4, v_6$, or $v_7$. 
Clearly it cannot be $v_7$ because that would imply $v_4v_7$ is an edge of $U(X)$, 
contradicting the fact that $v_4$ is not adjacent to every vertex except for $v_1$.
Suppose the second arc is between $v_3$ and $v_4$. Then, $v_6$ must be an 
arc-balancing vertex. Clearly, $v_6$ cannot balance the arc between $v_5$ and 
$v_4$, so it must balance the arc between $v_3$ and $v_4$. It follows that $v_3v_6$
is not an edge of $U(X)$, so $X$ or its dual is Figure~\ref{1non-dividing7}(vii). 
On the other hand, suppose the second arc is between $v_3$ and $v_6$. 
In this case, $X$ or its dual is Figure~\ref{1non-dividing7}(iv) or (v) 
depending whether $v_2v_6$ is an edge of $U(X)$. 

Finally, suppose there is an arc between $v_5$ and $v_7$. Clearly, none of 
$v_3, v_4, v_6$ can be an arc-balancing vertex for this arc. Hence 
$\{v_3, v_4, v_6\}$ is an arc-balancing triple. By Lemma \ref{pigbalance}, 
$v_6$ must balance an arc between $v_3$ and $v_4$. It follows that $v_3v_6$ is not 
an edge of $U(X)$, so $X$ or its dual is Figure~\ref{1non-dividing7}(viii).
\qed

\begin{theorem} \label{non-dividing8}
Let $X$ be an obstruction and $\prec: v_1, v_2, \dots, v_n$ be a straight 
enumeration of $U(X)$. Suppose that $v_2$ is the only cut-vertex and it is 
non-dividing. If $n=8$, then $X$ or its dual is one of the graphs in
Figure~\ref{1non-dividing8}.
	\begin{figure}[H]
		\captionsetup[subfigure]{labelformat=empty}
		\begin{subfigure}[b]{.5\textwidth}
			\centering
			\begin{tikzpicture}
			\node[unlabelled] (a1) at (0,0)	{};
			\node[unlabelled] (x1) at (1,0)	{};
			\node[unlabelled] (b1) at (2,0)	{};
			\node[unlabelled] (b2) at (3,0)	{};
			\node[unlabelled] (b3) at (4,0)	{};
			\node[unlabelled] (b4) at (5,0)	{};
			\node[unlabelled] (b5) at (6,0)	{};
			\node[unlabelled] (b6) at (7,0)	{};
			\draw	(a1) -- (x1)	(x1) -- (b1)	
			(b1) -- (b2)	(b3) -- (b4)
			(b5) -- (b6)
			(x1) edge[bend left = 45] (b2)
			(b1) edge[bend left = 45] (b3)
			(b2) edge[bend left = 45] (b4)
			(b3) edge[bend left = 45] (b5)
			(b4) edge[bend left = 45] (b6)
			(x1) edge[bend left = 55] (b3)
			(b1) edge[bend left = 55] (b4)
			(b2) edge[bend left = 55] (b5)
			(b3) edge[bend left = 55] (b6);
			\draw[->-=.5]	(b2) to (b3);
			\draw[->-=.5]	(b5) to (b4);
			\end{tikzpicture}
			\subcaption{(i)}
		\end{subfigure}%
		\begin{subfigure}[b]{.5\textwidth}
			\centering
			\begin{tikzpicture}
			\node[unlabelled] (a1) at (0,0)	{};
			\node[unlabelled] (x1) at (1,0)	{};
			\node[unlabelled] (b1) at (2,0)	{};
			\node[unlabelled] (b2) at (3,0)	{};
			\node[unlabelled] (b3) at (4,0)	{};
			\node[unlabelled] (b4) at (5,0)	{};
			\node[unlabelled] (b5) at (6,0)	{};
			\node[unlabelled] (b6) at (7,0)	{};
			\draw	(a1) -- (x1)	(x1) -- (b1)	(b2) -- (b3)
			(b1) -- (b2)	(b3) -- (b4)	(b4) -- (b5)
			(b5) -- (b6)
			(x1) edge[bend left = 45] (b2)
			(b3) edge[bend left = 45] (b5)
			(b4) edge[bend left = 45] (b6)
			(x1) edge[bend left = 55] (b3)
			(b1) edge[bend left = 55] (b4)
			(b2) edge[bend left = 55] (b5)
			(x1) edge[bend left = 65] (b4);
			\draw[->-=.5]	(b1) to [bend left = 45] (b3);
			\draw[->-=.5]	(b4) to [bend right = 45] (b2);
			\end{tikzpicture}
			\subcaption{(ii)}
		\end{subfigure}
\caption{Obstructions with a unique non-dividing cut-vertex on 8 vertices.
\label{1non-dividing8}}
	\end{figure}
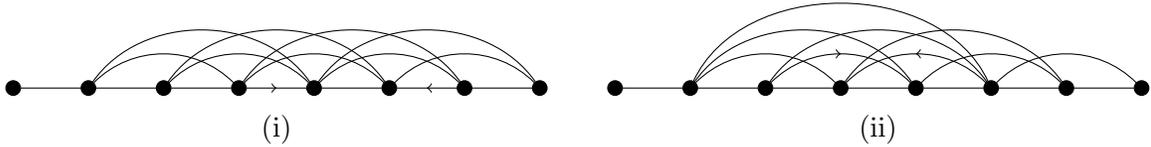
\end{theorem}
\pf Since there are six vertices succeeding $v_2$, exactly two of them are 
arc-balancing vertices and the other four are incident with arcs by 
Lemma \ref{only 1 cut-vertex lemma}(a). By Lemma \ref{only 1 cut-vertex lemma}(b), 
there exists a vertex succeeding $v_2$ that is adjacent to every vertex except for 
$v_1$. If any of $v_3, v_4, v_7, v_8$ is adjacent to every vertex except for $v_1$,
then $U(X)$ contains a copy of $K_5$ among the vertices $v_i$ with $i \geq 3$, 
contradicting Lemma \ref{at least 7}. Hence, only $v_5$ and $v_6$ can be adjacent 
to every vertex except for $v_1$.

Suppose $v_5$ is adjacent to every vertex except for $v_1$. By Lemma 
\ref{only 1 cut-vertex lemma}(b) again, $v_5$ is incident with exactly one arc. 
By Lemma \ref{pigbalance}, $v_8$ balances an arc between $v_5$ and one of 
$v_3, v_4$. If $v_3$ is an endvertex of this arc, then $v_3v_7$ would be an edge 
in $U(X)$, contradicting Lemma \ref{at least 7}. Hence, $v_8$ balances an arc 
between $v_5$ and $v_4$. It follows that $v_4v_7$ is an edge of $U(X)$. Moreover, 
there is an arc with both endvertices and arc-balancing vertex among 
$v_3, v_6, v_7$. If $v_7$ balances an arc between $v_3$ and $v_6$, then $v_3v_8$ 
is an edge of $U(X)$, contradiction Lemma \ref{only 1 cut-vertex lemma}(b). Hence 
$v_3$ balances an arc between $v_6$ and $v_7$. It follows that $X$ or its dual is 
Figure~\ref{1non-dividing8}(i).

On the other hand, suppose $v_5$ is not adjacent to every vertex except for $v_1$. 
By the previous discussion, $v_6$ must be the unique such vertex. 
By Lemma \ref{only 1 cut-vertex lemma}(b), $v_6$ is incident with an arc. 
By Lemma \ref{pigbalance}, $v_8$ balances an arc between $v_6$ and one of 
$v_3, v_4, v_5$. If $v_3$ is an endvertex of this arc, then $v_3v_7$ is an edge of 
$U(X)$, contradicting Lemma \ref{at least 7}. Suppose $v_8$ balances an arc 
between $v_6$ and $v_4$. Then, $v_4v_7$ is an edge of $U(X)$. Moreover, 
$\{v_3, v_5, v_7\}$ is an arc-balancing triple. By Lemma \ref{pigbalance}, 
$v_7$ balances an arc between $v_3$ and $v_5$. If $v_5v_8$ is an edge of $U(X)$, 
then $v_3v_8$ is also an edge, contradicting Lemma \ref{at least 7}. Hence $v_5v_8$
is not an edge and so $X$ or its dual is Figure~\ref{1non-dividing8}(ii). 
Suppose instead that $v_8$ balances an arc between $v_6$ and $v_5$. In this case, 
$\{v_3, v_4, v_7\}$ is an arc-balancing triple. By Lemma \ref{pigbalance}, $v_7$ 
balances an arc between $v_3$ and $v_4$. Since $v_5v_8$ and $v_3v_7$ are not edges 
of $U(X)$, $X$ or its dual is Figure~\ref{1non-dividing8}(ii).
\qed

\section{Obstructions without cut-vertices}
	
We now examine obstructions $X$ that do not contain cut-vertices. We shall consider
the complements $\overline{U(X)}$ of the underlying graphs $U(X)$. 
All theorems and proofs, including the drawings of obstructions $X$, in this 
section will be presented in terms of $\overline{U(X)}$ instead of $U(X)$.  

\begin{lemma}\label{2 non-adj}
Suppose that $X$ is an obstruction that contains no cut-vertices. Then, in
$\overline{U(X)}$, each vertex has at least two non-neighbours.
\end{lemma}
\pf Note that $X$ has at least three vertices. Since $X$ has no cut-vertices and 
$U(X)$ is connected, in $U(X)$ each vertex has at least two neighbours and hence 
in $\overline{U(X)}$ each vertex has at least two non-neighbours.
\qed

Recall from Corollary \ref{at most 6 non-cut-vertices} that if an obstruction $X$
has no cut-vertices then $\overline{U(X)}$ has at most six non-cut-vertices. 
We show this holds for every connected subgraph of $\overline{U(X)}$.

\begin{lemma}\label{non-cut-vertices}
Let $X$ be an obstruction that contains no cut-vertices and $H$ be a connected 
subgraph of $\overline{U(X)}$. Then $\overline{U(X)}$ contains at least as many 
non-cut-vertices as $H$. In particular, $H$ has at most six non-cut-vertices.
\end{lemma}
\pf Since adding edges does not decrease the number of non-cut-vertices, we may 
assume $H$ is an induced subgraph of $\overline{U(X)}$. Thus $H$ can
be obtained from $\overline{U(X)}$ by successively deleting non-cut-vertices. 
Since each deletion of a non-cut-vertex does not increase the number of 
non-cut-vertices, $\overline{U(X)}$ contains at least as many non-cut-vertices as 
$H$. By Corollary \ref{at most 6 non-cut-vertices}, $\overline{U(X)}$ has at most 
six non-cut-vertices. So $H$ has at most six non-cut-vertices.
\qed

\begin{lemma}\label{largest size of cycles}
If $X$ is an obstruction that contains no cut-vertices, then $\overline{U(X)}$ 
contains no induced cycle of length $\geq 6$.
\end{lemma}
\pf By Lemma \ref{non-cut-vertices}, any connected subgraph of $\overline{U(X)}$
has at most six non-cut-vertices. Thus $\overline{U(X)}$ contains no induced cycle 
of length $\geq 7$. Theorem \ref{Tucker} ensures that $\overline{U(X)}$ does not
contain an induced cycle of length 6. Therefore $\overline{U(X)}$ contains no 
induced cycle of length $\geq 6$. 
\qed

Lemma \ref{largest size of cycles} implies that any induced cycle in 
$\overline{U(X)}$ has length 3, 4 or 5. We show that $\overline{U(X)}$ contains 
at most one $C_3$ and at most one induced $C_5$ and moreover, if $\overline{U(X)}$ 
contains an induced $C_5$, then it does not contain an induced $C_3$ or $C_4$. 

\begin{lemma}\label{odd cycles}
Let $X$ be an obstruction. Suppose $C$ is an odd cycle (not necessarily induced) 
in $\overline{U(X)}$. Then, in $\overline{U(X)}$, each vertex is either in $C$ or 
adjacent to a vertex of $C$. In particular, each cut-vertex of $\overline{U(X)}$ 
is in $C$.
\end{lemma}
\pf Since $C$ is an odd cycle, $\overline{U(X)}$ contains an induced odd cycle
$C_{2k+1}$ on some vertices in $C$. By Theorem \ref{Tucker}, $\overline{U(X)}$ does
not contain $C_{2k+1} + K_1$ as an induced subgraph. Thus each vertex is either in 
$C_{2k+1}$ or adjacent to a vertex of $C_{2k+1}$. Since the vertices of $C_{2k+1}$ 
are all in $C$, each vertex is either in $C$ or adjacent to a vertex of $C$.  
Consequently, each cut-vertex of $\overline{U(X)}$ is in $C$. 
\qed

\begin{lemma}\label{onec3}
Suppose that $X$ is an obstruction that contains no cut-vertices. Then 
$\overline{U(X)}$ contains at most one $C_3$.
\end{lemma}
\pf Suppose that $C$ and $C'$ are two copies of $C_3$ in $\overline{U(X)}$. 
If $C$ and $C'$ share no common vertex, then every vertex of $\overline{U(X)}$ is 
either not in $C$ or not in $C'$ and hence by by Lemma \ref{odd cycles} is 
a non-cut-vertex. But $\overline{U(X)}$ has at most six non-cut-vertices by 
Corollary \ref{at most 6 non-cut-vertices}, so $\overline{U(X)}$ is a union of 
$C$ and $C'$. According to Proposition \ref{vertex classification} each vertex 
of $\overline{U(X)}$ is an endvertex of an arc or an arc-balancing vertex. 
There are at most four endvertices of arcs and at most two arc-balancing vertices. 
So, among the six vertices of $\overline{U(X)}$, four are the endvertices of arcs 
and the remaining two are arc-balancing vertices. 
Suppose that $(a,b)$ is an arc (of $X$) and $u$ is its balancing vertex that is 
adjacent to $a$ but not to $b$ in $\overline{U(X)}$. Then each of the remaining
three vertices is adjacent to $a$ or $b$ and thus to both $a, b$. Hence $b$ is 
the only non-neighbour of $a$ in $\overline{U(X)}$, a contradiction to Lemma 
\ref{2 non-adj}. Therefore any two copies of $C_3$ in $\overline{U(X)}$ must share 
a common vertex.

Suppose that $C$ and $C'$ share exactly one common vertex. Denote $C: v_1v_2v_3$ 
and $C': v_1v_4v_5$. Let $u, w$ be two non-neighbours of $v_1$ in $\overline{U(X)}$
guaranteed by Lemma \ref{2 non-adj}. Each vertex except $v_1$ is not in $C$ or 
$C'$ and hence by Lemma \ref{odd cycles} is a non-cut-vertex. 
Since $\overline{U(X)}$ has at most six non-cut-vertices, it consists of $C,C'$ and
$u,v$. A similar argument as above shows that, among the six non-cut-vertices 
$u, w, v_2, v_3, v_4, v_5$, four are the endvertices of arcs and the remaining two 
are arc-balancing vertices. We claim that the two arc-balancing vertices are $u,w$. 
Indeed, since $v_1$ is not an arc-balancing vertex, there is no arc between 
$u,w$ and $v_2, v_3, v_4, v_5$. Suppose that there is an arc between $u$ and $w$.
Assume without loss of generality that this arc is balanced by $v_2$ which is 
adjacent to $u$ but not $w$. By Lemma \ref{odd cycles}, $w$ is adjacent to a 
vertex in $C$. Since $w$ is not adjacent to $v_1$ or $v_2$, it is adjacent to 
$v_3$. Since $v_3$ does not balance the arc between $u$ and $w$, $v_3$ is adjacent
to $u$. But then $uv_2v_3$ and $C'$ are vertex-disjoint copies of $C_3$, 
a contradiction. Hence neither of $u, w$ is an endvertex of an arc so both
are arc-balancing vertices. Without loss of generality assume that $u$ balances
an arc between $v_2$ and $v_4$ and is adjacent to $v_2$ but not $v_4$. 
Since $v_3$ is adjacent to $v_2$, it must be adjacent to $v_4$. Similarly, $v_5$ 
must be adjacent to $v_2$. By Lemma \ref{odd cycles}, $u$ is adjacent to a vertex
in $C'$ which can only be $v_5$. Hence $uv_2v_5$ and $v_1v_3v_4$ are vertex-disjoint
copies of $C_3$, a contradiction. Therefore any two copies of $C_3$ in 
$\overline{U(X)}$ must share at least two common vertices. 

Suppose that $C$ and $C'$ share exactly two vertices. Denote $C: v_1v_2v_3$ and 
$C': v_1v_2v_4$. We claim that in $\overline{U(X)}$ any vertex 
$v \notin C \cup C'$ that is adjacent to one of $v_3, v_4$ must be adjacent to both 
$v_3, v_4$ and neither of $v_1, v_2$. Without loss of generality, suppose 
$v \notin C \cup C'$ is adjacent to $v_3$. If it is also adjacent to $v_1$, then 
$v_1v_3v$ and $v_1v_2v_4$ would be two distinct copies of $C_3$ in 
$\overline{U(X)}$ that share exactly one common vertex, a contradiction to the 
above. Hence, $v$ is not adjacent to $v_1$. Similarly, $v$ is not adjacent to $v_2$.
By Lemma \ref{odd cycles} $v$ must be adjacent to a vertex in $C'$ so it is 
adjacent to $v_4$. 

Now, we show that $v_3$ and $v_4$ are incident with different arcs. Suppose $v_3$ 
balances an arc between $a$ and $b$ and is adjacent to $a$ but not $b$. 
If $a=v_1$, then since $v_2$ and $v_4$ are adjacent to $v_1$, they must also be 
adjacent to $b$. So, $bv_2v_4$ and $C$ are two distinct copies of $C_3$ in 
$\overline{U(X)}$ that share exactly one common vertex, a contradiction. Thus, 
$a \neq v_1$. Similarly, $a \neq v_2$. Suppose $a = v_4$. Since $v_1$ and $v_2$ are adjacent to $a$, they must be adjacent to $b$ as well. Thus $v_1v_2b$ and 
$v_3v_4v_2$ are two copies of $C_3$ in $\overline{U(X)}$ that share exactly one 
common vertex, a contradiction. 

It follows that $a \notin C \cup C'$. Since $a \notin C\cup C'$ and $a$ is adjacent
to $v_3$, it is adjacent to both of $v_3, v_4$ by the above claim. Since $v_4$ is 
adjacent to $a$, it must also be adjacent to $b$. Moreover, $b \notin C \cup C'$ 
because it is adjacent to $v_4$ but not $a$. Since $b \notin C \cup C'$ and $b$ is 
adjacent to $v_4$, the above claim implies $b$ is adjacent to both of $v_3, v_4$, 
a contradiction because $v_3$ balances the arc $(a,b)$. Thus, $v_3$ is not an 
arc-balancing vertex. Since $v_3$ is not in $C'$, it is not a cut-vertex of 
$\overline{U(X)}$ by Lemma \ref{odd cycles}. By Proposition 
\ref{vertex classification}, $v_3$ is incident with an arc. Similarly, $v_4$ is 
incident with an arc. If $v_3$ and $v_4$ are incident with the same arc, then 
there must be a vertex $u$ that is adjacent to exactly one of $v_3, v_4$ because 
arcs in $X$ are not balanced. Clearly, $u \notin C \cup C'$. This is a 
contradiction because any vertex not in $C \cup C'$ is adjacent to either both of 
$v_3, v_4$ or neither, by above claim. Thus, $v_3$ and $v_4$ are each incident 
with a different arc.

Suppose there exists a vertex $v \notin C \cup C'$ that is adjacent to either of 
$v_3, v_4$. By the above claim, we know that $v$ is adjacent to both of $v_3, v_4$ 
and neither of $v_1, v_2$. Since $v$ is adjacent to both of $v_3, v_4$, which are 
each incident with a different arc, $v$ is not incident with an arc. Moreover, 
since $v$ is not on the odd cycle $C$, Lemma \ref{odd cycles} implies $v$ is not 
a cut-vertex of $\overline{U(X)}$. So by Proposition \ref{vertex classification}, 
$v$ is an arc-balancing vertex. Since $v_3$ and $v_4$ are each incident with a 
different arc, we may assume without loss of generality that $v$ balances an arc 
incident with $v_3$. Let $w$ denote the other endvertex of this arc. Then, $w$ is 
adjacent to $v_1$ and $v_2$, so $v_1v_2w$ is a triangle. 
By Lemma \ref{odd cycles}, $v$ is adjacent to a vertex in $v_1v_2w$, which must be 
$w$, a contradiction because $v$ balances the arc between $v_3$ and $w$. 
It follows $\overline{U(X)}$ does not contain a vertex $v \notin C \cup C'$ 
that is adjacent to either of $v_3, v_4$.

By Lemma \ref{2 non-adj}, $v_1$ has at least two non-neighbours, say $u$ and $w$. 
Clearly, $u, w \notin C \cup C'$. So by the above, neither of $u, w$ is adjacent to
either of $v_3, v_4$. By Lemma \ref{odd cycles}, each of $u, w$ is adjacent to 
a vertex in $C$, which must be $v_2$. Similarly, $v_2$ has at least two 
non-neighbours, say $x$ and $y$, and each is adjacent to $v_1$. 
By Lemma \ref{odd cycles}, each of $v_3, v_4, u, w, x, y$ is a non-cut-vertex of 
$\overline{U(X)}$ and so by Corollary \ref{at most 6 non-cut-vertices} they form
two disjoint arc-balancing triples. Since $v_3$ and $v_4$ are each incident with 
a different arc,
exactly two of $u, w, x, y$ are arc-balancing vertices. Without loss of generality,
assume $u$ is an arc-balancing vertex for an arc incident with $v_3$. Since $v_3$ 
is adjacent to both $v_1$ and $v_2$, the other endvertex must also be adjacent to 
both $v_1$ and $v_2$. This is a contradiction because none of $w, x, y$ is adjacent
to both $v_1$ and $v_2$ by assumption. It follows that $C$ and $C'$ cannot share 
two common vertices. Therefore $\overline{U(X)}$ contains at most one $C_3$.
\qed

\begin{lemma}\label{one C5}
Suppose that $X$ is an obstruction that contains no cut-vertices. Then 
$\overline{U(X)}$ contains at most one induced $C_5$.
\end{lemma}
\pf Suppose that $C$ and $C'$ are induced copies of $C_5$ contained in 
$\overline{U(X)}$. By Lemma \ref{odd cycles}, any vertex not in $C$ or $C'$ is
a non-cut-vertex of $\overline{U(X)}$ and hence by Corollary 
\ref{at most 6 non-cut-vertices} there can be at most six such vertices. Thus $C$ 
and $C'$ must share at least two common vertices. If $C$ and $C'$ share less 
two or three common vertices, then the subgraph of $\overline{U(X)}$ induced by 
$C \cup C'$ is connected and has at least seven non-cut-vertices, contradicting 
Lemma \ref{non-cut-vertices}. Hence, $C$ and $C'$ must share exactly four vertices.

Denote $C: v_1v_2v_3v_4v_5$ and $C': v_2v_3v_4v_5v_6$. Then $v_1v_6$ is not an edge
in $\overline{U(X)}$ as otherwise $v_1v_2v_6$ and $v_1v_5v_6$ are two copies of 
$C_3$ in $\overline{U(X)}$, a contradiction to Lemma \ref{onec3}. 
We claim that $v_1, v_6$ are endvertices of arcs in $X$. By symmetry we only prove
that $v_1$ is an endvertex of an arc in $X$. We prove it by contradiction. 
So assume that $v_1$ is not an endvertex of an arc in $X$. Since $v_1$ is not in
$C'$, by Lemma \ref{odd cycles} it is not a cut-vertex of $\overline{U(X)}$.
Hence $v_1$ is an arc-balancing vertex for some arc according to Proposition 
\ref{vertex classification}. Suppose that $v_1$ balances the arc between vertices
$a, b$ and is adjacent to $a$ but not to $b$ in $\overline{U(X)}$. 
If $a$ is not in $C'$, then $a$ must be adjacent to a vertex of $C'$ by
Lemma \ref{odd cycles}. But then the subgraph of $\overline{U(X)}$ induced by 
$C \cup C' \cup \{a\}$ is connected and has seven non-cut-vertices, contradicting 
Lemma \ref{non-cut-vertices}. Hence $a$ is a vertex of $C'$ and therefore it is 
$v_2$ or $v_5$. Assume by symmetry $a = v_2$. Since $v_1$ balances the arc between
$a, b$, every vertex not in $\{v_1,a,b\}$ is either adjacent to both $a, b$ or 
neither. It follows that $b$ cannot be in $C \cup C'$. Thus the subgraph of 
$\overline{U(X)}$ induced by $C \cup C' \cup \{b\}$ is connected and has seven
non-cut-vertices, a contradiction to Lemma \ref{non-cut-vertices}. 
Therefore $v_1, v_6$ are both endvertices of arcs of $X$. 
We claim that there is no arc between $v_1, v_6$. Suppose not; there is an arc
between $v_1, v_6$. Then there must exist a vertex $u$ adjacent to exactly one of 
$v_1, v_6$. 
A similar argument as above shows that $u$ is not in $C \cup C'$ but adjacent to
a vertex in $C \cup C'$. Thus the subgraph of $\overline{U(X)}$ induced by 
$C \cup C' \cup \{u\}$ is connected and contains seven non-cut-vertices, 
a contradiction. Thus, $v_1, v_6$ are endvertices of different arcs.

The subgraph of $\overline{U(X)}$ induced by $C \cup C'$ contains six 
non-cut-vertices, so $\overline{U(X)}$ contains six non-cut-vertices by 
Lemma \ref{non-cut-vertices}. It follows from Proposition 
\ref{vertex classification} that $X$ contains exactly four vertices incident 
to arcs and exactly two arc-balancing vertices. In particular, both arcs have 
an arc-balancing vertex.

By Proposition \ref{vertex classification}, $v_3$ is a cut-vertex of 
$\overline{U(X)}$, an arc-balancing vertex, or is incident with an arc. 
We claim it must be a cut-vertex of $\overline{U(X)}$. Suppose instead $v_3$ is an 
arc-balancing vertex. Without loss of generality, assume it balances the arc 
incident with $v_1$. Then, the other endvertex must be adjacent to each of 
$v_2, v_5$. Clearly, the subgraph of $\overline{U(X)}$ induced by $C \cup C'$ 
together with this endvertex contains seven non-cut-vertices, a contradiction. 
On the other hand, suppose $v_3$ is incident with an arc. The other endvertex is 
one of $v_1, v_6$. Without loss of generality, assume it is $v_1$. Then, $v_4$ and 
$v_5$ are both vertices adjacent to exactly one of the endvertices of this arc, 
so the arc between $v_1$ and $v_3$ has no arc-balancing vertex, a contradiction. 
Thus, $v_3$ is a cut-vertex of $\overline{U(X)}$.

Let $v_7$ be a neighbour of $v_3$ belong to a different component of 
$\overline{U(X-v_3)}$ as the vertices in $(C \cup C') \setminus \{v_3\}$. 
By Lemma \ref{odd cycles}, $v_7$ cannot be a cut-vertex of $\overline{U(X)}$. 
On the other hand, suppose $v_7$ is incident with an arc. 
Without loss of generality, assume the other endvertex is $v_1$. Since $v_2$ and 
$v_3$ are both vertices adjacent to exactly one of $v_1, v_7$, there is 
no corresponding arc-balancing vertex for this arc, a contradiction. 
Thus, $v_7$ cannot be incident with an arc. 
By Proposition \ref{vertex classification}, $v_7$ is an arc-balancing vertex for 
either the arc incident with $v_1$ or the arc incident with $v_6$. 
In either case, the other endvertex must be adjacent to both $v_2$ and $v_5$. 
Clearly, the subgraph of $\overline{U(X)}$ induced by $C \cup C'$ together with 
this endvertex contains seven non-cut-vertices, a contradiction.
\qed

\begin{lemma}\label{C5lemma}
Let $X$ be an obstruction that contains no cut-vertices. If $\overline{U(X)}$ 
contains an induced $C_5$, then it contains neither $C_3$ nor induced $C_4$.
\end{lemma}
\pf Let $C: v_1v_2v_3v_4v_5$ be an induced $C_5$ in $\overline{U(X)}$. We first 
show that $\overline{U(X)}$ does not contain $C_3$. Suppose otherwise and let 
$C'$ be a $C_3$ in $\overline{U(X)}$. A similar argument as the one in 
Lemma \ref{one C5} shows that $C$ and $C'$ have exactly two common vertices. 
Without loss of generality let $C': v_1v_2v_6$. By Lemma \ref{odd cycles}, $v_4$ 
must be adjacent to a vertex in $C'$, which clearly must be $v_6$. 
The subgraph induced by $C \cup C'$ contains six non-cut-vertices, so 
$\overline{U(X)}$ contains six non-cut-vertices by Lemma \ref{non-cut-vertices}. 
Each of these six non-cut-vertices is an arc-balancing vertex or incident with 
an arc by Proposition \ref{vertex classification}. Hence, each arc has 
an arc-balancing vertex. If both endvertices of some arc are in $C$, then $C$ 
contains two other vertices which are both adjacent to exactly one of 
the endvertices, contradicting the fact that each arc has a unique arc-balancing 
vertex. It follows that each arc has at most one endvertex in $C$. In particular, 
at most two vertices in $C$ are incident with arcs. On the other hand, at most two 
vertices in $C$ are arc-balancing. It follows from 
Proposition \ref{vertex classification} that $C$ has a cut-vertex. 
By Lemma \ref{odd cycles}, each cut-vertex belongs to $C \cap C'$, so only 
$v_1$ and $v_2$ can be cut-vertices.

We claim that if $v_1$ is a cut-vertex, then there exists a vertex $v_7$ that is 
adjacent only to $v_1$ and an arc between $v_5$ and $v_7$ that is balanced by $v_4$.
Suppose $v_1$ is a cut-vertex. Let $v_7$ be a vertex adjacent to $v_1$ that belongs
to a different component of $\overline{U(X-v_1)}$ as the vertices in 
$(C \cup C') \setminus \{v_1\}$. If $v_7$ is adjacent to a vertex other than 
$v_1$, then that vertex must be adjacent to a vertex in $C \cup C'$ by 
Lemma \ref{odd cycles}, contradicting the choice of $v_7$. Hence, $v_7$ is adjacent
only to $v_1$. By Lemma \ref{odd cycles}, $v_7$ is not a cut-vertex. 
Suppose $v_7$ is an arc-balancing vertex. Then, $v_1$ is incident with an arc, and 
the other endvertex of this arc must be adjacent to $v_2, v_5, v_6$. 
Clearly, this endvertex is none of the vertices in $C \cup C'$, so the subgraph of 
$\overline{U(X)}$ induced by $C \cup C'$ together with this endvertex contains 
seven non-cut-vertices, contradicting Lemma \ref{non-cut-vertices}. 
Hence, $v_7$ is not an arc-balancing vertex. 
By Proposition \ref{vertex classification}, $v_7$ is incident with an arc. 
Since $v_7$ is a leaf, the other endvertex $u$ of the arc incident with $v_7$ has
exactly two neighbours, one of which is $v_1$. Since any subgraph of 
$\overline{U(X)}$ contains at most six non-cut-vertices, $u$ must be in
$C \cup C'$. Clearly, $u = v_5$ and the vertex which balances the arc between 
$v_5$ and $v_7$ is $v_4$. This proves our claim.

Recall that at least one of $v_1, v_2$ is a cut-vertex. Without loss of generality,
assume $v_1$ is a cut-vertex. By the above, there exist a vertex $v_7$ that is 
adjacent only to $v_1$ and an arc between $v_5$ and $v_7$ that is balanced by $v_4$.
In particular, $v_6$ is not an arc-balancing vertex for this arc. 
By Lemma \ref{odd cycles}, $v_6$ is not a cut-vertex. Hence 
by Proposition \ref{vertex classification}, $v_6$ is arc-balancing for 
the other arc or incident with it. By symmetry, if $v_2$ is also a cut-vertex, 
then there exist a vertex $v_8$ that is adjacent only to $v_2$ and an arc between 
$v_3$ and $v_8$, a contradiction because $v_6$ is neither 
arc-balancing for this arc nor incident with this arc. Thus, $v_2$ is not 
a cut-vertex. Since $v_4$ balances an arc between $v_5$ and $v_7$ and 
$v_2, v_3, v_6$ are non-cut-vertices, $\{v_2, v_3, v_6\}$ is an arc-balancing 
triple. It follows that $v_2$ balances an arc between $v_3$ and $v_6$, 
a contradiction. Thus, $\overline{U(X)}$ does not contain an induced $C_3$.

It remains to show that $\overline{U(X)}$ does not contain an induced $C_4$. 
Suppose otherwise, and let $C'$ be such a cycle. A similar argument as above shows 
that $C$ and $C'$ have exactly three common vertices. 
Let $C': v_1v_2v_3v_6$. If $v_6$ is adjacent to neither $v_4$ or $v_5$, then 
$v_1v_6v_3v_4v_5$ is an induced $C_5$, contradicting Lemma \ref{one C5}. Hence, 
$v_6$ is adjacent to one of $v_4, v_5$. It follows that $\overline{U(X)}$ 
contains an induced $C_3$, a contradiction.
\qed

\subsection{$\overline{U(X)}$ is disconnected}

We first examine obstructions $X$ that do not contain cut-vertices for which
$\overline{U(X)}$ is disconnected. These obstructions have a simple structure
as described in the following theorem. 

\begin{theorem}\label{disconnected theorem}
Let $X$ be an obstruction that does not contain cut-vertices. Suppose that 
$\overline{U(X)}$ is disconnected. Then the following statements hold:
\begin{itemize}
\item $\overline{U(X)}$ is the union of two disjoint paths 
      $P: p_1, \dots, p_k$ and $Q: q_1, \dots, q_\ell$;
\item $X$ or its dual contains the arcs $(p_1,q_1), (q_\ell,p_k)$ if $k+\ell$ is 
      even, and $(p_1,q_1), (p_k,q_\ell)$ otherwise.
\end{itemize}
That is, $U(X)$ is one of the graphs in Figure~\ref{disconnectedfigure} and $X$ or 
its dual contains the dotted arcs. 
\begin{figure}[H]
        \captionsetup[subfigure]{labelformat=empty}
        \begin{subfigure}[b]{.5\textwidth}
        	\centering
        	\begin{tikzpicture}
        	\node[unlabelled] (p1) at (0,1) {};
        	\node[unlabelled] (p2) at (1,1) {};
        	\node (p3) at (2,1) {$\dots$};
        	\node[unlabelled] (p4) at (3,1) {};
        	\node[unlabelled] (q1) at (0,0) {};
        	\node[unlabelled] (q2) at (1,0) {};
        	\node (q3) at (2,0) {$\dots$};
        	\node[unlabelled] (q4) at (3,0) {};
        	\node (lp1) at (0,1.3) {$p_1$};
        	\node (lp2) at (1,1.3) {$p_2$};
        	\node (lp4) at (3,1.3) {$p_k$};
        	\node (lq1) at (0,-0.3) {$q_1$};
        	\node (lq2) at (1,-0.3) {$q_2$};
        	\node (lq4) at (3,-0.3) {$q_\ell$};
        	\draw (p1) -- (p2) (p2) -- (p3) (p3) -- (p4) (q1) -- (q2) (q2) -- (q3) (q3) -- (q4);
        	\draw[dashed, ->-=0.5] (p1) to (q1);
        	\draw[dashed, ->-=0.5] (q4) to (p4);
        	\end{tikzpicture}
        	\subcaption{(i) $k + \ell$ is even}
        \end{subfigure}%
    	\begin{subfigure}[b]{.5\textwidth}
    		\centering
    		\begin{tikzpicture}
    		\node[unlabelled] (p1) at (0,1) {};
    		\node[unlabelled] (p2) at (1,1) {};
    		\node (p3) at (2,1) {$\dots$};
    		\node[unlabelled] (p4) at (3,1) {};
    		\node[unlabelled] (q1) at (0,0) {};
    		\node[unlabelled] (q2) at (1,0) {};
    		\node (q3) at (2,0) {$\dots$};
    		\node[unlabelled] (q4) at (3,0) {};
    		\node (lp1) at (0,1.3) {$p_1$};
    		\node (lp2) at (1,1.3) {$p_2$};
    		\node (lp4) at (3,1.3) {$p_k$};
    		\node (lq1) at (0,-0.3) {$q_1$};
    		\node (lq2) at (1,-0.3) {$q_2$};
    		\node (lq4) at (3,-0.3) {$q_\ell$};
    		\draw (p1) -- (p2) (p2) -- (p3) (p3) -- (p4) (q1) -- (q2) (q2) -- (q3) (q3) -- (q4);
    		\draw[dashed, ->-=0.5] (p1) to (q1);
    		\draw[dashed, ->-=0.5] (p4) to (q4);
    		\end{tikzpicture}
    		\subcaption{(ii) $k + \ell$ is odd}
    	\end{subfigure}
\caption{Obstructions $X$ for which $\overline{U(X)}$ is disconnected.
        \label{disconnectedfigure}}
\end{figure}
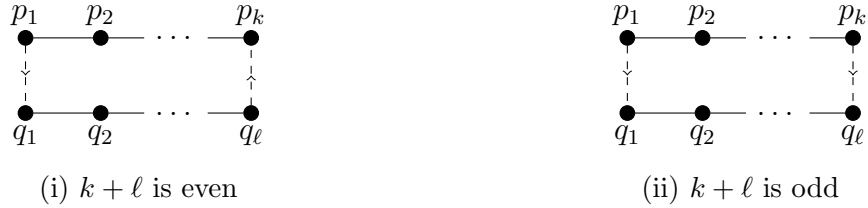
\end{theorem}
\pf Let $(a,b)$ and $(c,d)$ be the two arcs of $X$. Then $ab$ and $cd$ belong to
the same implication class of $U(X)$. By Theorem \ref{structure}, either $ab$ and
$cd$ are unbalanced edges of $U(X)$ within a component of $\overline{U(X)}$ or they
are edges between two components of $\overline{U(X)}$. Since $\overline{U(X)}$ is
disconnected, it has at least two components. If some component of 
$\overline{U(X)}$ does not contain any of $a, b, c, d$, then any non-cut-vertex of
that component is not a cut-vertex of $U(X)$ by assumption and is not an 
arc-balancing vertex because it is not adjacent to any of $a,b,c,d$ in 
$\overline{U(X)}$. This contradicts Proposition \ref{vertex classification}. 
Thus $\overline{U(X)}$ has exactly two components and $ab, cd$ are edges between 
them.

Consider a component $H$ of $\overline{U(X)}$ and let $P$ be a shortest path in 
$H$ between some two of $a, b, c, d$. If $H$ contains a vertex $v$ that is
not in $P$ then it follows from Proposition \ref{complement gamma} that $ab$ and 
$cd$ are still in the same implication class of $U(X-v)$, which is a contradiction 
to Theorem \ref{different implication classes theorem}. This shows that each 
component of $\overline{U(X)}$ is a path connecting two vertices of $a, b, c, d$ 
and $\overline{U(X)}$ is the union of two disjoint paths.

Let $P: p_1, \dots, p_k$ and $Q: q_1, \dots, q_\ell$ be the two paths in 
$\overline{U(X)}$. The two arcs are between $p_1$ and $q_1$ and between 
$p_k$ and $q_\ell$ respectively. Without loss of generality, assume $(p_1, q_1)$ is
an arc. Suppose $k+\ell$ is even. If $k, \ell$ are both even, then 
$(p_1,q_1)\Gamma^*(q_\ell,p_1)$ and $(q_\ell,p_1)\Gamma^*(p_k, q_\ell)$ by 
Proposition \ref{complement gamma}. Since the arcs must be opposing, the other arc 
is $(q_\ell, p_k)$. Otherwise, $k, \ell$ are both odd. In this case, we have 
$(p_1,q_1)\Gamma^*(p_1,q_\ell)$ and $(p_1,q_\ell)\Gamma^*(p_k, q_\ell)$, so 
the other arc is $(q_\ell, p_k)$. Hence $X$ or its dual is 
Figure~\ref{disconnectedfigure}(i). A similar proof shows that, when $k+\ell$ is 
odd, $X$ or its dual is Figure~\ref{disconnectedfigure}(ii).  
\qed

\begin{corollary} \label{disconnected}
If $X$ is an obstruction that does not contain cut-vertices and for which
$\overline{U(X)}$ is disconnected, then $\overline{U(X)}$ contains no cycle.
\qed
\end{corollary}

\subsection{$\overline{U(X)}$ is a tree}

We next examine obstructions $X$ that do not contain cut-vertices and for which
$\overline{U(X)}$ is a tree. We begin with a useful lemma.

\begin{lemma}\label{tree structure}
Let $X$ be an obstruction that contains no cut-vertices. If $\overline{U(X)}$ is 
a tree, then it is a caterpillar and has at most four leaves. Moreover, suppose 
$P: p_1, p_2, \dots, p_k$ is a longest path in $\overline{U(X)}$. 
If $p_1$ is an arc-balancing vertex, then $p_2$ has only two neighbours (namely, 
$p_1,p_3$) and $p_1$ balances an arc between $p_2$ and a leaf adjacent to $p_3$ 
but not in $P$.
\end{lemma}
\pf Since $U(X)$ is a proper circular-arc graph, $\overline{U(X)}$ does not contain
the fifth graph in Figure \ref{tuckerlist} by Theorem \ref{Tucker} and hence is 
a caterpillar. If $v$ is a leaf of $\overline{U(X)}$ that is not incident with 
an arc of $X$, then by Proposition \ref{vertex classification} $v$ is an 
arc-balancing vertex and hence adjacent to a vertex that is incident with an arc. 
Clearly, the vertex adjacent to $v$ cannot be adjacent to any other leaf.
Since there are at most four vertices incident with arcs, $\overline{U(X)}$ has 
at most four leaves. 

Since $P$ is a longest path, $p_1$ is a leaf. If $p_1$ is an arc-balancing vertex,
then $p_2$ is incident with an arc balanced by $p_1$. Let $u$ be the other endvertex
of the arc. Every vertex other than $p_1$ is adjacent either to both $p_2, u$ or 
neither. Since $p_3$ is adjacent to $p_2$, it is adjacent to $u$.  Since 
$\overline{U(X)}$ is a tree, $p_3$ is the only neighbour of $p_2$ other than 
$p_1$ and the only neighbour of $u$. It follows that $p_1,p_3$ are the only 
neighbours of $p_2$. If $u$ is in $P$ then $u = p_4$ and $k = 4$. Thus each vertex 
not in $P$ can only be adjacent to $p_3$ in $\overline{U(X)}$, which implies that 
$p_1$ is a cut-vertex of $U(X)$, a contradiction to the assumption. Therefore $u$ 
is a leaf of $\overline{U(X)}$ adjacent to $p_3$ but not in $P$.
\qed

\begin{theorem}\label{treecase}
Let $X$ be an obstruction that contains no cut-vertices and for which  
$\overline{U(X)}$ is a tree. Let $P: p_1, p_2, \dots, p_k$ be a longest path in
$\overline{U(X)}$. Then $\overline{U(X)}$ consists of $P$ and $u, v$ (possibly 
$u = v$) where $u$ is either a leaf adjacent to some $p_\ell$ but not in $P$ or 
$u=p_\ell$ and $v$ is either a leaf adjacent to some $p_j$ but not in $P$ or 
$v=p_j$, and one of the following statements holds:
\begin{description}
\item{(i)} $u$ is not in $P$ and $\ell = 3$, $v$ is not in $P$ and $j = k-2$, and 
         $X$ or its dual has arcs $(p_2,u), (p_{k-1},v)$ 
         (See Figure~\ref{treefigure}(i)); 
\item{(ii)} $u$ is not in $P$ and $\ell = 3$, $1 \leq j \leq k-2$ with $j> 2$ when 
         $v$ is not in $P$, and $X$ or its dual has arcs $(p_2,u), (v,p_k)$ if 
         either $k+j$ is even and $v$ is not in $P$ or $k+j$ is odd and $v$ is 
         in $P$; otherwise $X$ or its dual has arcs $(p_2,u), (p_k,v)$
         (See Figure~\ref{treefigure}(ii));

\begin{figure}[H]
	\centering
	\captionsetup[subfigure]{labelformat=empty}
	\begin{subfigure}[b]{.5\textwidth}
		\centering
		\begin{tikzpicture}
		\node[unlabelled] (p1) at (0,0) {};
		\node[unlabelled] (p2) at (1,0) {};
		\node[unlabelled] (p3) at (2,0) {};
		\node (p) at (3,0) {$\dots$};
		\node[unlabelled] (pk2) at (4,0) {};
		\node[unlabelled] (pk1) at (5,0) {};
		\node[unlabelled] (pk) at (6,0) {};
		\node[unlabelled] (u) at (2,1) {};
		\node[unlabelled] (v) at (4,1) {};
		\node (lp1) at (0,-0.3) {$p_1$};
		\node (lp2) at (1,-0.3) {$p_2$};
		\node (lp3) at (2,-0.3) {$p_3$};
		\node (lpk2) at (4,-0.3) {$p_{k-2}$};
		\node (lpk1) at (5,-0.3) {$p_{k-1}$};
		\node (lpk) at (6,-0.3) {$p_k$};
		\node (lu) at (2,1.3) {$u$};
		\node (lv) at (4,1.3) {$v$};
		\draw (p1) -- (p2) (p2) -- (p3) (p3) -- (p) (p) -- (pk2) (pk2) -- (pk1) (pk1) -- (pk) (p3) -- (u) (pk2) -- (v);
		\draw[dashed, ->-=0.5] (p2) to (u);
		\draw[dashed, ->-=0.5] (pk1) to (v);
		\end{tikzpicture}
		\subcaption{(i)}
	\end{subfigure}

	\begin{subfigure}[t]{.45\textwidth}
		\centering
		\begin{tikzpicture}
		\node[unlabelled] (p1) at (0,0) {};
		\node[unlabelled] (p2) at (1,0) {};
		\node[unlabelled] (p3) at (2,0) {};
		\node (p) at (2.75,0) {$\dots$};
		\node[unlabelled] (pj) at (3.5,0) {};
		\node (p') at (4.25,0) {$\dots$};
		\node[unlabelled] (pk1) at (5,0) {};
		\node[unlabelled] (pk) at (6,0) {};
		\node[unlabelled] (u) at (2,1) {};
		\node[unlabelled] (v) at (3.5,1) {};
		\node (lp1) at (0,-0.3) {$p_1$};
		\node (lp2) at (1,-0.3) {$p_2$};
		\node (lp3) at (2,-0.3) {$p_3$};
		\node (lpj) at (3.5,-0.3) {$p_j$};
		\node (lpk1) at (5,-0.3) {$p_{k-1}$};
		\node (lpk) at (6,-0.3) {$p_k$};
		\node (lu) at (2,1.3) {$u$};
		\node (lv) at (3.5,1.3) {$v$};
		\draw (p1) -- (p2) (p2) -- (p3) (p3) -- (p) (p) -- (pj) (pj) -- (p') (p') -- (pk1) (pk1) -- (pk) (p3) -- (u) (pj) -- (v);
		\draw[dashed, ->-=0.5] (p2) to (u);
		\end{tikzpicture}
\subcaption{(ii.a): The second arc is $(v, p_k)$ if $k+j$ is even and 
               $(p_k, v)$ otherwise where $2<j<k-1$.}
	\end{subfigure}\hfill%
	\begin{subfigure}[t]{.45\textwidth}
		\centering
		\begin{tikzpicture}
		\node[unlabelled] (p1) at (0,0) {};
		\node[unlabelled] (p2) at (1,0) {};
		\node[unlabelled] (p3) at (2,0) {};
		\node (p) at (2.75,0) {$\dots$};
		\node[unlabelled] (pk1) at (3.5,0) {};
		\node[unlabelled] (pk) at (4.5,0) {};
		\node[unlabelled] (u) at (2,1) {};
		\node (lp1) at (0,-0.3) {$p_1$};
		\node (lp2) at (1,-0.3) {$p_2$};
		\node (lp3) at (2,-0.3) {$p_3$};
		\node (lpk1) at (3.5,-0.3) {$p_{k-1}$};
		\node (lpk) at (4.5,-0.3) {$p_k$};
		\node (lu) at (2,1.3) {$u$};
		\draw (p1) -- (p2) (p2) -- (p3) (p3) -- (p) (p) -- (pk1) (pk1) -- (pk) (p3) -- (u);
		\draw[dashed, ->-=0.5] (p2) to (u);
		\end{tikzpicture}
\subcaption{(ii.b): The second arc is $(p_j,p_k)$ if $k+j$ is odd and 
               $(p_k,p_j)$ otherwise where $1 \leq j \leq k-2$.}
	\end{subfigure}
	
	\begin{subfigure}[t]{.45\textwidth}
		\centering
		\begin{tikzpicture}
		\node[unlabelled] (p1) at (0,0) {};
		\node[unlabelled] (p2) at (1,0) {};
		\node (p) at (1.75,0) {$\dots$};
		\node[unlabelled] (pl) at (2.5,0) {};
		\node[unlabelled] (pj) at (3.5,0) {};
		\node (p') at (4.25,0) {$\dots$};
		\node[unlabelled] (pk1) at (5,0) {};
		\node[unlabelled] (pk) at (6,0) {};
		\node[unlabelled] (u) at (2.5,1) {};
		\node[unlabelled] (v) at (3.5,1) {};
		\node (lp1) at (0,-0.3) {$p_1$};
		\node (lp2) at (1,-0.3) {$p_2$};
		\node (lpl) at (2.5,-0.3) {$p_\ell$};
		\node (lpj) at (3.5,-0.3) {$p_{\ell+1}$};
		\node (lpk1) at (5,-0.3) {$p_{k-1}$};
		\node (lpk) at (6,-0.3) {$p_k$};
		\node (lu) at (2.5,1.3) {$u$};
		\node (lv) at (3.5,1.3) {$v$};
		\draw (p1) -- (p2) (p2) -- (p) (p) -- (pl) (pl) -- (pj) (pj) -- (p') (p') -- (pk1) (pk1) -- (pk) (pl) -- (u) (pj) -- (v);
		\draw[dashed, ->-=0.5] (p1) to [bend left = 65] (pk);
		\end{tikzpicture}
		\subcaption{(iii.a): The second arc is $(v, u)$ if $k$ is even and $(u, v)$ otherwise.}
	\end{subfigure}\hfill%
	\begin{subfigure}[t]{.45\textwidth}
		\centering
		\begin{tikzpicture}
		\node[unlabelled] (p1) at (0,0) {};
		\node[unlabelled] (p2) at (1,0) {};
		\node (p) at (1.75,0) {$\dots$};
		\node[unlabelled] (pl) at (2.5,0) {};
		\node[unlabelled] (pj) at (3.5,0) {};
		\node (p') at (4.25,0) {$\dots$};
		\node[unlabelled] (pk1) at (5,0) {};
		\node[unlabelled] (pk) at (6,0) {};
		\node[unlabelled] (u) at (2.5,1) {};
		\node (lp1) at (0,-0.3) {$p_1$};
		\node (lp2) at (1,-0.3) {$p_2$};
		\node (lpl) at (2.5,-0.3) {$p_\ell$};
		\node (lpj) at (3.5,-0.3) {$p_{\ell+1}$};
		\node (lpk1) at (5,-0.3) {$p_{k-1}$};
		\node (lpk) at (6,-0.3) {$p_k$};
		\node (lu) at (2.5,1.3) {$u$};
		\draw (p1) -- (p2) (p2) -- (p) (p) -- (pl) (pl) -- (pj) (pj) -- (p') (p') -- (pk1) (pk1) -- (pk) (pl) -- (u);
		\draw[dashed, ->-=0.5] (p1) to [bend left = 65] (pk);
		\end{tikzpicture}
\subcaption{(iii.b): The second arc is $(p_{\ell+1},u)$ if $k$ is odd and 
      $(u,p_{\ell+1})$ otherwise.}
	\end{subfigure}
	
	\begin{subfigure}[b]{.45\textwidth}
		\centering
		\begin{tikzpicture}
		\node[unlabelled] (p1) at (0,0) {};
		\node[unlabelled] (p2) at (1,0) {};
		\node (p) at (1.75,0) {$\dots$};
		\node[unlabelled] (pj) at (2.5,0) {};
		\node[unlabelled] (pl) at (3.5,0) {};
		\node (p') at (4.25,0) {$\dots$};
		\node[unlabelled] (pk1) at (5,0) {};
		\node[unlabelled] (pk) at (6,0) {};
		\node[unlabelled] (v) at (2.5,1) {};
		\node[unlabelled] (u) at (3.5,1) {};
		\node (lp1) at (0,-0.3) {$p_1$};
		\node (lp2) at (1,-0.3) {$p_2$};
		\node (lpj) at (2.5,-0.3) {$p_{\ell-1}$};
		\node (lpl) at (3.5,-0.3) {$p_\ell$};
		\node (lpk1) at (5,-0.3) {$p_{k-1}$};
		\node (lpk) at (6,-0.3) {$p_k$};
		\node (lv) at (2.5,1.3) {$v$};
		\node (lu) at (3.5,1.3) {$u$};
		\draw (p1) -- (p2) (p2) -- (p) (p) -- (pl) (pl) -- (pj) (pj) -- (p') (p') -- (pk1) (pk1) -- (pk) (pl) -- (u) (pj) -- (v);
		\draw[dashed, ->-=0.5] (p1) to [bend left = 65] (u);
		\end{tikzpicture}
		\subcaption{(iv.a): The second arc is $(p_k, v)$ if $k+\ell+j$ is even and $(v, p_k)$ otherwise.}
	\end{subfigure}\hfill%
	\begin{subfigure}[b]{.45\textwidth}
		\centering
		\begin{tikzpicture}
		\node[unlabelled] (p1) at (0,0) {};
		\node[unlabelled] (p2) at (1,0) {};
		\node (p) at (1.75,0) {$\dots$};
		\node[unlabelled] (pl) at (2.5,0) {};
		\node (p') at (3.25,0) {$\dots$};
		\node[unlabelled] (pj) at (4,0) {};
		\node (p'') at (4.75,0) {$\dots$};
		\node[unlabelled] (pk1) at (5.5,0) {};
		\node[unlabelled] (pk) at (6.5,0) {};
		\node[unlabelled] (u) at (2.5,1) {};
		\node[unlabelled] (v) at (4,1) {};
		\node (lp1) at (0,-0.3) {$p_1$};
		\node (lp2) at (1,-0.3) {$p_2$};
		\node (lpl) at (2.5,-0.3) {$p_\ell$};
		\node (lpj) at (4,-0.3) {$p_j$};
		\node (lpk1) at (5.5,-0.3) {$p_{k-1}$};
		\node (lpk) at (6.5,-0.3) {$p_k$};
		\node (lu) at (2.5,1.3) {$u$};
		\node (lv) at (4,1.3) {$v$};
		\draw (p1) -- (p2) (p2) -- (p) (p) -- (pl) (pl) -- (p') (p') -- (pj) (pj) -- (p'') (p'') -- (pk1) (pk1) -- (pk) (pl) -- (u) (pj) -- (v);
		\draw[dashed, ->-=0.5] (p1) to (u);
		\end{tikzpicture}
		\subcaption{(iv.b): The second arc is $(p_k, v)$ if $k+\ell+j$ is even and $(v, p_k)$ otherwise.}
	\end{subfigure}
	
	\begin{subfigure}[b]{.45\textwidth}
		\centering
		\begin{tikzpicture}
		\node[unlabelled] (p1) at (0,0) {};
		\node[unlabelled] (p2) at (1,0) {};
		\node (p) at (1.75,0) {$\dots$};
		\node[unlabelled] (pj) at (2.5,0) {};
		\node[unlabelled] (pl) at (3.5,0) {};
		\node (p') at (4.25,0) {$\dots$};
		\node[unlabelled] (pk1) at (5,0) {};
		\node[unlabelled] (pk) at (6,0) {};
		\node (lp1) at (0,-0.3) {$p_1$};
		\node (lp2) at (1,-0.3) {$p_2$};
		\node (lpj) at (2.5,-0.3) {$p_{\ell-1}$};
		\node (lpl) at (3.5,-0.3) {$p_{\ell}$};
		\node (lpk1) at (5,-0.3) {$p_{k-1}$};
		\node (lpk) at (6,-0.3) {$p_k$};
		\draw (p1) -- (p2) (p2) -- (p) (p) -- (pl) (pl) -- (pj) (pj) -- (p') (p') -- (pk1) (pk1) -- (pk);
		\draw[dashed, ->-=0.5] (p1) to [bend left = 65] (pl);
		\end{tikzpicture}
\subcaption{(iv.c): The second arc is $(p_k,p_{\ell-1})$ if $k+\ell+j$ is even and 
           $(p_{\ell-1}, p_k)$ otherwise.}
	\end{subfigure}\hfill%
	\begin{subfigure}[b]{.45\textwidth}
		\centering
		\begin{tikzpicture}
		\node[unlabelled] (p1) at (0,0) {};
		\node[unlabelled] (p2) at (1,0) {};
		\node (p) at (1.75,0) {$\dots$};
		\node[unlabelled] (pl) at (2.5,0) {};
		\node (p') at (3.25,0) {$\dots$};
		\node[unlabelled] (pj) at (4,0) {};
		\node (p'') at (4.75,0) {$\dots$};
		\node[unlabelled] (pk1) at (5.5,0) {};
		\node[unlabelled] (pk) at (6.5,0) {};
		\node (lp1) at (0,-0.3) {$p_1$};
		\node (lp2) at (1,-0.3) {$p_2$};
		\node (lpl) at (2.5,-0.3) {$p_{\ell}$};
		\node (lpj) at (4,-0.3) {$p_j$};
		\node (lpk1) at (5.5,-0.3) {$p_{k-1}$};
		\node (lpk) at (6.5,-0.3) {$p_k$};
		\draw (p1) -- (p2) (p2) -- (p) (p) -- (pl) (pl) -- (p') (p') -- (pj) (pj) -- (p'') (p'') -- (pk1) (pk1) -- (pk);
		\draw[dashed, ->-=0.5] (p1) to [bend left = 55] (pl);
		\end{tikzpicture}
\subcaption{(iv.d): The second arc is $(p_k,p_j)$ if $k+\ell+j$ is even and 
                   $(p_j, p_k)$ otherwise.}
	\end{subfigure}
	
	\begin{subfigure}[b]{.45\textwidth}
		\centering
		\begin{tikzpicture}
		\node[unlabelled] (p1) at (0,0) {};
		\node[unlabelled] (p2) at (1,0) {};
		\node (p) at (1.75,0) {$\dots$};
		\node[unlabelled] (pj) at (2.5,0) {};
		\node[unlabelled] (pl) at (3.5,0) {};
		\node (p') at (4.25,0) {$\dots$};
		\node[unlabelled] (pk1) at (5,0) {};
		\node[unlabelled] (pk) at (6,0) {};
		\node[unlabelled] (u) at (3.5,1) {};
		\node (lp1) at (0,-0.3) {$p_1$};
		\node (lp2) at (1,-0.3) {$p_2$};
		\node (lpj) at (2.5,-0.3) {$p_{\ell-1}$};
		\node (lpl) at (3.5,-0.3) {$p_\ell$};
		\node (lpk1) at (5,-0.3) {$p_{k-1}$};
		\node (lpk) at (6,-0.3) {$p_k$};
		\node (lu) at (3.5,1.3) {$u$};
		\draw (p1) -- (p2) (p2) -- (p) (p) -- (pl) (pl) -- (pj) (pj) -- (p') (p') -- (pk1) (pk1) -- (pk) (pl) -- (u);
		\draw[dashed, ->-=0.5] (p1) to (u);
		\end{tikzpicture}
		\subcaption{(iv.e): The second arc is $(p_k,p_{\ell-1})$ if $k+\ell+j$ is odd and $(p_{\ell-1}, p_k)$ otherwise.}
	\end{subfigure}\hfill%
	\begin{subfigure}[b]{.45\textwidth}
		\centering
		\begin{tikzpicture}
		\node[unlabelled] (p1) at (0,0) {};
		\node[unlabelled] (p2) at (1,0) {};
		\node (p) at (1.75,0) {$\dots$};
		\node[unlabelled] (pl) at (2.5,0) {};
		\node[unlabelled] (u) at (2.5,1) {};
		\node (p') at (3.25,0) {$\dots$};
		\node[unlabelled] (pj) at (4,0) {};
		\node (p'') at (4.75,0) {$\dots$};
		\node[unlabelled] (pk1) at (5.5,0) {};
		\node[unlabelled] (pk) at (6.5,0) {};
		\node (lp1) at (0,-0.3) {$p_1$};
		\node (lp2) at (1,-0.3) {$p_2$};
		\node (lpl) at (2.5,-0.3) {$p_\ell$};
		\node (lu) at (2.5,1.3) {$u$};
		\node (lpj) at (4,-0.3) {$p_j$};
		\node (lpk1) at (5.5,-0.3) {$p_{k-1}$};
		\node (lpk) at (6.5,-0.3) {$p_k$};
		\draw (p1) -- (p2) (p2) -- (p) (p) -- (pl) (pl) -- (p') (p') -- (pj) (pj) -- (p'') (p'') -- (pk1) (pk1) -- (pk) (pl) -- (u);
		\draw[dashed, ->-=0.5] (p1) to (u);
		\end{tikzpicture}
		\subcaption{(iv.f): The second arc is $(p_k,p_j)$ if $k+\ell+j$ is odd and $(p_j, p_k)$ otherwise.}
	\end{subfigure}

	\caption{Obstructions $X$ for which $\overline{U(X)}$ is a tree.
		\label{treefigure}}
\end{figure}
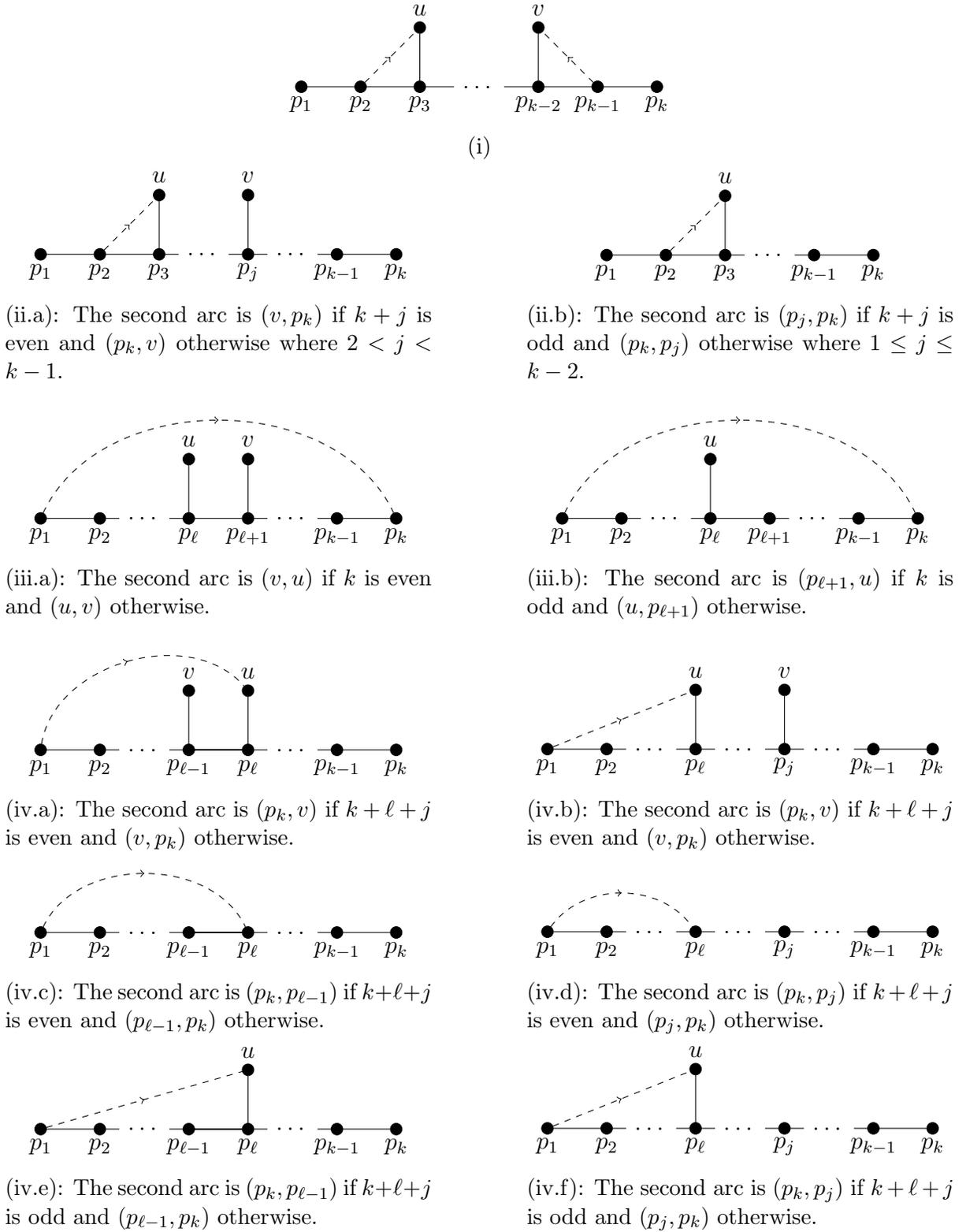

\item{(iii)} $u$ is not in $P$ and $2 \leq \ell \leq k-2$, $j=\ell+1$, $X$ or its
        dual has arcs $(p_1,p_k), (v,u)$ if either $k$ is even and $v$ is not in
        $P$ or $k$ is odd and $v$ is in $P$; otherwise $X$ or its dual has arcs
        $(p_1,p_k), (u,v)$ (See Figure~\ref{treefigure}(iii));
\item{(iv)} $3 \leq \ell \leq k-1$, $\ell-1 \leq j \leq k-2$, and $X$ or its dual
        has arcs $(p_1,u),(p_k,v)$ if either $k+\ell+j$ is even and $P$ contains
        both $u,v$ or neither, or $k+\ell+j$ is odd and $P$ contains exactly one
        of $u,v$; otherwise $X$ or its dual has arcs $(p_1,u),(v,p_k)$
        (See Figure~\ref{treefigure}(iv)).
\end{description}
\end{theorem}

\pf Suppose both $p_1,p_k$ are arc-balancing vertices.
By Lemma \ref{tree structure}, $p_1$ balances an arc between $p_2$ and a leaf $u$
adjacent to $p_3$ but not in $P$, and $p_k$ balances an arc between $p_{k-1}$ and
a leaf $v$ adjacent to $p_{k-2}$ but not in $P$. In the tree $\overline{U(X)}$ the
unique $(u,p_k)$-path avoids $p_1$ and the unique $(p_1,v)$-path avoids $p_k$.
The lengths of these two paths have the same parity so by
Proposition \ref{complement gamma}, we have either
$(u,p_1)\Gamma^*(p_k,p_1)\Gamma^*(p_k,v)$ or
$(u,p_1)\Gamma^*(p_1,p_k)\Gamma^*(p_k,v)$. In both cases
$(u,p_1)\Gamma^*(p_k,v)$ so
$(p_2,u)\Gamma(u,p_1)\Gamma^*(p_k,v)\Gamma(v,p_{k-1})$. Since the two arcs of $X$
are opposing, $X$ or its dual contains arcs $(p_2,u),(p_{k-1},v)$.
The minimality of $X$ ensures that $\overline{U(X)}$ consists of $P$ and $u, v$
and thus statement $(i)$ holds.

Suppose next that $p_1$ is an arc-balancing vertex but $p_k$ is not. By Lemma
\ref{tree structure} $p_1$ balances an arc between $p_2$ and a leaf $u$ adjacent
to $p_3$ but not in $P$. Since $p_k$ is not an arc-balancing vertex,
it is an endvertex of an arc. Let $v$ be the other endvertex. Then either
$v = p_j$ for some $1 \leq j \leq k-2$ or a leaf adjacent to some vertex in $P$.
Suppose that $v$ is a leaf adjacent to $p_j$. Then $j \notin \{1,k\}$ because $P$
is the longest path in $\overline{U(X)}$, and $j \neq 2$ because $p_2$ has no
neighbour other than $p_1,p_3$ according to Lemma \ref{tree structure}.
Moreover, $j \neq k-1$ as otherwise the arc between $p_j$ and $p_k$ is balanced,
which is not possible. So $2 < j < k-1$. In the tree $\overline{U(X)}$ the unique
$(u,p_k)$-path avoids $p_1$ and the unique $(p_1,v)$-path avoids $p_k$.
If $k+j$ is even and $v$ is not in $P$ or $k+j$ is odd and $v$ is in $P$, then
the lengths of these two paths have the same parity.
By Proposition \ref{complement gamma}, $(u,p_1)\Gamma^*(p_k,v)$ and so
$(p_2,u)\Gamma(u,p_1)\Gamma^*(p_k,v)$. Since the two arcs of $X$
are opposing, $X$ or its dual contains arcs $(p_2,u),(v,p_k)$.
Otherwise, the lengths of the two paths have the opposite parities and we have
$(p_2,u)\Gamma(u,p_1)\Gamma^*(v,p_k)$. Hence $X$ or its dual contains arcs
$(p_2,u),(p_k,v)$. The minimality of $X$ ensures that $\overline{U(X)}$ consists
of $P$ and $u, v$ and thus statement $(ii)$ holds.

It remains to consider the case when neither of $p_1,p_k$ is an arc-balancing
vertex. Suppose first that $X$ contains an arc between $p_1$ and $p_k$.
Let $u,v$ be the endvertices of the other arc.
We claim that at least one of $u,v$ is not in $P$. Indeed, if they are both in $P$
(say $u=p_i$ and $v = p_j$ where $i < j$) then $j > i+1$. It is easy to check that
$X-p_{i+1}$ cannot be completed to a local tournament, which contradicts the
minimality of $X$. So at least one of $u,v$ is not in $P$.
Assume without loss of generality that $u$ is not in $P$. Since $\overline{U(X)}$
is a caterpillar, $u$ is a leaf adjacent to some $p_\ell$ in $P$, and $v = p_j$ or
$v$ is a leaf adjacent to some $p_j$ in $P$. By reversing $\prec$ if needed we
assume that $j \geq \ell$. Since $uv$ is an unbalanced edge of $U(X)$,
$j \neq \ell$. In the tree $\overline{U(X)}$ the unique $(p_1,u)$-path avoids
$p_k$ and the unique $(v,p_k)$-path avoids $u$. If $k$ is even and $v$ is not in
$P$ or $k$ is odd and $v$ is in $P$, then the lengths of these two paths have
the same parity. By Proposition \ref{complement gamma}, $(p_1,p_k)\Gamma^*(u,v)$
and hence $X$ or its dual contains arcs $(p_1,p_k),(v,u)$.
Otherwise, the lengths of the two paths have opposite parities and we have
$(p_1,p_k)\Gamma^*(v,u)$ and $X$ or its dual contains arcs $(p_1,p_k),(u,v)$.
If $j > \ell+1$, then the two arcs are still opposing in $X-p_{\ell+1}$,
a contradiction to the assumption that $X$ is an obstruction. So $j = \ell+1$.
The minimality of $X$ ensures $\overline{U(X)}$ contains no other vertices.
Therefore statement $(iii)$ holds.

Suppose now that $X$ does not contain an arc between $p_1,p_k$. Then $p_1,p_k$ are 
incident with different arcs. Let $u,v$ be the other endvertices of the arcs 
incident with $p_1,p_k$ respectively. Then $u=p_\ell$ or is a leaf adjacent to 
some $p_\ell$ in $P$ and $v=p_j$ or is a leaf adjacent to some $p_j$ in $P$.
Since $X$ has no arc between $p_1,p_k$ and $p_1u$ is an unbalanced edge of $U(X)$,
$3 \leq \ell \leq k-1$. Similarly, $2 \leq j \leq k-2$. 
In $\overline{U(X)}$ the unique $(p_1,v)$-path avoids $p_k$ and the unique 
$(u,p_k)$-path avoids $p_1$. If $k+\ell+j$ is even and $P$ contains either both 
$u,v$ or neither, or $k+\ell+j$ is odd and $P$ contains exactly one of $u,v$, then 
the lengths of these two paths have opposite parities.
By Proposition \ref{complement gamma}, $(p_1,u)\Gamma^*(v,p_k)$ and hence 
$X$ or its dual contains arcs $(p_1,u), (p_k,v)$.
Otherwise, the lengths of the two paths have the same parity and
$(p_1,u)\Gamma^*(p_k,v)$ and $X$ or its dual contains arcs $(p_1,u), (v,p_k)$. 
If $j < \ell-1$, then the two arcs are opposing in $X-p_{\ell-1}$, contradicting
that $X$ is an obstruction. So $j \geq \ell-1$. The minimality of $X$ ensures that 
$\overline{U(X)}$ contains no other vertices. Therefore statement $(iv)$ holds.
\qed

\subsection{$\overline{U(X)}$ contains a $C_3$ but no induced $C_4$}

We now examine obstructions $X$ that do not contain cut-vertices and for which 
$\overline{U(X)}$ contains cycles. By Corollary \ref{disconnected}, 
$\overline{U(X)}$ is connected. We know from Lemma \ref{largest size of cycles}
that any induced cycle in $\overline{U(X)}$ is of length 3, 4 or 5, and also
from Lemma \ref{C5lemma} that if $\overline{U(X)}$ contains an induced cycle 
of length 5 then it does not contain an induced cycle of length 3 or 4.
  
We divide our discussion into four cases: $\overline{U(X)}$ contains a $C_3$
but no induced $C_4$; $\overline{U(X)}$ contains an induced $C_4$ but no $C_3$;
$\overline{U(X)}$ contains both a $C_3$ and an induced $C_4$; and 
$\overline{U(X)}$ contains an induced $C_5$. These four cases will be treated 
separately.

\begin{lemma}\label{c3noc4lemma}
Let $X$ be an obstruction that contains no cut-vertices. Suppose $\overline{U(X)}$
contains a $C_3$ but no induced $C_4$. Then the $C_3$ is the only cycle in 
$\overline{U(X)}$ and any vertex not on the $C_3$ is a leaf adjacent to a vertex on 
the $C_3$ and incident with an arc. Moreover, any vertex on the $C_3$ is adjacent 
to a vertex not on it.
\end{lemma}
\pf Since $\overline{U(X)}$ contains a $C_3$ but no induced $C_4$, by Lemmas 
\ref{largest size of cycles}, \ref{onec3}, and \ref{C5lemma}, the $C_3$ is the 
unique cycle in $\overline{U(X)}$. Let $C: v_1v_2v_3$ the unique cycle. Consider 
a vertex $u$ that is not on $C$. By Lemma \ref{odd cycles}, $u$ is adjacent 
to a vertex on $C$. Since $C$ is the unique cycle in $\overline{U(X)}$, $u$ 
must be a leaf. Clearly, $u$ is not a cut-vertex of $\overline{U(X)}$ and by 
assumption is not a cut-vertex of $U(X)$. If $u$ is an arc-balancing vertex, then 
$u$ balances an arc incident with a vertex on $C$. Thus the other two 
vertices of $C$ must be adjacent to both endvertices of the arc, 
a contradiction to the fact $C$ is the unique cycle in $\overline{U(X)}$. 
So $u$ is not an arc-balancing vertex and therefore by Lemma 
\ref{vertex classification} it is incident with an arc.
 
It remains to show that each vertex on $C$ is adjacent to a vertex not on it.
Suppose on the contrary that $v_1$ is not adjacent to a vertex not on $C$. 
By Lemma \ref{2 non-adj}, $v_2$ and $v_3$ each has two non-neighbours.
Clearly, the non-neighbours of $v_2$ and of $v_3$ are not on $C$. 
We know from the above they are endvertices of arcs. Since $v_1$ is adjacent to 
none of them, $v_1$ is not an arc-balancing vertex. By assumption $v_1$ is not
a cut-vertex of $U(X)$. It cannot be a cut-vertex of $\overline{U(X)}$ because it
is adjacent only to $v_2,v_3$ (which are adjacent). This is a contradiction to
Lemma \ref{vertex classification}.
\qed

\begin{theorem} \label{c3only}
Let $X$ be an obstruction that contains no cut-vertices. Suppose $\overline{U(X)}$ 
contains a $C_3$ but no induced $C_4$. Then $\overline{U(X)}$ is one of the graphs 
in Figure~\ref{onlyc3} and $X$ or its dual contains the dotted arcs.
\begin{figure}[H]
	\captionsetup[subfigure]{labelformat=empty}
	\begin{subfigure}[b]{.5\textwidth}
		\centering
		\begin{tikzpicture}
			\node[unlabelled] (1) at (0.5,0.866) {};
			\node[unlabelled] (3) at (0,0) {};
			\node[unlabelled] (2) at (1,0) {};
			\node[unlabelled] (u) at (0.5,1.866) {};
			\node[unlabelled] (w) at (-1,0) {};
			\node[unlabelled] (v) at (2,0) {};
			\node (l1) at (0.2,0.866) {$v_1$};
			\node (l3) at (0,-0.3) {$v_3$};
			\node (l2) at (1,-0.3) {$v_2$};
			\node (lu) at (0.5,2.166) {$u$};
			\node (lw) at (-1,-0.3) {$w$};
			\node (lv) at (2,-0.3) {$v$};
			\draw (1) -- (2) (2) -- (3) (3) -- (1) (1) -- (u) (2) -- (v) (3) -- (w);
			\draw[dashed, ->-=0.5] (u) to (v);
			\draw[dashed, ->-=0.5] (w) to (1);
		\end{tikzpicture}
		\subcaption{(i)}
	\end{subfigure}%
	\begin{subfigure}[b]{.5\textwidth}
		\centering
		\begin{tikzpicture}
			\node[unlabelled] (1) at (0.5,0.866) {};
			\node[unlabelled] (3) at (0,0) {};
			\node[unlabelled] (2) at (1,0) {};
			\node[unlabelled] (u) at (1,1.732) {};
			\node[unlabelled] (z) at (0,1.732) {};
			\node[unlabelled] (w) at (-1,0) {};
			\node[unlabelled] (v) at (2,0) {};
			\node (l1) at (0.2,0.866) {$v_1$};
			\node (l3) at (0,-0.3) {$v_3$};
			\node (l2) at (1,-0.3) {$v_2$};
			\node (lu) at (1,2.032) {$u$};
			\node (lz) at (0,2.032) {$z$};
			\node (lw) at (-1,-0.3) {$w$};
			\node (lv) at (2,-0.3) {$v$};
			\draw (1) -- (2) (2) -- (3) (3) -- (1) (1) -- (u) (2) -- (v) (3) -- (w)	(1) -- (z);
			\draw[dashed, ->-=0.5] (u) to (v);
			\draw[dashed, ->-=0.5] (z) to (w);
		\end{tikzpicture}
		\subcaption{(ii)}
	\end{subfigure}
\caption{Obstructions $X$ for which $\overline{U(X)}$ contains a $C_3$ but no
induced $C_4$. \label{onlyc3}}
\end{figure}
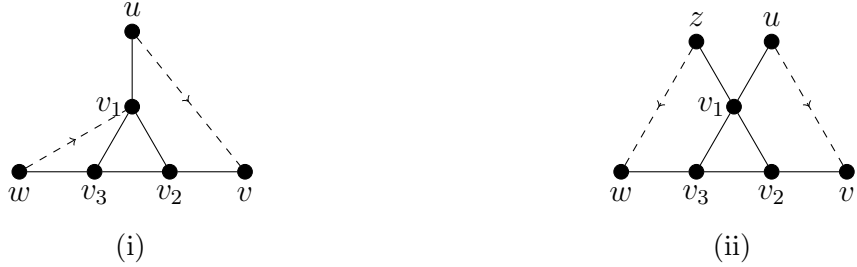
\end{theorem}
\pf Suppose $X$ is an obstruction. Let $C:\ v_1v_2v_3$ be the unique $C_3$ in 
$\overline{U(X)}$. By Lemma \ref{c3noc4lemma}, each vertex of $C$ is 
adjacent to a vertex not on it and each vertex not on $C$ is a leaf adjacent
to a vertex of $C$. Let $u, v, w$ be vertices adjacent to $v_1, v_2, v_3$ 
respectively but not on $C$. By Lemma \ref{c3noc4lemma}, each of $u, v, w$ is incident with an arc. Since $X$ contains exactly two arcs, there must be an arc 
with both endvertices among $u, v, w$. Without loss of generality, assume there is 
an arc between $u$ and $v$. By possibly considering the dual of $X$, let $(u,v)$ be
an arc. On the other hand, let $z$ denote the other endvertex of the arc incident 
with $w$. First suppose $z$ is on $C$. Then $z \in \{v_1, v_2\}$. 
Without loss of generality, assume $z=v_1$. Then 
$(u,v)\Gamma (v_2,u)\Gamma (u,v_3)\Gamma (w, u) \Gamma(v_1, w)=(z, w)$. 
Since the two arcs in $X$ are opposing, the second arc must be $(w, z)=(w,v_1)$. 
Thus $\overline{U(X)}$ is Figure~\ref{onlyc3}(i). Otherwise, $z$ is not on the $C$.
By Lemma \ref{c3noc4lemma}, $z$ is a leaf adjacent to a vertex on $C$. 
Clearly, $z$ cannot be adjacent to $v_3$ because otherwise the arc between $w$ and 
$z$ would be balanced. Hence, assume without loss of generality that $z$ is 
adjacent to $v_1$. If $z=u$, then $u, z$ belong to one component of 
$\overline{U(X-v_1)}$ and $v, w$ belong to another, so $uv$ and $wz$ belong to 
the same implication class of $U(X-v_1)$, contradicting Theorem 
\ref{different implication classes theorem}. Hence, $z \neq u$. In this case, 
we have $(u, v)\Gamma^*(w, u)\Gamma(v_1, w)\Gamma(w, z)$. Hence, the second arc is 
$(z, w)$. Thus $\overline{U(X)}$ is Figure~\ref{onlyc3}(ii).
\qed

\subsection{$\overline{U(X)}$ contains an induced $C_4$ but no $C_3$}

We consider next the case when $\overline{U(X)}$ contains an induced $C_4$ but no 
$C_3$. Since $U(X)$ is a proper circular-arc graph, by Theorem \ref{Tucker} any 
induced $C_4$ in $\overline{U(X)}$ contains at most two cut-vertices of 
$\overline{U(X)}$.

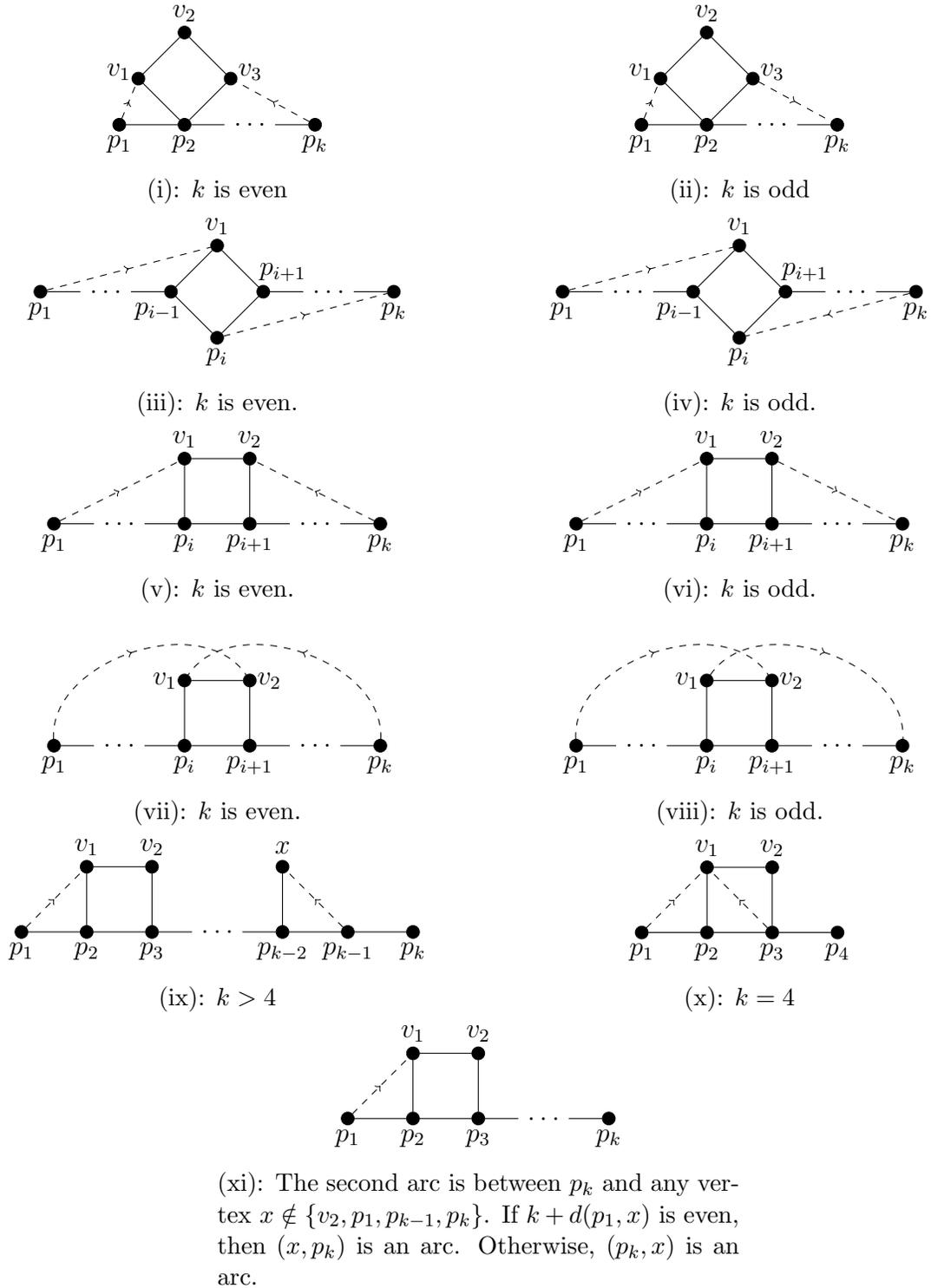
\begin{figure}[H]
		\centering
		\captionsetup[subfigure]{labelformat=empty}
		\begin{subfigure}[b]{.5\textwidth}
			\centering
			\begin{tikzpicture}
			\node[unlabelled] (1) at (-0.7071,0.7071) {};
			\node[unlabelled] (2) at (0,1.4142) {};
			\node[unlabelled] (3) at (0.7071,0.7071) {};
			\node[unlabelled] (4) at (0,0) {};
			\node[unlabelled] (p1) at (-1,0) {};
			\node (p) at (1,0) {$\ldots$};
			\node[unlabelled] (pk) at (2,0) {};
			\node (l1) at (-1.0071,0.8071) {$v_1$};
			\node (l2) at (0,1.7142) {$v_2$};
			\node (l3) at (1.0071,0.8071) {$v_3$};
			\node (l4) at (0,-0.3) {$p_2$};
			\node (lp1) at (-1,-0.3) {$p_1$};
			\node (lpk) at (2,-0.3) {$p_k$};
			\draw (1) -- (2) (2) -- (3) (3) -- (4) (4) -- (1) (4) -- (p1) (4) -- (p) (p) -- (pk);
			\draw[dashed, ->-=0.5] (p1) to (1);
			\draw[dashed, ->-=0.5] (pk) to (3);
			\end{tikzpicture}
			\subcaption{(i): $k$ is even}
		\end{subfigure}%
		\begin{subfigure}[b]{.5\textwidth}
			\centering
			\begin{tikzpicture}
			\node[unlabelled] (1) at (-0.7071,0.7071) {};
			\node[unlabelled] (2) at (0,1.4142) {};
			\node[unlabelled] (3) at (0.7071,0.7071) {};
			\node[unlabelled] (4) at (0,0) {};
			\node[unlabelled] (p1) at (-1,0) {};
			\node (p) at (1,0) {$\ldots$};
			\node[unlabelled] (pk) at (2,0) {};
			\node (l1) at (-1.0071,0.8071) {$v_1$};
			\node (l2) at (0,1.7142) {$v_2$};
			\node (l3) at (1.0071,0.8071) {$v_3$};
			\node (l4) at (0,-0.3) {$p_2$};
			\node (lp1) at (-1,-0.3) {$p_1$};
			\node (lpk) at (2,-0.3) {$p_k$};
			\draw (1) -- (2) (2) -- (3) (3) -- (4) (4) -- (1) (4) -- (p1) (4) -- (p) (p) -- (pk);
			\draw[dashed, ->-=0.5] (p1) to (1);
			\draw[dashed, ->-=0.5] (3) to (pk);
			\end{tikzpicture}
			\subcaption{(ii): $k$ is odd}
		\end{subfigure}
		
		\begin{subfigure}[b]{.5\textwidth}
			\centering
			\begin{tikzpicture}
			\node[unlabelled] (1) at (0.7071, 0.7071) {};
			\node[unlabelled] (2) at (1.4142, 0) {};
			\node[unlabelled] (3) at (0.7071, -0.7071) {};
			\node[unlabelled] (4) at (0,0) {};
			\node[unlabelled] (p1) at (-2,0) {};
			\node (p) at (-1,0) {$\ldots$};
			\node (q) at (2.4142,0) {$\ldots$};
			\node[unlabelled] (pk) at (3.4142,0) {};
			\node (l1) at (0.7071, 1.0071) {$v_1$};
			\node (l2) at (1.7142, 0.3) {$p_{i+1}$};
			\node (l3) at (0.7071, -1.0071) {$p_i$};
			\node (l4) at (-0.2,-0.3) {$p_{i-1}$};
			\node (lp1) at (-2,-0.3) {$p_1$};
			\node (lpk) at (3.4142,-0.3) {$p_k$};
			\draw (1) -- (2) (2) -- (3) (3) -- (4) (4) -- (1) (4) -- (p) (p) -- (p1) (2) -- (q) (q) -- (pk);
			\draw[dashed, ->-=0.5] (p1) to (1);
			\draw[dashed, ->-=0.5] (3) to (pk);
			\end{tikzpicture}
			\subcaption{(iii): $k$ is even.}
		\end{subfigure}%
		\begin{subfigure}[b]{.5\textwidth}
			\centering
			\begin{tikzpicture}
			\node[unlabelled] (1) at (0.7071, 0.7071) {};
			\node[unlabelled] (2) at (1.4142, 0) {};
			\node[unlabelled] (3) at (0.7071, -0.7071) {};
			\node[unlabelled] (4) at (0,0) {};
			\node[unlabelled] (p1) at (-2,0) {};
			\node (p) at (-1,0) {$\ldots$};
			\node (q) at (2.4142,0) {$\ldots$};
			\node[unlabelled] (pk) at (3.4142,0) {};
			\node (l1) at (0.7071, 1.0071) {$v_1$};
			\node (l2) at (1.7142, 0.3) {$p_{i+1}$};
			\node (l3) at (0.7071, -1.0071) {$p_i$};
			\node (l4) at (-0.2,-0.3) {$p_{i-1}$};
			\node (lp1) at (-2,-0.3) {$p_1$};
			\node (lpk) at (3.4142,-0.3) {$p_k$};
			\draw (1) -- (2) (2) -- (3) (3) -- (4) (4) -- (1) (4) -- (p) (p) -- (p1) (2) -- (q) (q) -- (pk);
			\draw[dashed, ->-=0.5] (p1) to (1);
			\draw[dashed, ->-=0.5] (pk) to (3);
			\end{tikzpicture}
			\subcaption{(iv): $k$ is odd.}
		\end{subfigure}
	
		\begin{subfigure}[b]{.5\textwidth}
			\centering
			\begin{tikzpicture}
			\node[unlabelled] (1) at (0,1) {};
			\node[unlabelled] (2) at (1,1) {};
			\node[unlabelled] (3) at (1,0) {};
			\node[unlabelled] (4) at (0,0) {};
			\node[unlabelled] (p1) at (-2,0) {};
			\node (p) at (-1,0) {$\ldots$};
			\node (q) at (2,0) {$\ldots$};
			\node[unlabelled] (pk) at (3,0) {};
			\node (l1) at (0,1.3) {$v_1$};
			\node (l2) at (1,1.3) {$v_2$};
			\node (l3) at (1,-0.3) {$p_{i+1}$};
			\node (l4) at (0,-0.3) {$p_i$};
			\node (lp1) at (-2,-0.3) {$p_1$};
			\node (lpk) at (3,-0.3) {$p_k$};
			\draw (1) -- (2) (2) -- (3) (3) -- (4) (4) -- (1) (4) -- (p) (p) -- (p1) (3) -- (q) (q) -- (pk);
			\draw[dashed, ->-=0.5] (p1) to (1);
			\draw[dashed, ->-=0.5] (pk) to (2);
			\end{tikzpicture}
			\subcaption{(v): $k$ is even.}
		\end{subfigure}%
		\begin{subfigure}[b]{.5\textwidth}
			\centering
			\begin{tikzpicture}
			\node[unlabelled] (1) at (0,1) {};
			\node[unlabelled] (2) at (1,1) {};
			\node[unlabelled] (3) at (1,0) {};
			\node[unlabelled] (4) at (0,0) {};
			\node[unlabelled] (p1) at (-2,0) {};
			\node (p) at (-1,0) {$\ldots$};
			\node (q) at (2,0) {$\ldots$};
			\node[unlabelled] (pk) at (3,0) {};
			\node (l1) at (0,1.3) {$v_1$};
			\node (l2) at (1,1.3) {$v_2$};
			\node (l3) at (1,-0.3) {$p_{i+1}$};
			\node (l4) at (0,-0.3) {$p_i$};
			\node (lp1) at (-2,-0.3) {$p_1$};
			\node (lpk) at (3,-0.3) {$p_k$};
			\draw (1) -- (2) (2) -- (3) (3) -- (4) (4) -- (1) (4) -- (p) (p) -- (p1) (3) -- (q) (q) -- (pk);
			\draw[dashed, ->-=0.5] (p1) to (1);
			\draw[dashed, ->-=0.5] (2) to (pk);
			\end{tikzpicture}
			\subcaption{(vi): $k$ is odd.}
		\end{subfigure}
		
		\begin{subfigure}[b]{.5\textwidth}
			\centering
			\begin{tikzpicture}
			\node[unlabelled] (1) at (0,1) {};
			\node[unlabelled] (2) at (1,1) {};
			\node[unlabelled] (3) at (1,0) {};
			\node[unlabelled] (4) at (0,0) {};
			\node[unlabelled] (p1) at (-2,0) {};
			\node (p) at (-1,0) {$\ldots$};
			\node (q) at (2,0) {$\ldots$};
			\node[unlabelled] (pk) at (3,0) {};
			\node (l1) at (-0.3,1) {$v_1$};
			\node (l2) at (1.3,1) {$v_2$};
			\node (l3) at (1,-0.3) {$p_{i+1}$};
			\node (l4) at (0,-0.3) {$p_i$};
			\node (lp1) at (-2,-0.3) {$p_1$};
			\node (lpk) at (3,-0.3) {$p_k$};
			\draw (1) -- (2) (2) -- (3) (3) -- (4) (4) -- (1) (4) -- (p) (p) -- (p1) (3) -- (q) (q) -- (pk);
			\draw[dashed, ->-=0.5, bend left = 75] (p1) to (2);
			\draw[dashed, ->-=0.5, bend right = 75] (pk) to (1);
			\end{tikzpicture}
			\subcaption{(vii): $k$ is even.}
		\end{subfigure}%
		\begin{subfigure}[b]{.5\textwidth}
			\centering
			\begin{tikzpicture}
			\node[unlabelled] (1) at (0,1) {};
			\node[unlabelled] (2) at (1,1) {};
			\node[unlabelled] (3) at (1,0) {};
			\node[unlabelled] (4) at (0,0) {};
			\node[unlabelled] (p1) at (-2,0) {};
			\node (p) at (-1,0) {$\ldots$};
			\node (q) at (2,0) {$\ldots$};
			\node[unlabelled] (pk) at (3,0) {};
			\node (l1) at (-0.3,1) {$v_1$};
			\node (l2) at (1.3,1) {$v_2$};
			\node (l3) at (1,-0.3) {$p_{i+1}$};
			\node (l4) at (0,-0.3) {$p_i$};
			\node (lp1) at (-2,-0.3) {$p_1$};
			\node (lpk) at (3,-0.3) {$p_k$};
			\draw (1) -- (2) (2) -- (3) (3) -- (4) (4) -- (1) (4) -- (p) (p) -- (p1) (3) -- (q) (q) -- (pk);
			\draw[dashed, ->-=0.5, bend left = 75] (p1) to (2);
			\draw[dashed, ->-=0.5, bend left = 75] (1) to (pk);
			\end{tikzpicture}
			\subcaption{(viii): $k$ is odd.}
		\end{subfigure}

		\begin{subfigure}[b]{.5\textwidth}
			\centering
			\begin{tikzpicture}
			\node[unlabelled] (1) at (0,1) {};
			\node[unlabelled] (2) at (1,1) {};
			\node[unlabelled] (3) at (1,0) {};
			\node[unlabelled] (4) at (0,0) {};
			\node[unlabelled] (p1) at (-1,0) {};
			\node (p) at (2,0) {$\ldots$};
			\node[unlabelled] (p3) at (3,0) {};
			\node[unlabelled] (p4) at (4,0) {};
			\node[unlabelled] (p5) at (5,0) {};
			\node[unlabelled] (x) at (3,1) {};
			\node (l1) at (0,1.3) {$v_1$};
			\node (l2) at (1,1.3) {$v_2$};
			\node (l3) at (1,-0.3) {$p_3$};
			\node (l4) at (0,-0.3) {$p_2$};
			\node (lp1) at (-1,-0.3) {$p_1$};
			\node (lp3) at (3,-0.3) {$p_{k-2}$};
			\node (lp4) at (4,-0.3) {$p_{k-1}$};
			\node (lp5) at (5,-0.3) {$p_k$};
			\node (lx) at (3,1.3) {$x$};
			\draw (1) -- (2) (2) -- (3) (3) -- (4) (4) -- (1) (4) -- (p1) (3) -- (p) (p) -- (p3) (p3) -- (p4) (p4) -- (p5) (p3) -- (x);
			\draw[dashed, ->-=0.5] (p1) to (1);
			\draw[dashed, ->-=0.5] (p4) to (x);
			\end{tikzpicture}
			\subcaption{(ix): $k > 4$}
		\end{subfigure}%
		\begin{subfigure}[b]{.5\textwidth}
			\centering
			\begin{tikzpicture}
			\node[unlabelled] (1) at (0,1) {};
			\node[unlabelled] (2) at (1,1) {};
			\node[unlabelled] (3) at (1,0) {};
			\node[unlabelled] (4) at (0,0) {};
			\node[unlabelled] (p1) at (-1,0) {};
			\node[unlabelled] (pk) at (2,0) {};
			\node (l1) at (0,1.3) {$v_1$};
			\node (l2) at (1,1.3) {$v_2$};
			\node (l3) at (1,-0.3) {$p_3$};
			\node (l4) at (0,-0.3) {$p_2$};
			\node (lp1) at (-1,-0.3) {$p_1$};
			\node (lpk) at (2,-0.3) {$p_4$};
			\draw (1) -- (2) (2) -- (3) (3) -- (4) (4) -- (1) (4) -- (p1) (3) -- (pk);
			\draw[dashed, ->-=0.5] (p1) to (1);
			\draw[dashed, ->-=0.5] (3) to (1);
			\end{tikzpicture}
			\subcaption{(x): $k = 4$}
		\end{subfigure}
	
		\begin{subfigure}[b]{.5\textwidth}
			\centering
			\begin{tikzpicture}
			\node[unlabelled] (1) at (0,1) {};
			\node[unlabelled] (2) at (1,1) {};
			\node[unlabelled] (3) at (1,0) {};
			\node[unlabelled] (4) at (0,0) {};
			\node[unlabelled] (p1) at (-1,0) {};
			\node (p) at (2,0) {$\ldots$};
			\node[unlabelled] (pk) at (3,0) {};
			\node (l1) at (0,1.3) {$v_1$};
			\node (l2) at (1,1.3) {$v_2$};
			\node (l3) at (1,-0.3) {$p_3$};
			\node (l4) at (0,-0.3) {$p_2$};
			\node (lp1) at (-1,-0.3) {$p_1$};
			\node (lpk) at (3,-0.3) {$p_k$};
			\draw (1) -- (2) (2) -- (3) (3) -- (4) (4) -- (1) (4) -- (p1) (3) -- (p) (p) -- (pk);
			\draw[dashed, ->-=0.5] (p1) to (1);
			\end{tikzpicture}
			\subcaption{(xi): The second arc is between $p_k$ and any vertex $x \notin \{v_2, p_1, p_{k-1}, p_k\}$. If $k + d(p_1, x)$ is even, then $(x, p_k)$ is an arc. Otherwise, $(p_k, x)$ is an arc.}
		\end{subfigure}

		\caption{Obstructions $X$ for which $\overline{U(X)}$ contains a unique induced $C_4$ but no $C_3$. \label{onec4only}}
	\end{figure}

\begin{theorem} \label{onec4}
Let $X$ be an obstruction that contains no cut-vertices. Suppose $\overline{U(X)}$ 
contains a unique induced $C_4$ but no $C_3$. Then $\overline{U(X)}$ is one of 
the graphs in Figure~\ref{onec4only} and $X$ or its dual contains the dotted arcs.
\end{theorem}

\pf Let $C: v_1v_2v_3v_4$ be the unique induced
$C_4$ in $\overline{U(X)}$. Since $U(X)$ and $\overline{U(X)}$ are both connected,
at least one vertex on $C$ is adjacent to a vertex not on $C$. Moreover, since $C$
is the unique cycle in $\overline{U(X)}$, any vertex on $C$ that is adjacent to
a vertex not on $C$ is a cut-vertex of $\overline{U(X)}$. So $C$ contains at least
one cut-vertex.

Suppose that only one vertex on $C$ is a cut-vertex of $\overline{U(X)}$.
Without loss of generality assume $v_4$ is such a vertex. We claim that $v_1, v_3$
are incident with different arcs. Indeed, since $v_1, v_3$ are not cut-vertices,
by Proposition \ref{vertex classification} they are either incident with arcs or
arc-balancing vertices. If $v_1$ is an arc-balancing vertex, then it balances an arc
incident with $v_2$ or $v_4$. Note that $v_3$ is adjacent to both $v_2$ and $v_4$
so $v_3$ must be adjacent to the other endvertex of the arc balanced by $v_1$,
which is not possible. Hence $v_1$ is not an arc-balancing vertex. By symmetry
$v_3$ is not an arc-balancing vertex either. Therefore each of $v_1, v_3$ is
incident with an arc. Since $v_1, v_3$ have the same neighbourhood, there cannot
be an arc between $v_1, v_3$, which implies that $v_1, v_3$ are incident with
different arcs as claimed.

Let $u, w$ denote the other endvertices of the arcs incident with $v_1, v_3$
respectively. Clearly, $u, w$ are not on $C$. Since $v_4$ is the unique cut-vertex
in $C$, each of $u, w$ belongs to a component of $\overline{U(X-v_4)}$ that does
not contain a vertex of $C$. According to Theorem
\ref{different implication classes theorem}, $v_1u$ and $v_3w$ belong to different
implication classes of $U(X-v_4)$. Since $v_1,v_3$ are in the same component of
$\overline{U(X-v_4)}$, $u, w$ are in different components of $\overline{U(X-v_4)}$
by Theorem \ref{structure}. The vertex $v_2$ is not a cut-vertex so it is
an arc-balancing vertex by Proposition \ref{vertex classification}. Without loss of
generality, assume $v_2$ balances the arc between $v_1$ and $u$. Then $u$ must be
a leaf adjacent to $v_4$, as otherwise there is a vertex adjacent to $u$ but not to
$v_1$, a contradiction to the fact $v_2$ balances the arc between $v_1$ and $u$.
Let $P:u=p_1, p_2, \dots, p_k=w$ be a shortest $(u,w)$-path. Such a path exists
because $\overline{U(X)}$ is connected. It is easy to see that $p_2 = v_4$.
By possibly considering the dual of $X$, assume $(p_1,v_1)$ is an arc of $X$.
If $k$ is even, then
$(p_1,v_1)\Gamma(v_2,p_1)\Gamma(p_2,v_2)\Gamma^*(p_k,v_2)\Gamma(v_3,p_k)$
by Proposition \ref{complement gamma}. The two arcs of $X$ are opposing, so
the second arc is $(p_k,v_3)$. The minimality of $X$ ensures $\overline{U(X)}$ is
Figure~\ref{onec4only}(i) and $X$ contains the dotted arcs.
Otherwise, $k$ is odd and the second arc is $(v_3,p_k)$, so $\overline{U(X)}$ is
Figure~\ref{onec4only}(ii) and $X$ contains the dotted arcs.

Suppose that exactly two vertices of $C$ are cut-vertices of $\overline{U(X)}$.
We consider first the case when the two cut-vertices of $\overline{U(X)}$ on $C$ 
are non-consecutive, say $v_2$ and $v_4$.
We claim that $v_1, v_3$ are incident with different arcs. Since $v_1, v_3$ are not
cut-vertices, neither of them is adjacent to any vertex not on $C$. In particular, 
if $v_1$ is an arc-balancing vertex, it must balance an arc incident with $v_2$ or 
$v_4$, and the other endvertex is adjacent to $v_3$ but not to $v_1$. Such a vertex
does not exist, so $v_1$ is not an arc-balancing vertex. Similarly, $v_3$ is not an
arc-balancing vertex. By Proposition \ref{vertex classification}, $v_1, v_3$ are 
incident with arcs. Moreover, since $v_1, v_3$ share the same neighbourhood, they 
must be incident with different arcs as claimed. Let $H_1$ denote a component of 
$\overline{U(X-v_2)}$ not containing vertices on $C$, and $H_2$ denote a component
of $\overline{U(X-v_4)}$ not containing vertices on $C$. 
Since $C$ is the unique cycle in $\overline{U(X)}$, $H_1, H_2$ are vertex-disjoint
trees. Let $u, w$ be leaves of $\overline{U(X)}$ in $H_1, H_2$ respectively. 
Clearly, neither $u$ nor $w$ can balance the arc incident with $v_1$ because 
otherwise the other endvertex would be adjacent to both of $v_2, v_4$ and thus 
would be $v_3$, a contradiction to the fact that $v_1, v_3$ are incident with 
different arcs. Similarly, neither $u$ nor $w$ can balance the arc incident with 
$v_3$. Hence each 
of $u, w$ is incident with an arc by Proposition \ref{vertex classification}. 
Without loss of generality, assume there is an arc between $u, v_3$ and an arc 
between $w, v_1$. By the choice of $u$ and $w$, there is a $(w,u)$-path that 
contains $v_3$ but not $v_1$. Let $P: w=p_1, \dots, p_k=u$ be a shortest 
$(w,u)$-path where $p_{i-1} = v_4$, $p_i = v_3$, and $p_{i+1} = v_2$ for some $i$. 
By possibly considering the dual of $X$, assume $(w,v_1) = (p_1,v_1)$ is an arc. 
Suppose $k$ is even. If $i$ is even, $(p_1,v_1)\Gamma^*(p_1,p_k)\Gamma^*(p_k,p_i)$ 
by Proposition \ref{complement gamma}. So the second arc is $(p_i,p_k) = (v_3,u)$. 
If instead $i$ is odd, then $(p_1,v_1)\Gamma^*(p_k,p_1)\Gamma^*(p_k,p_i)$. 
The second arc is again $(p_i,p_k) = (v_3,u)$. The minimality of $X$ ensures 
$\overline{U(X)}$ is Figure~\ref{onec4only}(iii) and $X$ contains the dotted arcs. 
Otherwise, $k$ is odd and $\overline{U(X)}$ is Figure~\ref{onec4only}(iv) and 
$X$ contains the dotted arcs.

We now consider the case when the two cut-vertices of $\overline{U(X)}$ on $C$ are 
consecutive, say $v_3$ and $v_4$. First suppose both $v_1$ and $v_2$ are incident 
with arcs. Clearly, $v_1, v_2$ are incident with different arcs. Let $u, w$ be the 
other two endvertices of the arcs. By a similar argument as above, $u, w$ are 
leaves in components of $\overline{U(X-p_3)}, \overline{U(X-p_4)}$ respectively.
Let $P: w=p_1, \dots, p_k=u$ be a shortest $(w,u)$-path where $p_i=v_4$ and 
$p_{i+1} = v_3$ for some $i$. There are two possibilities: either $w$ or $u$ is 
the endvertex of the arc incident with $v_1$. Suppose there is an arc between
$w$ and $v_1$. By possibly considering the dual of $X$, assume 
$(p_1,v_1) = (w,v_1)$ is an arc in $X$. Suppose $k$ is even. If $i$ is odd, then 
Proposition \ref{complement gamma} implies 
$(p_1,v_1)\Gamma^*(p_1,p_k)\Gamma^*(v_2,p_k)$. If $i$ is even, then 
$(p_1,v_1)\Gamma^*(p_k,p_1)\Gamma^*(v_2,p_k)$. In either case, the second arc is 
$(p_k,v_2)$, so $\overline{U(X)}$ is Figure~\ref{onec4only}(v) and $X$ contains 
the dotted arcs. Otherwise, $k$ is odd and $\overline{U(X)}$ is 
Figure~\ref{onec4only}(vi) and $X$ contains the dotted arcs. On the other hand, 
suppose $(p_1,v_2)$ is an arc in $X$. Suppose $k$ is even. Then we have
$(p_1,v_2)\Gamma^*(p_1,p_k)\Gamma^*(v_1,p_k)$ if $i$ is even, and
$(p_1,v_2)\Gamma^*(p_k,p_1)\Gamma^*(v_1,p_k)$ if $i$ is odd. 
In either case, the second arc is $(p_k,v_1)$, so $\overline{U(X)}$ is 
Figure~\ref{onec4only}(vii) and $X$ contains the dotted arcs. Otherwise, $k$ is odd
and $\overline{U(X)}$ is Figure~\ref{onec4only}(viii) and $X$ contains the dotted 
arcs.

Suppose that one of $v_1, v_2$ is not incident with an arc. Without loss of 
generality, assume it is $v_2$. Then $v_2$ is an arc-balancing vertex by 
Proposition \ref{vertex classification}. Since $v_3$ is a cut-vertex, it is 
adjacent to a vertex $x$ not on $C$. So, if $v_2$ balances an arc incident with 

$v_3$, then the other endvertex must be adjacent to both $v_4$ and $x$, 
contradicting the fact that $C$ is the unique cycle. Hence $v_2$ balances an arc 
incident with $v_1$. Since $v_1$ is adjacent only to $v_2$ and $v_4$, the other 
endvertex $w$ is a leaf adjacent to $v_4$. Without loss of generality, assume 
$(w,v_1)$ is an arc. Since $v_3$ is a cut-vertex, there is a component $H$ of 
$\overline{U(X-v_3)}$ not containing the vertices on $C$. Let $u$ be a vertex of 
maximal distance from $v_3$ in $H$, and let $P: w=p_1, \dots, p_k=u$ be a shortest 
$(w,u)$-path in $\overline{U(X)}$. Clearly, $p_2 = v_4$ and $p_3 = v_3$. Moreover, 
since $C$ is the unique cycle and $u$ is of maximal distance from $v_3$ in $H$, 
$u$ is a leaf. First suppose $u$ balances an arc incident with $p_{k-1}$. 
There are two cases depending on whether or not $k>4$. 
If $k>4$, then $p_{k-1} \neq v_3$, so the other endvertex is a leaf $x$ adjacent 
to $p_{k-2}$. If $k=4$, then $p_{k-1} = v_3$, so the other endvertex is $v_1$, 
because it must be adjacent to both $v_2$ and $v_4$ and $C$ is the unique cycle. 
In either case, we have $d(v_1,p_k) + d(p_1,x) = 2k - 3$, so one of 
$d(v_1,p_k), d(p_1,x)$ is even and the other is odd. If $d(v_1,p_k)$ is 
even and $d(p_1,x)$ is odd, then Proposition \ref{complement gamma} implies 
$(p_1,v_1)\Gamma^*(p_1,p_k)\Gamma^*(p_k,x)\Gamma(x,p_{k-1})$. On the other hand,
if $d(v_1,p_k)$ is odd and $d(p_1,x)$ is even, then 
$(p_1,v_1)\Gamma^*(p_k,p_1)\Gamma^*(p_k,x)\Gamma(x,p_{k-1})$. In either case, 
the second arc must be $(p_{k-1},x)$. Thus, $\overline{U(X)}$ is 
Figure~\ref{onec4only}(ix) if $k>4$ and is Figure~\ref{onec4only}(x) if $k=4$, and 
$X$ contains the dotted arcs.

Otherwise, $u$ is incident with an arc by Proposition \ref{vertex classification}. 
Let $x$ denote the other endvertex. Since $v_2$ is not incident with an arc, 
$x \neq v_2$. Theorem \ref{different implication classes theorem} implies $p_1v_1$ 
and $p_kx$ belong to different implication classes of $U(X-p_2)$, so $x \neq p_1$ 
by Theorem \ref{structure}. Thus, $x \notin \{v_2, p_1, p_{k-1}, p_k\}$. Suppose 
$k + d(p_1, x)$ is even. Since $d(v_1, p_k) + d(p_1, x) = (k-1) + d(p_1, x)$, one 
of $d(v_1, p_k)$ and $d(p_1, x)$ is even and the other is odd. If $d(v_1, p_k)$ is 
even and $d(p_1, x)$ is odd, then Proposition \ref{complement gamma} implies 
$(p_1,v_1)\Gamma^*(p_1,p_k)\Gamma^*(p_k,x)$. Otherwise if $d(v_1, p_k)$ is odd and 
$d(p_1,x)$ is even, then $(p_1,v_1)\Gamma^*(p_k,p_1)\Gamma^*(p_k,x)$. 
In either case, the second arc is $(x,p_k)$. Otherwise, $k + d(p_1, x)$ is odd and 
the second arc is $(p_k, x)$. So, $\overline{U(X)}$ is Figure~\ref{onec4only}(xi) 
and $X$ contains the dotted arcs.
\qed

\begin{figure}[H]
	\centering
	\captionsetup[subfigure]{labelformat=empty}
	\begin{subfigure}[b]{\textwidth}
		\centering
		\begin{tikzpicture}
		\node[unlabelled] (1) at (-1.7071,0.7071) {};
		\node[unlabelled] (2) at (-2.4142,0) {};
		\node[unlabelled] (3) at (-1.7071,-0.7071) {};
		\node[unlabelled] (x) at (-1.7071,-1.7071) {};
		\node[unlabelled] (p1) at (-1,0) {};
		\node (p) at (0,0) {$\ldots$};
		\node[unlabelled] (pk) at (1,0) {};
		\node[unlabelled] (4) at (1.7071,0.7071) {};
		\node[unlabelled] (5) at (2.4142,0) {};
		\node[unlabelled] (6) at (1.7071,-0.7071) {};
		\node[unlabelled] (y) at (1.7071,-1.7071) {};
		\node (l1) at (-1.7071,1.0071) {$v_1$};
		\node (l2) at (-2.7142,0) {$v_2$};
		\node (l3) at (-1.4071,-1.0071) {$v_3$};
		\node (lx) at (-1.7071,-2.0071) {$x$};
		\node (l4) at (1.7071,1.0071) {$v_4$};
		\node (l5) at (2.7142,0) {$v_5$};
		\node (l6) at (1.4071,-1.0071) {$v_6$};
		\node (ly) at (1.7071,-2.0071) {$y$};
		\node (lp1) at (-0.7,-0.3) {$p_1$};
		\node (lpk) at (0.7,-0.3) {$p_k$};
		\draw (1) -- (2) (1) -- (p1) (2) -- (3) (3) -- (x) (3) -- (p1) (p1) -- (p) (p) -- (pk) (4) -- (pk) (4) -- (5) (pk) -- (6) (5) -- (6) (6) -- (y);
		\draw[dashed, ->-=0.5] (2) to (x);
		\draw[dashed, ->-=0.5] (5) to (y);
		\end{tikzpicture}
		\subcaption{(i)}
	\end{subfigure}
	
	\begin{subfigure}[b]{.5\textwidth}
		\centering
		\begin{tikzpicture}
		\node[unlabelled] (1) at (0,0) {};
		\node[unlabelled] (4) at (1,0) {};
		\node[unlabelled] (5) at (2,0) {};
		\node[unlabelled] (2) at (0,1) {};
		\node[unlabelled] (3) at (1,1) {};
		\node[unlabelled] (6) at (2,1) {};
		\node[unlabelled] (x) at (-1,0) {};
		\node[unlabelled] (y) at (3,1) {};
		\node (l1) at (0,-0.3) {$v_1$};
		\node (l4) at (1,-0.3) {$v_4$};
		\node (l5) at (2,-0.3) {$v_5$};
		\node (l2) at (0,1.3) {$v_2$};
		\node (l3) at (1,1.3) {$v_3$};
		\node (l6) at (2,1.3) {$v_6$};
		\node (lx) at (-1,-0.3) {$x$};
		\node (ly) at (3,1.3) {$y$};
		\draw (1) -- (2) (2) -- (3) (3) -- (4) (4) -- (1) (4) -- (5) (5) -- (6) (6) -- (3) (1) -- (x) (6) -- (y);
		\draw[dashed, ->-=0.5] (x) to (2);
		\draw[dashed, ->-=0.5] (y) to (5);
		\end{tikzpicture}
		\subcaption{(ii)}
	\end{subfigure}%
	\begin{subfigure}[b]{.5\textwidth}
		\centering
		\begin{tikzpicture}
		\node[unlabelled] (1) at (0,0) {};
		\node[unlabelled] (4) at (1,0) {};
		\node[unlabelled] (5) at (2,0) {};
		\node[unlabelled] (2) at (0,1) {};
		\node[unlabelled] (3) at (1,1) {};
		\node[unlabelled] (6) at (2,1) {};
		\node (l1) at (0,-0.3) {$v_1$};
		\node (l4) at (1,-0.3) {$v_4$};
		\node (l5) at (2,-0.3) {$v_5$};
		\node (l2) at (0,1.3) {$v_2$};
		\node (l3) at (1,1.3) {$v_3$};
		\node (l6) at (2,1.3) {$v_6$};
		\draw (1) -- (2) (2) -- (3) (3) -- (4) (4) -- (1) (4) -- (5) (5) -- (6) (6) -- (3);
		\draw[dashed, ->-=0.8] (1) to (3);
		\draw[dashed, ->-=0.8] (4) to (2);
		\end{tikzpicture}
		\subcaption{(iii)}
	\end{subfigure}

 	\begin{subfigure}[b]{.25\textwidth}
 		\centering
 		\begin{tikzpicture}
 		\node[unlabelled] (1) at (0,2) {};
 		\node[unlabelled] (3) at (0,1) {};
 		\node[unlabelled] (5) at (0,0) {};
 		\node[unlabelled] (2) at (-1,1.5) {};
 		\node[unlabelled] (4) at (-1,0.5) {};
 		\node[unlabelled] (x) at (-2,0.5) {};
 		\node (p) at (.75,1) {$\ldots$};
 		\node[unlabelled] (pk) at (1.5,1) {};
 		\node (l2) at (-1.3, 1.5) {$v_2$};
 		\node (l4) at (-1.3, 0.8) {$v_4$};
 		\node (lx) at (-2, 0.2) {$u_4$};
 		\node (l5) at (0, -0.3) {$v_5$};
 		\node (l1) at (0, 2.3) {$v_1$};
 		\node (l3) at (0, 0.7) {$p_1$};
 		\node (lk) at (1.5,0.7) {$p_k$};
 		\draw (2) -- (1) (2) -- (3) (2) -- (5) (4) -- (1) (4) -- (3) (4) -- (5) (4) -- (x) (3) -- (p) (p) -- (pk);
 		\draw[dashed, ->-=0.5] (x) to (5);
 		\draw[dashed, ->-=0.5] (1) to (pk);
 		\end{tikzpicture}
 		\subcaption{(iv): $k$ is even.}
 	\end{subfigure}%
 	\begin{subfigure}[b]{.25\textwidth}
 		\centering
 		\begin{tikzpicture}
 		\node[unlabelled] (1) at (0,2) {};
 		\node[unlabelled] (3) at (0,1) {};
 		\node[unlabelled] (5) at (0,0) {};
 		\node[unlabelled] (2) at (-1,1.5) {};
 		\node[unlabelled] (4) at (-1,0.5) {};
 		\node[unlabelled] (x) at (-2,0.5) {};
 		\node (p) at (.75,1) {$\ldots$};
 		\node[unlabelled] (pk) at (1.5,1) {};
 		\node (l2) at (-1.3, 1.5) {$v_2$};
 		\node (l4) at (-1.3, 0.8) {$v_4$};
 		\node (lx) at (-2, 0.2) {$u_4$};
 		\node (l5) at (0, -0.3) {$v_5$};
 		\node (l1) at (0, 2.3) {$v_1$};
 		\node (l3) at (0, 0.7) {$p_1$};
 		\node (lk) at (1.5,0.7) {$p_k$};
 		\draw (2) -- (1) (2) -- (3) (2) -- (5) (4) -- (1) (4) -- (3) (4) -- (5) (4) -- (x) (3) -- (p) (p) -- (pk);
 		\draw[dashed, ->-=0.5] (x) to (5);
 		\draw[dashed, ->-=0.5] (pk) to (1);
 		\end{tikzpicture}
 		\subcaption{(v): $k \geq 3$ is odd.}
 	\end{subfigure}%
 	\begin{subfigure}[b]{.25\textwidth}
 		\centering
 		\begin{tikzpicture}
 		\node[unlabelled] (1) at (0,2) {};
 		\node[unlabelled] (3) at (0,1) {};
 		\node[unlabelled] (5) at (0,0) {};
 		\node[unlabelled] (2) at (-1,1.5) {};
 		\node[unlabelled] (4) at (-1,0.5) {};
 		\node[unlabelled] (x) at (-2,0.5) {};
 		\node[unlabelled] (y) at (1,1) {};
 		\node[unlabelled] (z) at (1,2) {};
 		\node (l2) at (-1.3, 1.5) {$v_2$};
 		\node (l4) at (-1.3, 0.8) {$v_4$};
 		\node (lx) at (-2, 0.2) {$u_4$};
 		\node (ly) at (1.3,1) {$u_3$};
 		\node (lz) at (1.3,2) {$z$};
 		\node (l5) at (0, -0.3) {$v_5$};
 		\node (l1) at (0, 2.3) {$v_1$};
 		\node (l3) at (0, 0.7) {$v_3$};
 		\draw (2) -- (1) (2) -- (3) (2) -- (5) (4) -- (1) (4) -- (3) (4) -- (5) (4) -- (x) (3) -- (z) (1) -- (z) (3) -- (y);
 		\draw[dashed, ->-=0.5] (x) to (5);
 		\draw[dashed, ->-=0.5] (y) to (z);
 		\end{tikzpicture}
 		\subcaption{(vi)}
 	\end{subfigure}%
 	\begin{subfigure}[b]{.25\textwidth}
 		\centering
 		\begin{tikzpicture}
 		\node[unlabelled] (1) at (0,2) {};
 		\node[unlabelled] (3) at (0,1) {};
 		\node[unlabelled] (5) at (0,0) {};
 		\node[unlabelled] (2) at (-1,1.5) {};
 		\node[unlabelled] (4) at (-1,0.5) {};
 		\node[unlabelled] (x) at (-2,0.5) {};
 		\node[unlabelled] (y) at (1,1) {};
 		\node (l2) at (-1.3, 1.5) {$v_2$};
 		\node (l4) at (-1.3, 0.8) {$v_4$};
 		\node (lx) at (-2, 0.2) {$u_4$};
 		\node (l5) at (0, -0.3) {$v_5$};
 		\node (l1) at (0, 2.3) {$v_1$};
 		\node (l3) at (0, 0.7) {$v_3$};
 		\node (ly) at (1.3,1) {$u_3$};
 		\draw (2) -- (1) (2) -- (3) (2) -- (5) (4) -- (1) (4) -- (3) (4) -- (5) (4) -- (x) (3) -- (y);
 		\draw[dashed, ->-=0.5] (x) to (5);
 		\draw[dashed, ->-=0.5] (3) to (1);
 		\end{tikzpicture}
 		\subcaption{(vii)}
 	\end{subfigure}%
\caption{Obstructions $X$ for which $\overline{U(X)}$ contains two induced $C_4$'s
      but no $C_3$. \label{twoc4}}
\end{figure}
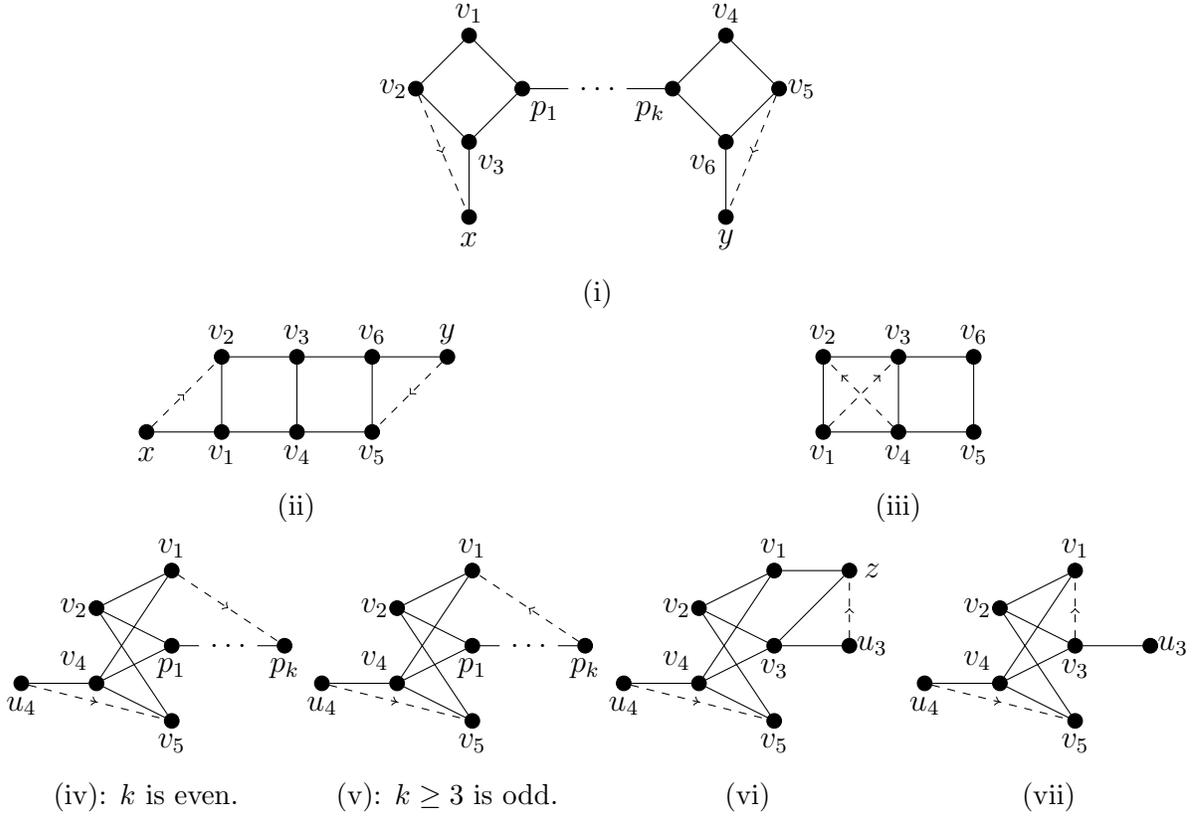

\begin{theorem} \label{twoc4s}
Let $X$ be an obstruction that has no cut-vertices. Suppose that $\overline{U(X)}$ 
contains two induced $C_4$'s but no $C_3$. Then $\overline{U(X)}$ is one of 
the graphs in Figure~\ref{twoc4} and $X$ or its dual contains the dotted arcs.
\end{theorem}
\pf Suppose there are two induced $C_4$'s in $\overline{U(X)}$ which share at most
one common vertex. Let $C$ and $C'$ be such induced $C_4$'s and let
$P: p_1,p_2, \dots, p_k$ be a shortest path between a vertex of $C$ and a vertex of
$C'$. By Lemma \ref{non-cut-vertices} any connected subgraph of $\overline{U(X)}$
has at most six non-cut-vertices. The (connected) subgraph of $\overline{U(X)}$
induced by $C \cup C' \cup P$ has at least six non-cut-vertices and thus has
exactly six non-cut-vertices. This implies that $P$ is the unique path between $C$
and $C'$ and each $p_i$ of $P$ is a cut-vertex of $\overline{U(X)}$.
Since the subgraph of $\overline{U(X)}$ induced by $C \cup C' \cup P$ has six
non-cut-vertices, $\overline{U(X)}$ also has six non-cut-vertices according to
Lemma \ref{non-cut-vertices}. Thus by Corollary \ref{at most 6 non-cut-vertices}
the six non-cut-vertices of $\overline{U(X)}$ form two disjoint arc-balancing
triples.

Denote $C: v_1v_2v_3p_1$ and $C': v_4v_5v_6p_k$. We first show that $v_1$ is not
incident with an arc. Suppose there is an arc between $v_1$ and a vertex $z$.
Since $p_1$ is a cut-vertex of $\overline{U(X)}$, it does not balance the arc
between $v_1$ and $z$. Since $p_1$ adjacent to $v_1$, it is adjacent to  $z$.
If $z$ is not in $C \cup C' \cup P$, then the subgraph induced by
$C \cup C' \cup P \cup \{z\}$ contains seven non-cut-vertices (i.e.,
$v_1,v_2,v_3,v_4,v_5,v_6,z$), which contradicts Lemma \ref{non-cut-vertices}.
So $z$ is in $C \cup C' \cup P$. Note that $z$ is adjacent to $p_1$.
If $z \neq v_3$, then $v_2$ is adjacent to $v_1$ but not $z$ and there is
another vertex in $C' \cup P$ adjacent to $z$ but not $v_1$, a contradiction to
the fact that $v_1$ and $z$ are in a an arc-balancing triple.
Thus $z = v_3$. But then the vertex $v$ which balances the arc between
$v_1$ and $z$ can not be in $C \cup C' \cup P$. Assume without loss
of generality that $v$ is adjacent to $v_1$ but not to $z$. Since $v_1$ is incident
with an arc, it is not a cut-vertex of $\overline{U(X)}$. So $\overline{U(X)}-v_1$
has a $(v,v_3)$-path $Q$. The connected subgraph of $\overline{U(X)}$ induced by
$C \cup C' \cup P \cup Q$ contains seven non-cut-vertices
(i.e., $v_1,v_2,v_3,v_4,v_5,v_6,v$), a contradiction to
Lemma \ref{non-cut-vertices}. Therefore $v_1$ is not incident with an arc.
By symmetry, none of $v_3,v_4,v_6$ is incident with an arc.

Since $p_1$ is a cut-vertex and any induced $C_4$ in $\overline{U(X)}$ contains 
at most two cut-vertices of $\overline{U(X)}$, $v_1, v_3$ cannot both be 
cut-vertices of $\overline{U(X)}$. Moreover, we know from above that neither of
$v_1,v_3$ is incident with an arc so Proposition \ref{vertex classification} 
implies that one of $v_1, v_3$ is an arc-balancing vertex. Similarly, one of 
$v_4,v_6$ is an arc-balancing vertex. Hence one of $v_1,v_3$ is an arc-balancing
vertex and the other is a cut-vertex of $\overline{U(X)}$. Without loss of 
generality, assume $v_1$ is an arc-balancing vertex and $v_3$ is a cut-vertex of 
$\overline{U(X)}$. The vertex $v_1$ is adjacent to exactly one endvertex $u$ of 
the arc it balances. We claim that $u = v_2$. 
If $u$ is not in $C \cup C' \cup P$, then $\overline{U(X)}$ contains a 
$(u,v_2)$-path $Q$ not containing $v_1$ because $v_1$ is a non-cut-vertex.
So the connected subgraph induced by $C \cup C' \cup P \cup Q$ contains 
seven non-cut-vertices, contradicting Lemma \ref{non-cut-vertices}. So $u$ is in 
$C \cup C' \cup P$. Since $p_1$ is a cut-vertex and no cut-vertex of 
$\overline{U(X)}$ is arc-balancing or incident with an arc, $u \neq p_1$. Hence 
$u = v_2$ as claimed, and $v_1$ balances an arc incident with $v_2$. 
The other endvertex $x$ must therefore be outside of $C \cup C' \cup P$ and 
adjacent to $v_3$.

By symmetry, $v_4$ balances an arc between $v_5$ and a vertex $y$ outside of 
$C \cup C' \cup P$ and adjacent to $v_6$. By possibly taking the dual
of $X$ assume $(v_2,x)$ is an arc in $X$. If $k$ is even, then 
Proposition \ref{complement gamma} implies 
$(p_1,x)\Gamma^*(x,p_k)$ and $(p_1,y)\Gamma^*(y,p_k)$. Hence 
$$(v_2,x)\Gamma(x,v_1)\Gamma(p_1,x)\Gamma^*(x,p_k)\Gamma(v_6,x)\Gamma(x,y)
\Gamma(y,v_3)\Gamma(p_1,y)\Gamma^*(y,p_k)\Gamma(v_4,y)\Gamma(y,v_5)$$ 
and so the second arc is $(v_5,y)$. Otherwise, $k$ is odd and in this case, we have
$$(v_2,x)\Gamma(x,v_1)\Gamma(p_1,x)\Gamma^*(p_k,x)\Gamma(x,v_6)\Gamma(y,x)
\Gamma(v_3,y)\Gamma(y,p_1)\Gamma^*(y,p_k)\Gamma(v_4,y)\Gamma(y,v_5)$$
and the second arc is again $(v_5,y)$. The minimality of $X$ ensure that
$\overline{U(X)}$ is the graph in Figure~\ref{twoc4}(i) and $X$ contains 
the dotted arcs.

Suppose next there are two induced $C_4$'s in $\overline{U(X)}$ which share 
two common vertices but no two induced $C_4$'s in $\overline{U(X)}$ share three
common vertices. Then such two $C_4$'s must share an edge. Let 
$C: v_1v_2v_3v_4$ and $C': v_3v_4v_5v_6$ be such induced $C_4$'s in 
$\overline{U(X)}$. Since $\overline{U(X)}$ contains no $C_3$ and no two induced 
$C_4$'s share three vertices, the subgraph induced by $C \cup C'$ has exactly 
seven edges belonging to the two $C_4$'s. We claim that no vertex outside of 
$C \cup C'$ is adjacent to $v_3$ or $v_4$. Indeed, if some vertex $z$ outside of 
$C \cup C'$ is adjacent to $v_3$ or $v_4$, then it must be adjacent to at least two
vertices in $C \cup C'$, because otherwise Theorem \ref{Tucker} would imply that 
$U(X)$ is not a proper circular-arc graph. But then $C \cup C' \cup \{z\}$ would
induced a connected subgraph in $\overline{U(X)}$ having seven non-cut-vertices, 
a contradiction to Lemma \ref{non-cut-vertices}. 

Suppose both $v_3$ and $v_4$ are arc-balancing vertices. If $v_3$ balances an arc 
incident with $v_4$, then the other endvertex $z$ of the arc is not in 
$C \cup C'$ that is adjacent to both $v_1$ and $v_5$. Thus the (connected) subgraph
of $\overline{U(X)}$ induced by $C \cup C' \cup \{z\}$ has seven non-cut-vertices, 
contradicting Lemma \ref{non-cut-vertices}. Hence $v_3$ does not balance an arc 
incident with $v_4$. Moreover, since no vertices outside of $C \cup C'$ is adjacent
to $v_3$, $v_3$ must balance an arc incident with $v_2$ or $v_6$. Similarly, $v_4$ 
must balance an arc incident with $v_1$ or $v_5$. Without loss of generality, 
assume $v_3$ balances an arc incident with $v_2$. Thus the other endvertex $x$ is 
a vertex whose only neighbour in $C \cup C'$ is $v_1$. We claim $v_4$ balances 
an arc incident with $v_5$. Otherwise, $v_4$ balances an arc incident with $v_1$, 
so the other endvertex of the arc has $v_2$ as the only neighbour in $C \cup C'$. 
Clearly, either $v_5$ is a non-cut-vertex or $\overline{U(X)}$ has a non-cut-vertex
in a component of $\overline{U(X-v_5)}$ not containing vertices in $C \cup C'$. 
In either case, $\overline{U(X)}$ contains a non-cut-vertex that is neither 
an endvertex of an arc nor an arc-balancing vertex, a contradiction by 
Proposition \ref{vertex classification}. Hence $v_4$ balances an arc incident with 
$v_5$ as claimed. The other endvertex $y$ of the arc has $v_6$ as the only 
neighbour in $C \cup C'$. By possibly considering the dual of $X$, assume 
$(x,v_2)$ is an arc in $X$. Since $(x,v_2)\Gamma(v_3,x)\Gamma(x,v_6)
\Gamma(y,x)\Gamma(v_1,y)\Gamma(y,v_4) \Gamma(v_5,y)$, the second arc is $(y,v_5)$. 
The minimality of $X$ ensures that $\overline{U(X)}$ is Figure~\ref{twoc4}(ii)
and $X$ contains the dotted arcs.

Suppose at least one of $v_3, v_4$ is not an arc-balancing vertex. Without loss of 
generality assume that $v_3$ is not an arc-balancing vertex. Then by Proposition 
\ref{vertex classification}, $v_3$ must be incident with an arc. The subgraph of 
$\overline{U(X)}$ induced by $C \cup C'$ has six non-cut-vertices so by Lemma 
\ref{non-cut-vertices} $\overline{U(X)}$ has six non-cut-vertices. Corollary 
\ref{at most 6 non-cut-vertices} implies that the six non-cut-vertices of 
$\overline{U(X)}$ form two disjoint arc-balancing triples. Since $v_3$ has three 
neighbours in $C \cup C'$ and the arc incident with $v_3$ has an arc-balancing 
vertex, the other endvertex must be adjacent to at least two of the three 
neighbours of $v_3$ in $C \cup C'$. Since any connected subgraph of 
$\overline{U(X)}$ has at most six non-cut-vertices by Lemma \ref{non-cut-vertices},
so any vertex not in $C \cup C'$ is adjacent to at most one vertex in $C \cup C'$. 
It follows that the other endvertex of the arc incident with $v_3$ is in 
$C \cup C'$. Without loss of generality, assume $(v_1,v_3)$ is an arc in $X$. 
Clearly, $v_6$ is the $(v_1,v_3)$-balancing vertex. We claim that $v_4$ is incident
with an arc. Otherwise, Proposition \ref{vertex classification} would imply $v_4$ 
is an arc-balancing vertex and hence balances an arc incident with $v_5$. 
By the above, the other endvertex $x$ of the arc incident with $v_5$ is adjacent to
$v_6$. Clearly, either $v_2$ is a non-cut-vertex or $\overline{U(X)}$ has 
a non-cut-vertex in a component of $\overline{U(X-v_2)}$ not containing vertices 
in $C \cup C'$. In either case, $\overline{U(X)}$ contains a non-cut-vertex that 
does not belong to either arc-balancing triple, a contradiction by 
Proposition \ref{vertex classification}. Therefore $v_4$ must be incident with 
an arc. By a similar argument as above, the other endvertex of the arc incident 
with $v_4$ is in $C \cup C'$. Since $\overline{U(X)}$ contains two disjoint 
arc-balancing triples and $v_6$ is the $(v_1,v_3)$-balancing vertex, the other 
endvertex cannot be $v_6$ and hence must be $v_2$. Since 
$(v_1,v_3)\Gamma(v_6,v_1)\Gamma(v_2,v_6)\Gamma(v_5,v_2)\Gamma(v_2,v_4)$, the second
arc is $(v_4,v_2)$. The minimality of $X$ ensures that $\overline{U(X)}$ is 
Figure~\ref{twoc4}(iii) and $X$ contains the dotted arcs.

Suppose now that there are two induced $C_4$'s in $\overline{U(X)}$ which share 
three common vertices. Let $C: v_1v_2v_3v_4$ and $C': v_2v_3v_4v_5$ be such induced 
$C_4$'s. Note that $C \cup C'$ induces a $K_{2,3}$ in $\overline{U(X)}$. 
Each $v_i$ with $1 \leq i \leq 5$ may or may not be a cut-vertex of 
$\overline{U(X)}$. If $v_i$ is a cut-vertex of $\overline{U(X)}$, then 
$\overline{U(X-v_i)}$ must contain a non-cut-vertex of $\overline{U(X)}$ that 
is not in $C \cup C'$. Let $u_i$ be such a vertex in $\overline{U(X-v_i)}$ when 
$v_i$ is a cut-vertex; otherwise let $u_i = v_i$ for each $1 \leq i \leq 5$. 
First note that $u_2, u_4$ are non-adjacent and that $u_1,u_3,u_5$ are pairwise
non-adjacent. Moreover, if $u_i \neq v_i$ then $u_i$ is not adjacent to $u_j$
for all $j \neq i$. Since each $u_i$ is a non-cut-vertex, it is an endvertex
of an arc or an arc-balancing vertex by Proposition \ref{vertex classification}. 
This implies that there is an arc-balancing triple $T$ contained in 
$\{u_1,u_2,\dots,u_5\}$. Since there is exactly one edge in $T$, we know from 
the above observation the only edge in $T$ has one endvertex in $\{u_2,u_4\}$ and 
the other in $\{u_1,u_3,u_5\}$. Without loss of generality assume that $u_2u_5$ is 
the edge in $T$. Then we must have $u_2 = v_2$ and $u_5 = v_5$ and thus
neither $v_2$ nor $v_5$ is a cut-vertex of $\overline{U(X)}$. It is easy to see
that the third vertex of $T$ is $u_4$ and $v_2$ balances the arc between $u_4$ and
$v_5$. Without loss of generality assume $(u_4,v_5)$ is an arc in $X$.
Clearly, $u_4 \neq v_4$ and $v_4$ is a cut-vertex of $\overline{U(X)}$. 

Since $C$ has at most two cut-vertices of $\overline{U(X)}$, at most one of 
$v_1,v_3$ can be a cut-vertex. If neither of $v_1,v_3$ is a cut-vertex, then one
of them is an endvertex of an arc which is balanced by the other vertex. Since
$v_2$ is adjacent to both $v_1,v_3$, it is adjacent to the other endvertex of 
the arc, which implies $U(X)$ contains a $C_3$, a contradiction to the assumption.  
So exactly one of $v_1,v_3$ is a cut-vertex of $\overline{U(X)}$ and we assume 
it is $v_3$. Suppose that there is an arc between $v_1$ and $u_3$.  Let 
$P: v_3=p_1, \dots, p_k=u_3$ be the shortest $(v_3,u_3)$-path in $\overline{U(X)}$. 
Suppose $k$ is even. Then by Proposition \ref{complement gamma}, 
$$(u_4,v_5)\Gamma(v_2,u_4)\Gamma(u_4,p_1)\Gamma^*(p_k,u_4)\Gamma(v_4, p_k)
\Gamma(p_k,v_5)\Gamma(v_2, p_k)\Gamma(p_k,v_1)$$
and hence the second arc is $(v_1,p_k)$. The minimality of $X$ ensures 
$\overline{U(X)}$ is Figure~\ref{twoc4}(iv) and $X$ contains the dotted arcs.
Otherwise, $k$ is odd and $\overline{U(X)}$ is Figure~\ref{twoc4}(v) and 
$X$ contains the dotted arcs. 

Suppose that there is no arc between $v_1$ and $u_3$. Then either $u_3$ is incident
with an arc balanced by $v_1$ or $v_1$ is incident with an arc balanced by $u_3$.  
Suppose it is the the former. Let $z$ denote the other endvertex. Since $v_1$ is 
not adjacent to $u_3$, it is adjacent to $z$. Moreover, since $u_3$ is not adjacent
to $v_5$, $z$ is also not adjacent to $v_5$. In particular, $z \notin C \cup C'$. 
If a vertex other than $v_1$ is adjacent to $z$, then it must also be adjacent to
$u_3$. In particular, the choice of $u_3$ implies that $v_3$ is the only vertex 
that is possibly adjacent to $u_3$. Since $\overline{U(X)}$ is connected, $v_3$ 
must be adjacent to $u_3$, so $v_3$ is adjacent to $z$ as well. Since 
$(u_4,v_5)\Gamma(v_2,u_4)\Gamma(u_4, v_3)\Gamma(u_3,u_4)\Gamma(v_4,u_3)
\Gamma(u_3,v_1)\Gamma(z,u_3)$, the second arc is $(u_3, z)$. The minimality of $X$ 
ensures $\overline{U(X)}$ is Figure~\ref{twoc4}(vi) and $X$ contains 
the dotted arcs.

Suppose instead that $v_1$ is incident with an arc balanced by $u_3$. Since $u_3$ 
is not adjacent to $v_1$, it is adjacent to the other endvertex. Moreover, since 
$v_1$ is adjacent to $v_2$ and $v_4$, the other endvertex is also adjacent to $v_2$
and $v_4$. By the choice of $u_3$, the only vertex that can be adjacent to all of 
$u_3, v_2, v_4$ is $v_3$, so the other endvertex is $v_3$ and $u_3$ is adjacent to 
$v_3$. Since $(u_4,v_5)\Gamma(v_2,u_4)\Gamma(u_4,v_3)\Gamma(u_3,u_4)
\Gamma(v_4,u_3)\Gamma(u_3,v_1)\Gamma(v_1,v_3)$, the second arc is $(v_3,v_1)$. 
The minimality of $X$ ensures $\overline{U(X)}$ is Figure~\ref{twoc4}(vii) 
and $X$ contains the dotted arcs.
\qed

\subsection{$\overline{U(X)}$ contains a $C_3$ and an induced $C_4$}

\begin{lemma} \label{c3c4lemma}
Let $X$ be an obstruction which has no cut-vertices. Suppose $\overline{U(X)}$ 
contains a $C_3$ and an induced $C_4$. Then $\overline{U(X)}$ contains a unique 
$C_3$ and a unique induced $C_4$, which share two common vertices. Moreover, each 
vertex not in any of the cycles is a leaf adjacent to a vertex on the $C_3$ and 
is incident with an arc.
\end{lemma}  
\pf By Lemmas \ref{odd cycles} and \ref{onec3}, $\overline{U(X)}$ contains 
a unique $C_3$ and each vertex not on the $C_3$ is adjacent to a vertex in
the $C_3$. It follows that each vertex not on the $C_3$ is adjacent to exactly one 
vertex on the $C_3$. We show that if $C$ is an induced $C_4$ in $\overline{U(X)}$ 
then $C$ shares exactly two vertices with the $C_3$. Clearly, $C$ share at most two
vertices with the $C_3$. The fact that every vertex not on the $C_3$ is adjacent 
to a vertex in the $C_3$ implies that $C$ cannot share exactly one vertex with 
the $C_3$. If $C$ shares no vertex with the $C_3$, then $C \cup C_3$ induces 
a connected subgraph in $\overline{U(X)}$ with seven non-cut-vertices, 
a contradiction to Lemma \ref{non-cut-vertices}. 

Denote the unique $C_3$ in $\overline{U(X)}$ by $v_1v_2v_3$ and without loss of 
generality assume that $v_2v_3v_4v_5$ is an induced $C_4$ in the graph.   
Let $u \notin \{v_1, v_2, \dots, v_5\}$. From the above we know that 
$u$ is adjacent to exactly one vertex in the $C_3$. Suppose that $u$ is adjacent to
$v_1$. Then $u$ cannot be adjacent to both $v_4,v_5$ as otherwise $uv_4v_5$ is 
another $C_3$ in $\overline{U(X)}$, a contradiction. If $u$ is adjacent to one of 
$v_4,v_5$ then $uv_1v_2v_5v_4$ or $uv_1v_3v_4v_5$ is an induced $C_5$ in 
$\overline{U(X)}$, which contradicts Lemma \ref{C5lemma}. If $u$ is adjacent to 
a vertex $w \notin \{v_1, v_2, \dots, v_5\}$, then $w$ is not adjacent to $v_1$ 
due to the uniqueness of the $C_3$ and so is adjacent to $v_2$ or $v_3$. But then 
$\{u,w,v_1, v_2, \dots, v_5\}$ induces a connected subgraph of $\overline{U(X)}$
with seven non-cut-vertices, a contradiction to Lemma \ref{non-cut-vertices}.
Hence $u$ is a leaf in $\overline{U(X)}$. Suppose that $u$ is not adjacent to 
$v_1$. Then it is adjacent to $v_2$ or $v_3$. By symmetry we assume $u$ is adjacent
to $v_2$. It is not adjacent to $v_5$ as otherwise $uv_2v_5$ is another $C_3$ 
in $\overline{U(X)}$. It is not adjacent to $v_4$ as otherwise $uv_4v_5v_2$
is an induced $C_4$ which share just one vertex (namely, $v_2$) with the $C_3$.
Suppose that $u$ is adjacent to a vertex $w \notin \{v_1, v_2, \dots, v_5\}$. Then
$w$ is not adjacent to $v_2$ due to the uniqueness of the $C_3$. It is not adjacent
to $v_1$ because from the above any such vertex is a leaf. So $w$ is adjacent to
$v_3$. But then $\{u,w,v_1, v_2, \dots, v_5\}$ induces a connected subgraph of 
$\overline{U(X)}$ with seven non-cut-vertices, a contradiction to 
Lemma \ref{non-cut-vertices}. Therefore any vertex 
$u \notin \{v_1, v_2, \dots, v_5\}$ is a leaf adjacent to a vertex in the $C_3$.   
It follows that $v_2v_3v_4v_5$ is the unique induced $C_4$ in $\overline{U(X)}$.
If such a vertex $u$ is an arc-balancing vertex, then it balances an arc incident
with one of $v_1, v_2, v_3$, so the other endvertex is adjacent to the remaining
two vertices among $v_1, v_2, v_3$, contradicting the uniqueness of the $C_3$. 
Moreover, such a $u$ is not a cut-vertex of $U(X)$ or of $\overline{U(X)}$.
Therefore by Proposition \ref{vertex classification} $u$ is incident with an arc.  
\qed

\begin{theorem}\label{c3c4theorem}
Let $X$ be an obstruction that contains no cut-vertices. Suppose $\overline{U(X)}$ 
contains a $C_3$ and an induced $C_4$. Then $\overline{U(X)}$ is one of the graphs 
in Figure~\ref{c3c4figure} and $X$ or its dual contains the dotted arcs.
	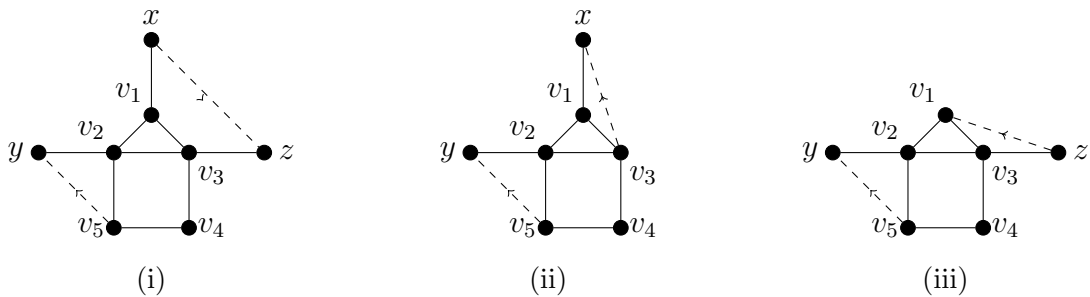
\begin{figure}[H]
		\centering
		\captionsetup[subfigure]{labelformat=empty}
		\begin{subfigure}[b]{.33\textwidth}
			\centering
			\begin{tikzpicture}
			\node[unlabelled]	(1) at (0.5, 1.5) {};
			\node[unlabelled]	(2) at (0, 1) {};
			\node[unlabelled]	(3) at (1, 1) {};
			\node[unlabelled]	(4) at (1, 0) {};
			\node[unlabelled]	(5) at (0, 0) {};
			\node[unlabelled]	(y) at (-1, 1) {};
			\node[unlabelled]	(x) at (0.5, 2.5) {};
			\node[unlabelled]	(z) at (2, 1) {};
			\node	(l1) at (0.2, 1.8) {$v_1$};
			\node	(l2) at (-0.3, 1.3) {$v_2$};
			\node	(l3) at (1.3, 0.7) {$v_3$};
			\node	(l5) at (-0.3, 0) {$v_5$};
			\node	(l4) at (1.3, 0) {$v_4$};
			\node	(ly) at (-1.3, 1) {$y$};
			\node	(lx) at (0.5, 2.8) {$x$};
			\node	(lz) at (2.3, 1) {$z$};
			\draw (1) -- (2) (1) -- (3) (2) -- (3) (3) -- (4) (4) -- (5) (5) -- (2) (2) -- (y) (3) -- (z) (1) -- (x);
			\draw[dashed, ->-=0.5] (5) to (y);
			\draw[dashed, ->-=0.5] (x) to (z);
			\end{tikzpicture}
			\subcaption{(i)}
		\end{subfigure}%
		\begin{subfigure}[b]{.33\textwidth}
			\centering
			\begin{tikzpicture}
			\node[unlabelled]	(1) at (0.5, 1.5) {};
			\node[unlabelled]	(2) at (0, 1) {};
			\node[unlabelled]	(3) at (1, 1) {};
			\node[unlabelled]	(4) at (1, 0) {};
			\node[unlabelled]	(5) at (0, 0) {};
			\node[unlabelled]	(y) at (-1, 1) {};
			\node[unlabelled]	(x) at (0.5, 2.5) {};
			\node	(l1) at (0.2, 1.8) {$v_1$};
			\node	(l2) at (-0.3, 1.3) {$v_2$};
			\node	(l3) at (1.3, 0.7) {$v_3$};
			\node	(l5) at (-0.3, 0) {$v_5$};
			\node	(l4) at (1.3, 0) {$v_4$};
			\node	(ly) at (-1.3, 1) {$y$};
			\node	(lx) at (0.5, 2.8) {$x$};
			\draw (1) -- (2) (1) -- (3) (2) -- (3) (3) -- (4) (4) -- (5) (5) -- (2) (2) -- (y) (1) -- (x);
			\draw[dashed, ->-=0.5] (5) to (y);
			\draw[dashed, ->-=0.5] (3) to (x);
			\end{tikzpicture}
			\subcaption{(ii)}
		\end{subfigure}%
		\begin{subfigure}[b]{.33\textwidth}
			\centering
			\begin{tikzpicture}
			\node[unlabelled]	(1) at (0.5, 1.5) {};
			\node[unlabelled]	(2) at (0, 1) {};
			\node[unlabelled]	(3) at (1, 1) {};
			\node[unlabelled]	(4) at (1, 0) {};
			\node[unlabelled]	(5) at (0, 0) {};
			\node[unlabelled]	(y) at (-1, 1) {};
			\node[unlabelled]	(z) at (2, 1) {};
			\node	(l1) at (0.2, 1.8) {$v_1$};
			\node	(l2) at (-0.3, 1.3) {$v_2$};
			\node	(l3) at (1.3, 0.7) {$v_3$};
			\node	(l5) at (-0.3, 0) {$v_5$};
			\node	(l4) at (1.3, 0) {$v_4$};
			\node	(ly) at (-1.3, 1) {$y$};
			\node	(lz) at (2.3, 1) {$z$};
			\draw (1) -- (2) (1) -- (3) (2) -- (3) (3) -- (4) (4) -- (5) (5) -- (2) (2) -- (y) (3) -- (z);
			\draw[dashed, ->-=0.5] (5) to (y);
			\draw[dashed, ->-=0.5] (z) to (1);
			\end{tikzpicture}
			\subcaption{(iii)}
		\end{subfigure}
\caption{Obstructions $X$ for which $\overline{U(X)}$ contains 
a $C_3$ and an induced $C_4$. \label{c3c4figure}}
	\end{figure}
\end{theorem}
\pf By Lemma \ref{c3c4lemma}, $\overline{U(X)}$ contains a unique $C_3$ and a unique
$C_4$ sharing two vertices. Denote the $C_3$ and the induced $C_4$ by 
$v_1v_2v_3$ and $v_2v_3v_4v_5$ respectively. We also know from the lemma that
every vertex not in any of the cycles is a leaf adjacent to a vertex in the $C_3$
and is incident with an arc. We claim that at least two vertices of the $C_3$ 
are neighbours of leaves. Indeed, since $U(X)$ does not contain cut-vertices and 
$v_1, v_2, \dots, v_5$ induce a path in $U(X)$, there must be at least one vertex 
not on the cycles, which implies at least one vertex on the $C_3$ is adjacent to
a leaf. If $v_2$ is the only vertex in the $C_3$ adjacent to a leaf. Then
$v_2$ is adjacent to every vertex except $v_4$ in $\overline{U(X)}$, 
a contradiction to Lemma \ref{2 non-adj}. So $v_2$ cannot be the only vertex in 
the $C_3$ adjacent to a leaf. By symmetry, $v_3$ cannot be the only vertex in 
the $C_3$ adjacent to a leaf. 

Suppose $v_1$ is the only vertex in the $C_3$ adjacent to a leaf. Let $u$ be 
a leaf of $\overline{U(X)}$ adjacent to $v_1$. By Lemma \ref{c3c4lemma} $u$ is 
incident with an arc. The other endvertex of this arc cannot be another leaf $v$ 
as otherwise $uv$ is a balanced edge in $U(X)$, a contradiction. 
So the other endvertex of this arc must be among $v_2, v_3, v_4, v_5$. 
Note first that none of $v_2, v_3, v_4, v_5$ is a cut-vertex.
If the arc is between $u$ and $v_2$, then $v_4$ does not balance this arc as 
it is adjacent to neither of the endvertices. So $v_4$ either balances or is 
incident with the second arc. If $v_4$ is incident with the second arc then $v_3$ 
and $v_5$ are arc-balancing vertices, which is not possible. If $v_4$ balances
the second arc, then the second arc must be incident with exactly one of 
$v_3, v_5$. But then $v_2$ is also a vertex adjacent to exactly one of the 
endvertices of the second arc, which is again impossible. This shows there is no
arc between $u$ and $v_2$. A similar argument shows that there is no arc between
$u$ and any of $v_3, v_4, v_5$. Therefore at least two vertices of the $C_3$
are neighbours of leaves.   

Suppose that all three vertices of the $C_3$ are neighbours of leaves. 
Let $x, y, z$ be leaves adjacent to $v_1, v_2, v_3$ respectively.
By Lemma \ref{c3c4lemma} each of $x, y, z$ is incident with an arc. So there is 
an arc between two of $x, y, z$. Suppose there is an arc between $y$ and $z$. 
By Lemma \ref{odd cycles}, neither of $v_4, v_5$ can be a cut-vertex
so each of them is an arc-balancing vertex or incident with an arc. 
At least one of $v_4, v_5$ must be an arc-balancing vertex because otherwise 
$x, y, z, v_4, v_5$ would be five vertices incident with arcs. Without loss of 
generality, assume $v_4$ is an arc-balancing vertex. Since $v_4$ is adjacent to 
neither of $y, z$, it cannot balance the arc between $y$ and $z$. Thus $v_4$ 
balances an arc between $x$ and $x_3$ or between $x$ and $v_5$. This is 
a contradiction because $v_2$ is another vertex adjacent to exactly one of 
endvertices of the arc balanced by $v_4$. Thus there is no arc between $y$ and $z$.
So there is an arc between $x$ and $y$ or between $x$ and $z$. By symmetry and 
taking the dual of $X$ if necessary we may assume that $(x,z)$ is an arc. 
Since $v_1$ and $v_3$ are two vertices adjacent to exactly one of $x, z$, there 
cannot be an $(x,z)$-balancing vertex. By Lemma \ref{odd cycles}, 
$v_4, v_5$ are non-cut-vertices of $\overline{U(X)}$ so each of them is 
arc-balancing or incident with an arc by Proposition \ref{vertex classification}. 
Clearly, none of them can be an $(x,z)$-balancing vertex. So either $v_4$ balances
the arc between $y$ and $v_5$ or $v_5$ balances the arc between $y$ and $v_4$.
The latter case is not possible because $v_3$ is adjacent to $v_4$ but not to $y$.
Hence there is an arc between $y$ and $v_5$. Since the two arcs of $X$ are 
opposing and 
$(x, z)\Gamma(z, v_1)\Gamma(v_2, z)\Gamma(z, y)\Gamma(y,v_3)\Gamma(v_4, y)
\Gamma(y, v_5)$, $(v_5,y)$ is an arc. The minimality of $X$ ensure that 
$\overline{U(X)}$ is the graph in Figure~\ref{c3c4figure}(i).

Suppose now that exactly two of $v_1, v_2, v_3$ are neighbours of leaves.
First consider the case when $v_1$ and $v_2$ are neighbours of leaves. 
Let $x, y$ be leaves adjacent to $v_1, v_2$ respectively. Clearly, $v_3$ is not 
a cut-vertex. Since $v_4, v_5$ are not the $C_3$, by Lemma \ref{odd cycles} 
they are not cut-vertices. Hence, each of $v_3, v_4, v_5$ is an arc-balancing
vertex or incident with an arc by Proposition \ref{vertex classification}. 
By Lemma \ref{c3c4lemma}, $x, y$ are incident with arcs. So at least one of
$v_3, v_4, v_5$ is an arc-balancing vertex.   
If $v_3$ is an arc-balancing vertex, then it must balance an arc incident with 
$v_2$ or $v_4$. But then $v_5$ is another vertex adjacent to exactly one endvertex
of this arc, a contradiction. Hence $v_3$ is not the arc-balancing vertex. 
For a similar reason, $v_5$ is also not an arc-balancing vertex. Thus $v_4$ is
an arc-balancing vertex. It is easy check an arc balanced by $v_4$ cannot be 
incident with $v_3$. So $v_4$ balances an arc incident with $v_5$. 
The other endvertex of this arc cannot be a leaf adjacent to $v_1$. Hence $v_4$
balances an arc between $v_5$ and a leaf adjacent to $v_2$. Without loss of 
generality assume it is between $v_5$ and $y$. By taking the dual of $X$ if 
necessary we may assume $(v_5, y)$ is an arc. Thus the second arc is between
$x$ and $v_3$.  Since 
$(v_5,y)\Gamma(y,v_4)\Gamma(v_4,v_2)\Gamma(v_1,v_4)\Gamma(v_4,x)\Gamma(x,v_3)$
and the two arcs of $X$ are opposing, the second arc is $(v_3, x)$. So 
$\overline{U(X)}$ is Figure~\ref{c3c4figure}(ii).

The case when $v_1$ and $v_3$ are neighbours of leaves is symmetric to the case 
when $v_1$ and $v_2$ are cut-vertices. So we now consider the case where $v_2$ and 
$v_3$ are neighbours of leaves. Let $y, z$ be leaves adjacent to $v_2, v_3$ 
respectively. By assumption, $v_1$ is not a cut-vertex. 
Since $v_4, v_5$ are not on the $C_3$, they are not cut-vertices by 
Lemma \ref{odd cycles}. Thus each of $v_1, v_4, v_5$ is an arc-balancing vertex or 
incident with an arc by Proposition \ref{vertex classification}. It follows that
at least one of them is an arc-balancing vertex. A similar proof as above shows 
that $v_4$ balances an arc between $v_5$ and $y$. Without loss of generality, 
assume $(v_5,y)$ is an arc. It is easy to see that $v_1$ is not an arc-balancing
vertex so it is incident with an arc. So the second arc is between $v_1$ and 
$z$. Since $(v_5,y)\Gamma(y,v_4)\Gamma(v_3,y)\Gamma(y,z)\Gamma(z,v_2)\Gamma(v_1,z)$
and the two arcs are opposing, the second arc is $(z,v_1)$. Hence 
$\overline{U(X)}$ is Figure~\ref{c3c4figure}(iii).
\qed

\subsection{$\overline{U(X)}$ contains an induced $C_5$}

\begin{lemma}\label{c5caselemma}
Let $X$ be an obstruction which has no cut-vertices. If $C$ is an induced $C_5$ 
in $\overline{U(X)}$, then the following statements hold:
\begin{enumerate}[label=(\alph*)]
\item $C$ is the unique cycle in $\overline{U(X)}$;
\item Each vertex not in $C$ is a leaf adjacent to a vertex in $C$ and 
      incident with an arc;
\item $C$ contains an arc-balancing vertex that is not a cut-vertex of 
      $\overline{U(X)}$;
\item If $v$ is an arc-balancing vertex in $C$, then $v$ balances an arc 
      between a neighbour of $v$ in $C$ and a leaf.
\end{enumerate}
\end{lemma}
\pf Let $C: v_1v_2v_3v_4v_5$ be an induced $C_5$ in $\overline{U(X)}$. 
By Lemma \ref{largest size of cycles}, $C$ is a longest induced cycle in 
$\overline{U(X)}$. According to Lemmas \ref{one C5} and \ref{C5lemma}, 
$\overline{U(X)}$ contains at most one induced $C_5$ but neither $C_3$ nor induced 
$C_4$. Thus $C$ is the unique cycle in $\overline{U(X)}$.

Suppose that $u$ is a vertex not in $C$. Then by Lemma \ref{odd cycles} $u$ is 
adjacent to a vertex in $C$. Since $C$ the unique cycle in $\overline{U(X)}$, $u$
is a leaf and hence not a cut-vertex of $\overline{U(X)}$. Let $v_i$ be the
neighbour of $u$. If $u$ is an arc-balancing vertex, then it balances an arc 
between $v_i$ and some vertex $w$. Since $v_i$ is adjacent to both 
$v_{i-1},v_{i+1}$ which do not balance the arc between $v_i$ and $w$, $w$ must be
adjacent to both $v_{i-1},v_{i+1}$. Thus $v_iv_{i+1}wv_{i-1}$ is a $C_4$,
a contradiction to the fact $C$ is the unique cycle in $\overline{U(X)}$.
Hence $u$ is incident with an arc by Proposition \ref{vertex classification}.

Suppose there are $k$ cut-vertices in $C$. We know from above that each such vertex
is adjacent to a leaf that is incident with an arc. Among the $5-k$ 
non-cut-vertices in $C$, at most $4-k$ can be incident with arcs because there are
at most four vertices incident with arcs. It follows that $C$ contains at least
one vertex that is not a cut-vertex of $\overline{U(X)}$ and not incident with an 
arc. Such a vertex must be an arc-balancing vertex by 
Proposition \ref{vertex classification}. 
Hence $C$ contains an arc-balancing vertex that is not a cut-vertex of 
$\overline{U(X)}$. Without loss of generality, assume $v_1$ is an arc-balancing
vertex and it balances an arc incident with $v_2$. Clearly the other endvertex 
cannot be in $C$ so it is a leaf of $\overline{U(X)}$.

Suppose that $v$ a vertex $C$ which balances an arc between $u$ and $w$. If one of 
$u,w$ is a leaf neighbour of $v$, then the other endvertex is an isolated vertex,
contradicting the fact that $\overline{U(X)}$ is connected. So neither of $u,w$ can
be a leaf neighbour of $v$. Since $v$ is adjacent to one of $u,w$, at least one of 
$u,w$ is in $C$. If the other vertex is also in $C$, then there is a vertex in $C$ 
which is not $v$ but is adjacent to exactly one of $u,w$, a contradiction. 
Therefore exactly one of $u,w$ is a neighbour of $v$ in $C$ and the other is 
a leaf. 
\qed

\begin{theorem}\label{c5theorem}
Let $X$ be an obstruction that has no cut-vertices. Suppose $\overline{U(X)}$ 
contains an induced $C_5$. Then $\overline{U(X)}$ is one of the graphs in 
Figure~\ref{c5} and $X$ or its dual contains the dotted arcs.
	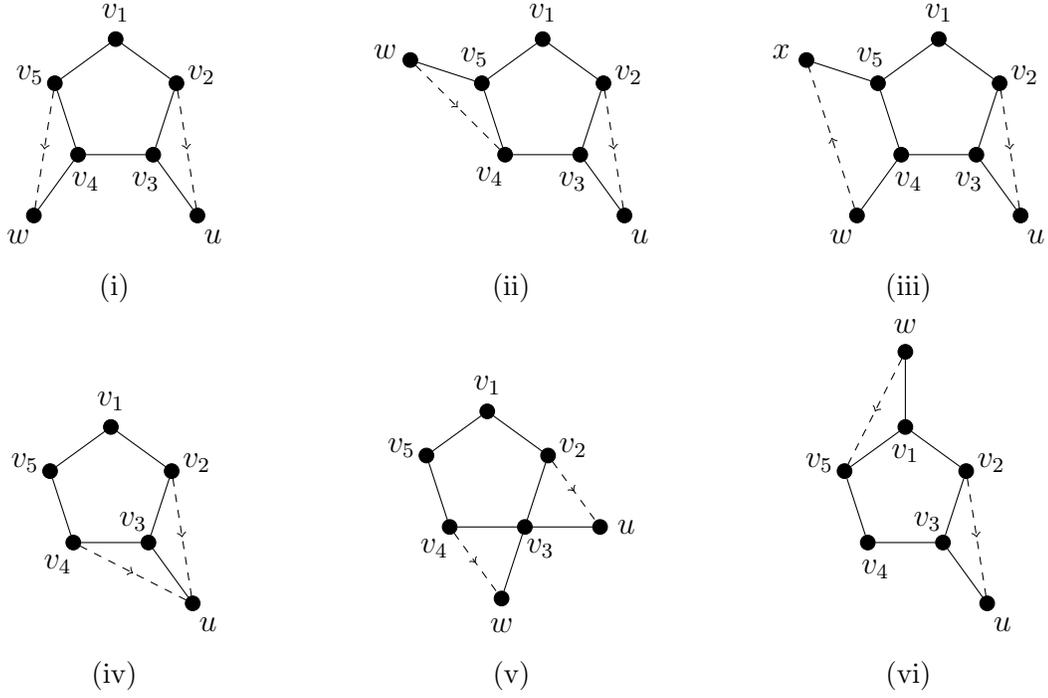
\begin{figure}[H]
		\centering
		\captionsetup[subfigure]{labelformat=empty}
		\begin{subfigure}[b]{.33\textwidth}
			\centering
			\begin{tikzpicture}
			\node[unlabelled]	(1) at (0.85*0, 0.85*1)	{};
			\node[unlabelled]	(5) at (0.85*-0.9511, 0.85*0.3090)	{};
			\node[unlabelled]	(4) at (0.85*-0.5878, 0.85*-0.8090)	{};
			\node[unlabelled]	(3) at (0.85*0.5878, 0.85*-0.8090)	{};
			\node[unlabelled]	(2) at (0.85*0.9511, 0.85*0.3090)	{};
			\node[unlabelled]	(w) at (1.85*-0.5878, 1.85*-0.8090)	{};
			\node[unlabelled]	(u) at (1.85*0.5878, 1.85*-0.8090)	{};
			\node	(l1) at (1.20*0, 1.20*1)	{$v_1$};
			\node	(l5) at (1.20*-0.9511, 1.20*0.3090)	{$v_5$};
			\node	(l4) at (0.85*-0.5878+0.1, 0.85*-0.8090-0.35)	{$v_4$};
			\node	(l3) at (0.85*0.5878-0.1, 0.85*-0.8090-0.35)	{$v_3$};
			\node	(l2) at (1.20*0.9511, 1.20*0.3090)	{$v_2$};
			\node	(lw) at (2.2*-0.5878, 2.2*-0.8090)	{$w$};
			\node	(lu) at (2.2*0.5878, 2.2*-0.8090)	{$u$};
			\draw	(1) -- (2)	(2) -- (3)	(3) -- (4)	(4) -- (5) (5) -- (1) (4) -- (w) (3) -- (u);
			\draw[dashed, ->-=0.5] (2) to (u);
			\draw[dashed, ->-=0.5] (5) to (w);
			\end{tikzpicture}
			\subcaption{(i)}
		\end{subfigure}%
		\begin{subfigure}[b]{.33\textwidth}
			\centering
			\begin{tikzpicture}
			\node[unlabelled]	(1) at (0.85*0, 0.85*1)	{};
			\node[unlabelled]	(5) at (0.85*-0.9511, 0.85*0.3090)	{};
			\node[unlabelled]	(4) at (0.85*-0.5878, 0.85*-0.8090)	{};
			\node[unlabelled]	(3) at (0.85*0.5878, 0.85*-0.8090)	{};
			\node[unlabelled]	(2) at (0.85*0.9511, 0.85*0.3090)	{};
			\node[unlabelled]	(u) at (1.85*0.5878, 1.85*-0.8090)	{};
			\node[unlabelled]	(w) at (1.85*-0.9511, 1.85*0.3090)	{};
			\node	(l1) at (1.20*0, 1.20*1)	{$v_1$};
			\node	(l5) at (0.85*-0.9511-0.1, 0.85*0.3090+0.35)	{$v_5$};
			\node	(l4) at (1.20*-0.5878, 1.20*-0.8090)	{$v_4$};
			\node	(l3) at (0.85*0.5878-0.1, 0.85*-0.8090-0.35)	{$v_3$};
			\node	(l2) at (1.20*0.9511, 1.20*0.3090)	{$v_2$};
			\node	(lu) at (2.20*0.5878, 2.20*-0.8090)	{$u$};
			\node	(lw) at (2.20*-0.9511, 2.20*0.3090)	{$w$};
			\draw	(1) -- (2)	(2) -- (3)	(3) -- (4)	(4) -- (5) (5) -- (1) (5) -- (w) (3) -- (u);
			\draw[dashed, ->-=0.5]	(w) to (4);
			\draw[dashed, ->-=0.5]	(2) to (u);
			\end{tikzpicture}
			\subcaption{(ii)}
		\end{subfigure}%
		\begin{subfigure}[b]{.33\textwidth}
			\centering
			\begin{tikzpicture}
			\node[unlabelled]	(1) at (0.85*0, 0.85*1)	{};
			\node[unlabelled]	(5) at (0.85*-0.9511, 0.85*0.3090)	{};
			\node[unlabelled]	(4) at (0.85*-0.5878, 0.85*-0.8090)	{};
			\node[unlabelled]	(3) at (0.85*0.5878, 0.85*-0.8090)	{};
			\node[unlabelled]	(2) at (0.85*0.9511, 0.85*0.3090)	{};
			\node[unlabelled]	(w) at (1.85*-0.5878, 1.85*-0.8090)	{};
			\node[unlabelled]	(u) at (1.85*0.5878, 1.85*-0.8090)	{};
			\node[unlabelled]	(x) at (1.85*-0.9511, 1.85*0.3090)	{};
			\node	(l1) at (1.20*0, 1.20*1)	{$v_1$};
			\node	(l5) at (0.85*-0.9511-0.1, 0.85*0.3090+0.35)	{$v_5$};
			\node	(l4) at (0.85*-0.5878+0.1, 0.85*-0.8090-0.35)	{$v_4$};
			\node	(l3) at (0.85*0.5878-0.1, 0.85*-0.8090-0.35)	{$v_3$};
			\node	(l2) at (1.20*0.9511, 1.20*0.3090)	{$v_2$};
			\node	(lw) at (2.2*-0.5878, 2.2*-0.8090)	{$w$};
			\node	(lu) at (2.2*0.5878, 2.2*-0.8090)	{$u$};
			\node	(lx) at (2.2*-0.9511, 2.2*0.3090)	{$x$};
			\draw	(1) -- (2)	(2) -- (3)	(3) -- (4)	(4) -- (5) (5) -- (1) (4) -- (w) (3) -- (u) (5) -- (x);
			\draw[dashed, ->-=0.5]	(w) to (x);
			\draw[dashed, ->-=0.5]	(2) to (u);
			\end{tikzpicture}
			\subcaption{(iii)}
		\end{subfigure}
		
		\begin{subfigure}[b]{.33\textwidth}
			\centering
			\begin{tikzpicture}
			\node[unlabelled]	(1) at (0.85*0, 0.85*1)	{};
			\node[unlabelled]	(5) at (0.85*-0.9511, 0.85*0.3090)	{};
			\node[unlabelled]	(4) at (0.85*-0.5878, 0.85*-0.8090)	{};
			\node[unlabelled]	(3) at (0.85*0.5878, 0.85*-0.8090)	{};
			\node[unlabelled]	(2) at (0.85*0.9511, 0.85*0.3090)	{};
			\node[unlabelled]	(u) at (1.85*0.5878, 1.85*-0.8090)	{};
			\node	(l1) at (1.20*0, 1.20*1)	{$v_1$};
			\node	(l5) at (1.20*-0.9511, 1.20*0.3090)	{$v_5$};
			\node	(l4) at (1.20*-0.5878, 1.20*-0.8090)	{$v_4$};
			\node	(l3) at (0.50*0.5878, 0.50*-0.8090)	{$v_3$};
			\node	(l2) at (1.20*0.9511, 1.20*0.3090)	{$v_2$};
			\node	(lu) at (2.20*0.5878, 2.20*-0.8090)	{$u$};
			\draw	(1) -- (2)	(2) -- (3)	(3) -- (4)	(4) -- (5) (5) -- (1) (3) -- (u);
			\draw[dashed, ->-=0.5]	(2) to (u);
			\draw[dashed, ->-=0.5]	(4) to (u);
			\end{tikzpicture}
			\subcaption{(iv)}
		\end{subfigure}%
		\begin{subfigure}[b]{.33\textwidth}
			\centering
			\begin{tikzpicture}
			\node[unlabelled]	(1) at (0.85*0, 0.85*1)	{};
			\node[unlabelled]	(5) at (0.85*-0.9511, 0.85*0.3090)	{};
			\node[unlabelled]	(4) at (0.85*-0.5878, 0.85*-0.8090)	{};
			\node[unlabelled]	(3) at (0.85*0.5878, 0.85*-0.8090)	{};
			\node[unlabelled]	(2) at (0.85*0.9511, 0.85*0.3090)	{};
			\node[unlabelled]	(u) at (0.85*0.5878+1, 0.85*-0.8090)	{};
			\node[unlabelled]	(w) at (0.19, -1.64)	{};
			\node	(l1) at (1.20*0, 1.20*1)	{$v_1$};
			\node	(l5) at (1.20*-0.9511, 1.20*0.3090)	{$v_5$};
			\node	(l4) at (1.20*-0.5878, 1.20*-0.8090)	{$v_4$};
			\node	(l3) at (1.20*0.5878, 1.20*-0.8090)	{$v_3$};
			\node	(l2) at (1.20*0.9511, 1.20*0.3090)	{$v_2$};
			\node	(lu) at (0.85*0.5878+1.35, 0.85*-0.8090)	{$u$};
			\node	(lw) at (0.19, -1.99)	{$w$};
			\draw	(1) -- (2)	(2) -- (3)	(3) -- (4)	(4) -- (5) (5) -- (1) (3) -- (u) (3) -- (w);
			\draw[dashed, ->-=0.5]	(2) to (u);
			\draw[dashed, ->-=0.5]	(4) to (w);
			\end{tikzpicture}
			\subcaption{(v)}
		\end{subfigure}%
		\begin{subfigure}[b]{.33\textwidth}
			\centering
			\begin{tikzpicture}
			\node[unlabelled]	(1) at (0.85*0, 0.85*1)	{};
			\node[unlabelled]	(5) at (0.85*-0.9511, 0.85*0.3090)	{};
			\node[unlabelled]	(4) at (0.85*-0.5878, 0.85*-0.8090)	{};
			\node[unlabelled]	(3) at (0.85*0.5878, 0.85*-0.8090)	{};
			\node[unlabelled]	(2) at (0.85*0.9511, 0.85*0.3090)	{};
			
			\node[unlabelled]	(w) at (1.85*0, 1.85*1)	{};
			\node[unlabelled]	(u) at (1.85*0.5878, 1.85*-0.8090)	{};
			
			\node	(l1) at (0.50*0, 0.50*1)	{$v_1$};
			\node	(l5) at (1.20*-0.9511, 1.20*0.3090)	{$v_5$};
			\node	(l4) at (0.85*-0.5878+0.1, 0.85*-0.8090-0.35)	{$v_4$};
			\node	(l3) at (0.50*0.5878, 0.50*-0.8090)	{$v_3$};
			\node	(l2) at (1.20*0.9511, 1.20*0.3090)	{$v_2$};
			
			\node	(lw) at (2.20*0, 2.20*1)	{$w$};
			\node	(lu) at (2.20*0.5878, 2.20*-0.8090)	{$u$};
			\draw	(1) -- (2)	(2) -- (3)	(3) -- (4)	(4) -- (5) (5) -- (1) (1) -- (w) (3) -- (u);
			\draw[dashed, ->-=0.5]	(2) to (u);
			\draw[dashed, ->-=0.5]	(w) to (5);
			\end{tikzpicture}
			\subcaption{(vi)}
		\end{subfigure}
\caption{Obstructions $X$ for which $\overline{U(X)}$ contains an induced $C_5$.
\label{c5}}
\end{figure}
\end{theorem}
\pf Let $C: v_1v_2v_3v_4v_5$ be an induced $C_5$ in $\overline{U(X)}$. 
By Lemma \ref{c5caselemma}, $C$ is the unique cycle in $\overline{U(X)}$ and has 
a vertex which is an arc-balancing but not a cut-vertex of $\overline{U(X)}$. 
Without loss of generality assume $v_1$ is such a vertex. 
By Lemma \ref{c5caselemma}, $v_1$ balances an arc between a neighbour of $v_1$
in $C$ and a leaf. Without loss of generality, assume $v_1$ balances an arc 
between $v_2$ and $u$ and $(v_2,u)$ is the arc. Since any vertex 
except $v_1$ that is adjacent to $u$ is also adjacent to $v_2$, $u$ is adjacent to 
$v_3$. Suppose that neither $v_4$ nor $v_5$ is an arc-balancing vertex. Then
by Proposition \ref{vertex classification}, each of $v_4, v_5$ is a cut-vertex of 
$\overline{U(X)}$ or incident with an arc. Since $v_4, v_5$ are adjacent in 
$\overline{U(X)}$, there is no arc between them. So one of $v_4, v_5$ is not 
incident with an arc and hence must be a cut-vertex of $\overline{U(X)}$. 

Suppose $v_4$ is a cut-vertex and $v_5$ is incident with an arc. By Lemma 
\ref{c5caselemma}(b), there is a leaf $w$ adjacent to $v_4$ and incident with 
an arc. Hence there is an arc between $v_5$ and $w$. Since we have 
$(v_2,u)\Gamma(u,v_1)\Gamma(v_1,v_3)\Gamma(v_4,v_1)\Gamma(v_1,w)\Gamma(w,v_5)$ and 
the two arcs are opposing in $X$, the second arc is $(v_5,w)$. Thus 
$\overline{U(X)}$ is Figure~\ref{c5}(i). 
Suppose instead that $v_4$ is incident with an arc and $v_5$ is a cut-vertex. 
By Lemma \ref{c5caselemma}, there is a leaf $w$ adjacent to $v_5$ and incident with
an arc. Hence there is an arc between $v_4$ and $w$. Since 
$(v_2,u)\Gamma(u,v_1)\Gamma(v_1,v_3)\Gamma(v_3,v_5)\Gamma(w,v_3)\Gamma(v_4,w)$,
the second arc is $(w,v_4)$ and $\overline{U(X)}$ is Figure~\ref{c5}(ii). 
Finally, suppose both $v_4$ and $v_5$ are cut-vertices. 
By Lemma \ref{c3noc4lemma}, there are leaves $w, x$ adjacent to $v_4, v_5$, 
respectively, and incident with arcs. Hence there is an arc between $w$ and $x$. 
Since we have 
$(v_2,u)\Gamma(u,v_1)\Gamma(v_5,u)\Gamma(u,x)\Gamma(x,v_3)\Gamma(v_4,x)\Gamma(x,w)$,
the second arc is $(w,x)$. Thus $\overline{U(X)}$ is Figure~\ref{c5}(iii).

Suppose exactly one of $v_4,v_5$ is an arc-balancing vertex. Clearly they cannot
both be arc-balancing vertices because $v_1$ is such a vertex and there are at most 
two arc-balancing vertices. Consider first the case when $v_5$ is an arc-balancing 
vertex. Then $v_5$ either balances an arc between $v_1$ and a leaf adjacent to
$v_2$ or an arc between $v_4$ and a leaf adjacent to $v_3$. However, the former is 
not possible, as otherwise $v_4$ is not an arc-balancing vertex and not incident 
with an arc so it is a cut-vertex by Proposition \ref{vertex classification}. 
But then a leaf adjacent to it is not incident with an arc, a contradiction to 
Lemma \ref{c5caselemma}. So $v_5$ balances an arc between $v_4$ and a leaf $w$ 
adjacent to $v_3$. Then we have 
$(v_2,u)\Gamma(u,v_1)\Gamma(v_1,v_3)\Gamma(w,v_1)\Gamma(v_5,w)\Gamma(w,v_4)$.
The second arc is $(v_4,w)$. When $w = u$, $\overline{U(X)}$ is 
Figure~\ref{c5}(iv); otherwise $\overline{U(X)}$ is Figure~\ref{c5}(v).

Consider now the case when $v_4$ is an arc-balancing vertex.
Then $v_4$ either balances an arc between $v_5$ and a leaf adjacent to $v_1$ or
an arc between $v_3$ and a leaf adjacent to $v_2$. If $v_4$ balances an arc 
between $v_5$ and a leaf $w$ adjacent to $v_1$, then
$(v_2,u)\Gamma(u,v_1)\Gamma(v_1,v_3)\Gamma(v_3,w)\Gamma(w,v_4)\Gamma(v_5,w)$.
The second arc is $(w,v_5)$ and $\overline{U(X)}$ is Figure~\ref{c5}(vi).
On the other hand, if $v_4$ balances an arc between $v_3$ and a leaf adjacent to 
$v_2$, then $v_5$ is not an arc-balancing vertex and not incident with an arc so
it is a cut-vertex by Proposition \ref{vertex classification}. But then a leaf
adjacent to it is not incident with an arc, a contradiction to Lemma 
\ref{c5caselemma}. Hence this is not possible.
\qed

\section{Conclusion}

Theorem \ref{main} follows immediately from Corollary~\ref{noarc} and Theorems 
\ref{div-cut}, \ref{2 non-dividing cut-vertices}, \ref{non-dividing4-5}, 
\ref{non-dividing6}, \ref{non-dividing7}, \ref{non-dividing8}, 
\ref{disconnected theorem}, \ref{treecase}, \ref{c3only}, \ref{onec4}, 
\ref{twoc4s}, \ref{c3c4theorem}, \ref{c5theorem}, and the fact that an obstruction
must be contained in a partially oriented graph which cannot be completed to
a local tournament.

Partially oriented graphs which cannot be completed to local tournaments and are 
minimal with respect to only vertex deletions can be derived from 
Theorem~\ref{main}. Indeed, suppose that $Y$ is such a graph, that is, $Y$ cannot 
be completed to a local tournament and for each vertex $v$ of $Y$, $Y-v$ can be 
completed to a local tournament. Since $Y$ cannot be completed to a local 
tournament, it contains an obstruction.
Since $Y$ is minimal with respect to vertex deletion, $Y$ is either an obstruction 
described in Theorem~\ref{main} or is obtained from an obstruction by orienting 
some edges.
 
An oriented graph $D = (V,A)$ is {\em transitive} if, for any three vertices 
$u, v, w \in V$, $uv \in A$ and $vw \in A$ imply $uw \in A$; it is {\em 
quasi-transitive} if, for any three vertices $u, v, w \in V$, $uv \in A$ and 
$vw \in A$ imply $uw \in A$ or $wu \in A$, cf. \cite{bh}. 
The orientation completion problems for the classes of transitive and 
quasi-transitive oriented graphs are polynomial time solvable, cf. \cite{bhz}.
The underlying graphs of transitive oriented graphs coincide with 
the underlying graphs of quasi-transitive oriented graphs, which are known as 
{\em comparability graphs}, cf. \cite{gallai,gh,golumbic}. 
It remains open problems to find obstructions for transitive orientation 
completions and quasi-transitive orientation completions.

\end{document}